\documentclass[journal,web]{ieeecolor}

\nonstopmode

\usepackage{generic}
\usepackage{cite}

\usepackage{amsmath,amssymb,amsfonts}
  \usepackage{amsthm} 
\usepackage{graphicx}
\usepackage{textcomp}
\def\BibTeX{{\rm B\kern-.05em{\sc i\kern-.025em b}\kern-.08em
    T\kern-.1667em\lower.7ex\hbox{E}\kern-.125emX}}
\usepackage{yfonts}

\usepackage{url}
\usepackage{physics}
\usepackage{soul}

\usepackage[svgnames]{xcolor}

\selectcolormodel{cmyk}

\definecolor{testing}{RGB}{251,128,114}
\definecolor{sensing}{RGB}{190,186,218}
\definecolor{dataAnalytics}{RGB}{179,222,105}
\definecolor{wirelessComm}{RGB}{253,180,98}
\definecolor{control}{RGB}{255,255,179}

\definecolor{slateBrown}{RGB}{158, 134, 120 }
\definecolor{slateBrownDull}{RGB}{215, 208, 198}
\definecolor{dirtyWhite}{RGB}{227, 227, 227}
\definecolor{slateGreen}{RGB}{137, 183, 97}

\definecolor{sageGreen}{RGB}{185, 199, 141}
\definecolor{jasmineOrange}{RGB}{253, 220, 139 }
\definecolor{cyanTeal}{RGB}{107, 179, 168}
\definecolor{lightSaffron}{RGB}{255, 184, 83}

\definecolor{chinaPink}{RGB}{233, 105, 173}
\definecolor{blueMetallic}{RGB}{40, 89, 119}
\definecolor{dandelionYellow}{RGB}{246, 228, 43}
\definecolor{skyBlue}{RGB}{0, 177, 255}
\definecolor{hawaiiBlue}{RGB}{0, 109, 176}
\definecolor{charcoalBlue}{RGB}{44, 68, 94}
\definecolor{darkElectricBlue}{RGB}{77, 97, 121}

\definecolor{tolRedForContrast}{HTML}{A50026}
\definecolor{tolBlueForContrast}{HTML}{B7DDEB}
\definecolor{tolYellowMain}{HTML}{DDAA33}
\definecolor{tolRedMain}{HTML}{BB5566}
\definecolor{tolBlueMain}{HTML}{004488}
\definecolor{tolPaleYellow}{HTML}{EEEEBB}
\definecolor{tolPaleRed}{HTML}{FFCCCC}
\definecolor{tolDarkRed}{HTML}{663333}
\definecolor{tolDarkYellow}{HTML}{666633}
\definecolor{tolDarkBlue}{HTML}{222255}
\definecolor{tolPaleBlue}{HTML}{BBCCEE}
\definecolor{tolDarkCyan}{HTML}{225555}
\definecolor{tolPaleCyan}{HTML}{CCEEFF}

\usepackage{graphicx,subcaption}

  \usepackage{rotating}
  \usepackage{floatpag}
  \rotfloatpagestyle{empty}

\theoremstyle{plain}

\usepackage{mathrsfs,multicol}
\usepackage[mathscr]{euscript}

\usepackage{caption}
 
\usepackage{wrapfig}

\usepackage{tikz}
\usepackage{pgfplots}
\usepackage{tikz-3dplot}

\usepackage{tkz-euclide}
\usetikzlibrary{quotes,angles}

\usepgfplotslibrary{patchplots}
\usetikzlibrary{decorations.pathmorphing} 
\usetikzlibrary{matrix} 
\usetikzlibrary{arrows.meta} 
\usetikzlibrary{shapes} 
\usetikzlibrary{calc} 
\usetikzlibrary{positioning} 
\usetikzlibrary{patterns}
\usetikzlibrary{intersections}
\usepgfplotslibrary{fillbetween}
\usepgfplotslibrary{polar}
\usepgfplotslibrary{colorbrewer}
\usepgfplotslibrary{groupplots}
\usetikzlibrary{decorations.markings}

\usepgfplotslibrary{fillbetween}
\usepgfplotslibrary{groupplots}

\usetikzlibrary{hobby}

\pgfdeclarelayer{foreground}
\pgfsetlayers{main,foreground}

\usetikzlibrary{shapes.geometric,shapes.symbols}
\tikzstyle{block} = [draw,rectangle,thick,minimum height=2em,minimum width=2em]
\tikzstyle{sum} = [draw,circle,inner sep=0mm,minimum size=4mm]
\tikzstyle{connector} = [->,thick]
\tikzstyle{line} = [thick]
\tikzstyle{branch} = [circle,inner sep=0pt,minimum size=1pt,fill=black,draw=black]
\tikzstyle{guide} = []
\tikzstyle{legendBlock} = [rectangle,minimum height=2em,minimum width=2em]

\renewcommand{\vec}[1]{\ensuremath{\boldsymbol{#1}}} 

\definecolor{yellow1}{RGB}{255,255,204}
\definecolor{blue2}{RGB}{161,218,180}
\definecolor{blue3}{RGB}{65,182,196}
\definecolor{blue4}{RGB}{44,127,184}
\definecolor{blue5}{RGB}{37,52,148}
\definecolor{blueForRed1}{RGB}{5,113,176}
\definecolor{brown1}{RGB}{166,97,26}
\definecolor{brown2}{RGB}{223,194,125}

\definecolor{red1}{RGB}{215,25,28}
\definecolor{red2}{RGB}{253,174,97}

\definecolor{redSeries1}{RGB}{254,240,217}
\definecolor{redSeries2}{RGB}{253,204,138}
\definecolor{redSeries3}{RGB}{252,141,89}
\definecolor{redSeries4}{RGB}{215,48,31}

\definecolor{grayMax}{cmyk}{0,0,0,85}
\definecolor{grayOne}{cmyk}{0,0,0,1}
\definecolor{grayTwo}{cmyk}{0,0,0,20}
\definecolor{grayTwoAndHalf}{cmyk}{0,0,0,28}
\definecolor{grayThree}{cmyk}{0,0,0,41}
\definecolor{grayFour}{cmyk}{0,0,0,61}
\definecolor{magentaForBlack}{RGB}{202,0,32}
\definecolor{peachForBlack}{RGB}{244,165,130}

\definecolor{darkViolet}{cmyk}{65,70,0,0}
\definecolor{darkSea}{cmyk}{85,30,0,0}
\definecolor{paleSea}{cmyk}{33,3,0,0}
\definecolor{paleViolet}{cmyk}{65,70,0,0}
\definecolor{paleGreen}{cmyk}{24,0,39,0}
\definecolor{paleOrange}{cmyk}{5,35,70,0}

\definecolor{c1}{RGB}{179,88,6}
\definecolor{c2}{RGB}{241,163,64}
\definecolor{c3}{RGB}{216,218,235}
\definecolor{c4}{RGB}{153,142,195}
\definecolor{c5}{RGB}{84,39,136}

\definecolor{tacGreen1}{cmyk}{50,0,17,0}
\definecolor{tacGreen2}{cmyk}{100,10,55,0}
\definecolor{tacBrown1}{cmyk}{12,20,45,0}
\definecolor{tacBrown2}{cmyk}{35,55,90,0}

\definecolor{darkTeal}{cmyk}{100,10,55,0}
\definecolor{lightTeal}{cmyk}{50,0,17,0}
\definecolor{toffee}{cmyk}{12,20,45,0}
\definecolor{darkToffee}{RGB}{166,97,26}

\definecolor{purDark}{RGB}{123,50,148}
\definecolor{purLight}{RGB}{194,165,207}
\definecolor{paleGreen}{RGB}{166,219,160}
\definecolor{darkGreen}{RGB}{0,136,55}

\definecolor{myCrimson}{RGB}{202,0,32}

   \newtheorem{theorem}{Theorem}
    \newtheorem{corollary}{Corollary}[theorem]
    \newtheorem{lemma}{Lemma}

    \newtheorem{definition}{Definition}
    
    \newtheorem{example}{Example}

 \pgfplotsset{compat=newest}

 \usepackage{mathtools}
 \usepackage{amsfonts}
 \usepackage{ucs}
 \usepackage[utf8x]{inputenc}
 \usepackage[T1]{fontenc}
 \usepackage{ae,aecompl}

 \usepackage[many]{tcolorbox}

\usepackage{xtab}
\usepackage{array, booktabs, caption, ragged2e}
\DeclareCaptionFont{blue}{\color{hawaiiBlue}}

\makeatletter
\let\mcnewpage\newpage
\newcommand{\changenewpage}{%
  \renewcommand\newpage{%
    \if@firstcolumn
      \hrule width\linewidth height0pt
      \columnbreak
    \else
      \mcnewpage
    \fi
}}
\makeatother

\usepackage{colortbl}

\makeatletter
\let\oldlt\longtable
\let\endoldlt\endlongtable
\def\longtable{\@ifnextchar[\longtable@i \longtable@ii}
\def\longtable@i[#1]{\begin{figure}[t]
\onecolumn
\begin{minipage}{0.5\textwidth}
\oldlt[#1]
}
\def\longtable@ii{\begin{figure}[t]
\onecolumn
\begin{minipage}{0.5\textwidth}
\oldlt
}
\def\endlongtable{\endoldlt
\end{minipage}
\twocolumn
\end{figure}}
\makeatother

\usepackage{fontawesome5}
\usepackage{physics}

\renewcommand{\theequation}{\arabic{equation}}

\begin{document}

\title{Predicting oscillations in relay feedback systems, using fixed points of Poincar{\'e} maps, and Hopf bifurcations}
\author{Maben Rabi
\thanks{submitted: ~2023}%
\thanks{Faculty of  Computer Sciences, Engineering and Economics, {\O}stfold University College, Halden, Norway (e-mail: firstname.lastname@hiof.no).}%
}
\maketitle

\begin{abstract}
    The relay autotuning method identifies plant parameters, from oscillations of the plant under relay feedback.  
    To predict the presence and nature of such oscillations, we apply the following two approaches: (a)~analysis of the switching dynamics, while using an ideal relay, and (b)~bifurcation analysis, while using a smooth 
    approximation of the relay.
For  stable plants with positive DC gains, 
 our analyses predict that: 
     (i)~a periodic orbit is guaranteed, for a class of non-minimum phase plants of relative degree one, whose step response starts with an inverse response, and (ii)~for a wider class of plants, whose root locus diagrams cross the imaginary axis at complex conjugate values, limit cycles are merely suggested. 
%
\end{abstract}

\begin{IEEEkeywords}
Relay feedback, limit cycle, self-oscillations, Relay auto-tuning, Hopf bifurcation, fixed point theorem. 
\end{IEEEkeywords}

\section{Introduction}
 \newcommand{\myrowcolourYellow}{\rowcolor[tolPaleYellow]{0.1}}
  \newcommand{\idealrelay}{
    \begin{tikzpicture}[yscale=0.25,xscale=0.20]
        \draw[ultra thick] (-2,-1)--(0,-1)--(0,1)--(2,1);
        \draw[thin, <->, black!50] (-2,0)--(2,0);
        \draw[thin, <->, black!50] (0,-2.5)--(0,2.5);
    \end{tikzpicture}
       }
\newcommand{\relayWithHysteresis}{
    \begin{tikzpicture}[yscale=0.25,xscale=0.34]
        \draw[ultra thick] (-2,-1)--(0.8,-1);
        \draw[ultra thick] (-0.8,1)--(2,1);
        \draw[thin, <->] (-2,0)--(2,0);
        \draw[thin, <->] (0,-2.5)--(0,2.5);
\begin{scope}[very thick,decoration={
               markings,
               mark=at position 0.5 with {\arrow{>}}}
             ] 
     \draw[postaction={decorate}] (0.8,-1)--(0.8,1);
     \draw[postaction={decorate}] (-0.8,1)--(-0.8,-1);
\end{scope}
    \end{tikzpicture}
       }
\IEEEPARstart{R}{elay} 
 auto-tuning~\cite{astromHagglund1984relayAutotuning}
 uses oscillations induced by feedback,
to identify a dynamical model for a given plant.
The auto-tuner runs a closed-loop experiment, by placing the plant  in
the Relay feedback system~(RFS) shown in Figure~\ref{fig:blockDiagram}. 
If  periodic oscillations occur, then the input to the plant
is a periodic rectangular wave. Then we can 
 use the period, amplitude and spectrum of the plant's output,
to learn aspects of the plant's dynamics.
For example, suppose that the RFS oscillations are symmetric, and are unimodal meaning that there are exactly two relay switches per period. Then the relay oscillation frequency is a good approximation for some 180-degree phase cross-over frequency of the plant's transfer function.
%
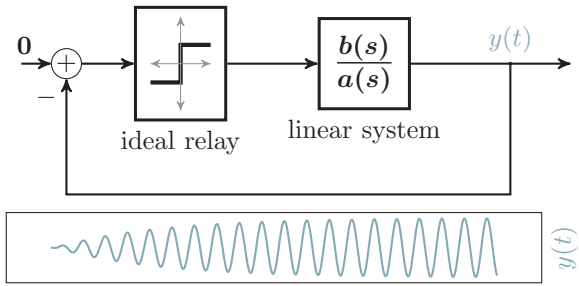
\begin{figure}
\begin{center}
\begin{tikzpicture}[scale=1, auto, >=stealth']
 \matrix[ampersand replacement=\&, row sep=0.3cm, column sep=0.4cm]
{
    \node(startOne) {};    \&
      \node[sum] (summer) {$\mathbf{+}$} ;  \&
    \node[block,minimum size=1.2cm,line width=0.4mm,label=below:{ideal relay}] (relay)  {\idealrelay }; \&
    \node[block,minimum size=1.2cm,line width=0.4mm,label=below:{linear system}] (plant)  {  {\mbox{\Large$\vec{{\frac{b(s)}{a(s)}}}$}}  }; \&
      \node[branch, label = above:{\textcolor{cyanTeal}{$ y(t) $}}] (bendOne) {} ; \& 
      \node (output) {};  
      \\
      \&
      \node (bendThree) {} ; \&
      \& 
     \&
      \node (bendTwo) {} ;  \&
  \\
    };
    \draw [connector] (startOne) -- node[pos=0.1] {$\vec{0}$} (summer);
    \draw [connector] (summer) -- (relay);
    \draw [connector] (relay)--(plant);     \draw [connector] (plant)--(output);
    \draw [connector] (bendOne.center)  -| (bendTwo.center) -- (bendThree.center) -|  node[pos=0.92]
  {$\mathbf{-}$} (summer);
\end{tikzpicture}
  \vspace{-7pt}
  \begin{tikzpicture}
\begin{axis}[
              yscale = 1.3,
              ylabel = {\textcolor{cyanTeal}{$ y(t) $}},
              yticklabel pos = right,
              title style={yshift=-1.6ex,},
              width=8.7cm,
              height=2.3cm,
              ticks = none ,
           ]
\addplot[thick, cyanTeal  ,title = {\emph{ $ {\frac{ - s + 1}{s^3 + 3 s^2 + 5 s + 6}} $}},
grid=both, every major grid/.style={gray, opacity=0.5}, 
    ]
    file {figuresAndData/relayOscillationFourthOrderNum-1+1+1.csv};
\end{axis}
\end{tikzpicture}
\end{center}
    \caption{Spontaneous oscillations under relay feedback.\label{fig:blockDiagram}} 
\end{figure}

     A rich variety of dynamical behaviours is  possible   for relay feedback systems. 
     These include settling into quiescence, asymmetric oscillations~\cite{boiko2021asymmetricOscillationsRelayFeedback}, chattering~\cite{flugge-lotz1953discontinuousControlSystems}, and even chaotic fluctuations~\cite{anosov1959relayChaos,cook1985relayChaos,holmberg1991thesisRelay,johanssonBarabanovAstrom2002chattering,sieber2010relayWithDelay,jeffrey2018discontinuous}. 
We do not know
yet 
of broad classes of plants for which oscillations, and their global stability can be guaranteed.





\subsection{Previous results on relay oscillations}
First order plants with positive DC gain, a stable pole, and a time delay have been shown by {\AA}str{\"o}m~\cite{astrom1995oscillationsRelay}, to self-oscillate. For such a plant, exactly one limit cycle is stable. This is also symmetric, and unimodal. There is an infinity of other periodic orbits, which are all unstable.

Second order plants with a positive DC gain, stable poles, and with one positive zero  have been shown to self-oscillate~\cite{holmberg1991thesisRelay,chernysheva2014hopfBifurcationForRelay,mabenRabi2021secondOrderRelayOscillations,rabi2018relay}. The resulting limit cycle is symmetric, unimodal, and globally stable. 




Time domain tools such as the Hamel locus~\cite{gillePelegrinDecaulne1959feedbackControlSystems,leMaitre1970graphicalAnalysisRelayOscillations} and frequency domain tools such as the Tsypkin locus~\cite{tsypkin1984relayControlSystems} give approximate necessary conditions for limit cycles. But these tools give no direct information about stability. The accuracies of these tools have been discussed, and improved by refinements such as the A-locus variant of the Tsypkin locus, and higher order Describing functions~\cite{bohn1961relay,bergen1962tsypkin,weischedel1973tsypkin,juddChirlian1974graphicalAnalysisRelayLimitCycle,judd1977errorBoundsRelaySystemAnalysis,atherton1981relayOscillations,mees1981dynamicsOfFeedbackSystems,kuznetsov2021tsypkinMethod}. Judd~\cite{judd1975relationshipHamelTsypkin} concludes that the Tsypkin locus gives more accurate predictions for plant transfer functions with fast roll-off at high frequencies, and that the Hamel locus gives more accurate predictions 
with slow roll-off.


The classic study of 
Andronov, Vitt and Khaikin~\cite{andronovVittKhaikin1966theoryOfOscillations,minorsky1962nonlinearOscillations} gave a method for analyzing individual second order systems, based on properties of ordinary differential equations on the plane. 
They proved the stability of some second order valve oscillators, which were significant at that time for Radio engineering.

A state space necessary condition for a symmetrical, periodic orbit is given in {\AA}str{\"o}m~\cite{astrom1995oscillationsRelay}. This paper also gives a linearized analysis for local stability in the neighbourhood of a limit cycle.
Goncalves~et~al.~\cite{goncalvesMegretskiDahleh2001relay} show how a Lyapunov analysis for global stability can be carried out.
A bifurcation analysis of the RFS, with a plant parameter being the bifurcation parameter is shown in~\cite{marioDiBernardoKarlJohansson2001ijbc,chernysheva2014hopfBifurcationForRelay}.

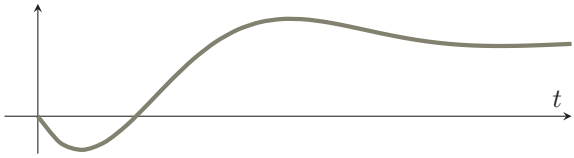
\begin{figure}
    \begin{center}
        \begin{tikzpicture} 
            \begin{axis}[axis lines = middle, ticks = none, xmin = -0.5, xmax = 8, ymin = -0.25, ymax = 0.75, xlabel=$t$,  
                x post scale = 1.1,  y post scale = 0.35, 
                ]					]
                \addplot[domain=0:8, smooth, tolDarkYellow, ultra thick] {0.5 - exp(-0.5*x)*( 0.5*cos(x*180*7/22) + sin(x*180*7/22) )};
                %
                %
            \end{axis}
        \end{tikzpicture}
    \end{center}
    \caption{\label{fig:megretskiTemplateForStepResponse}%
    Megretski's~\cite{megretski1996globalStabilityOfRelayOscillations} template for the plant's open loop step response, to make the closed loop RFS self-oscillate}
\end{figure}

Megretski~\cite{megretski1996globalStabilityOfRelayOscillations} gives a graphical template for the step response of a plant, for it to self-oscillate under relay feedback. This template is shown in Figure~\ref{fig:megretskiTemplateForStepResponse}. It shows a positive DC gain, is stable, and has an initial undershoot.
Megretski proves that if a given plant step response resembles this template, then self-oscillations are guaranteed, provided that a unique stationary solution exists, for a discrete time iteration involving 
inter-switching times.

The map from the state at one switching time to the state after a specified number of switches, has been called the Poincaré map. It has been suggested that RFS limit cycles be studied, by studying fixed points of such maps~\cite{blimanKrasnoselskii1997periodicSolutions,varigondaGeorgiou2001relay}, for example using Banach's contraction mapping theorem. 
But only for stable, 
second order plants has this map been shown to be contracting. 
Johansson and Rantzer~\cite{johanssonRantzer1996globalAnalysisOfThirdOrderRelay} treat the case of a third order plant with no zeroes, and three distinct, real, stable poles. They find conditions under which the said map obeys an `area contraction' property. 
When the plant is of higher order, we did not know properties required of
plants,
under which such maps have fixed points, and whether or not these points are globally or at  the least, locally attracting. 

\begin{figure}
\begin{center}
    \includegraphics[width=\linewidth]{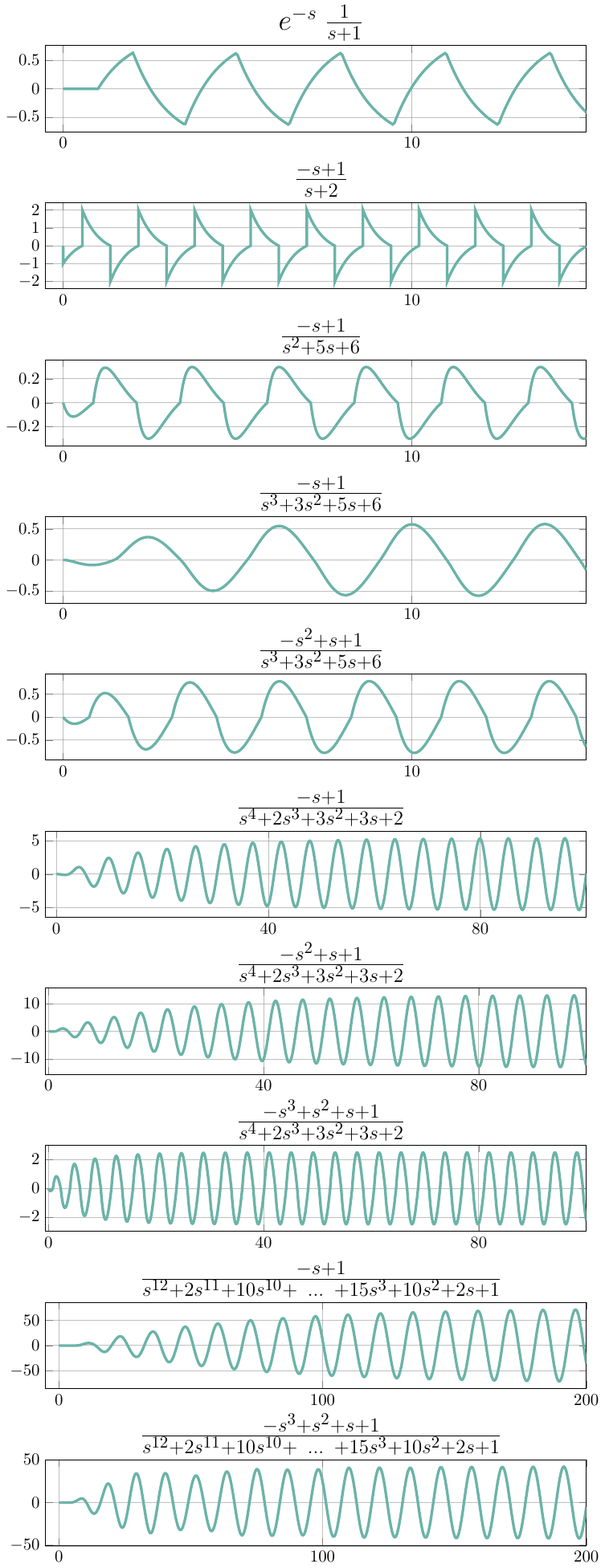}
\end{center}
\caption{RFS output for some nonminimum phase plants.\label{fig:oscillationPlots}}
\end{figure}

\subsection{Our results}
We study the switching dynamics of the RFS, by appealing to the topological and analytical (holomorphic) properties of trajectories of linear ODEs.
In specific, we show in Section~\ref{section:poincareMapIsAnalytic} that
%
for stable plants, the Poincaré map is not merely continuous, but analytic.

With this, we guarantee a regular periodic orbit, 
if the plant is stable, is rational, has relative degree of one, has a positive DC gain, and has an odd number of positive real zeros~(Section~\ref{section:existence}). Thus at least when the relative degree of the plant equals one, 
the RFS has a periodic orbit, 
if the plant step response matches Megretski's template.

In Section~\ref{section:globalStability} we provide a conditional guarantee of global stability. We consider the strong assumption that a Poincaré map is locally Schur stable, everywhere on a bounded, positively invariant set. Under this, we prove that the promised periodic orbit is in fact a symmetric, unimodal and
globally asymptotically stable limit cycle.



We change tools
in Section~\ref{section:smoothApproximation}. We approximate the relay by a cascade combination of a linear amplifier with high gain, and the hyperbolic tangent function. We take the resulting family of smooth nonlinear ODEs,
and look for bifurcations and structural stability.
 
 This suggests but does not confirm the presence or absence of limit cycles.
    We apply  the Hopf bifurcation theorem for systems in Lure feedback form~\cite{allwright1977harmonicBalanceAndHopfBifurcation,mees1981dynamicsOfFeedbackSystems,stanSepulchre2007analysisOfInterconnectedOscillators,muratArcakSalama2020ringOScillatorsTanh}, but also a stronger version due to Alexander and Yorke~\cite{alexanderYorke1978globalBifurcationsOfPeriodicOrbits}.
These suggest that the  number of limit cycles depends on the number of crossings of the imaginary axis, in the root locus diagram of the plant. 

Sections~\ref{section:liouville} and~\ref{section:monodromyMatrix} present the linearization around periodic orbits for the smooth approximation of the RFS. Here our calculations of local stability refine those in~\cite{balasubramanian1981stabilityOfLimitCycle,majhiAtherton1998stabilityOfLimitCyclesInRelay,johanssonRantzerAstrom1999fastSwitchesInRFS}.
We end up with an exact  expression for the monodromy matrix, 
 which improves that given in~\cite{majhiAtherton1998stabilityOfLimitCyclesInRelay}, when the plant's relative degree equals one.


\section{Preliminaries\label{section:introducingBRLURFplants}}
We model the RFS dynamics, first as an ODE with a  discontinuous right hand side, and then in the next section, as a nonlinear iteration in discrete-time. 

\subsection{Relay characteristic}
\noindent
We model the relay operator as an ideal binary switch, that has no bias, no delay, and no hysteresis. This is described using the hard signum function: 
\begin{align*}
    {\text{sign}} \bigl( e  \bigr)
    & = 
    \begin{cases}
        +1 , & {\text{if}} \ e > 0 , \\
        -1 , & {\text{if}} \ e < 0 , \\
        \left[ - 1 , + 1 \right] , & {\text{if}} \ e = 0 . 
    \end{cases}
\end{align*}
\subsection{The BRL-URF plants}
\noindent
Let the plant transfer function be: 
\begin{align}
    {\frac{b \left( s \right)}{a \left( s \right)}} & = {\frac{\phantom{s^{n} + \ \ } b_{n-1}s^{n-1} + \ldots + b_0 }
{s^n + a_{n-1}s^{n-1} + \ldots + a_0} }.
    \label{eqn:BRLURFtransferFunction}
\end{align}
\begin{definition}
    A transfer function is
    {\textbf{B}}ounded and
    {\textbf{R}}est{\textbf{L}}ess
    {\textbf{U}}nder
    {\textbf{R}}elay
    {\textbf{F}}eedback~(BRL-URF)
if it is rational, strictly proper, has a pole excess of one, is stable, has positive DC gain, and has an odd number of positive real zeroes.
\end{definition}
 There are constraints on the signs of some of the coefficients of the  transfer function~\eqref{eqn:BRLURFtransferFunction}, if it is BRL-URF.
 Firstly, Hurwitz stability requires all the denominator coefficients~$ a_i  $ to be positive, like the denominator's leading coefficient is. Then the DC gain being positive forces the numerator's constant term~$ b_0 $
  to be positive as well.  And finally, having an odd number of positive real zeros requires the numerator's leading coefficient to be 
  negative~\cite{vidyasagar1986undershoot,vidyasagar1987authorsReply,janMaciejowski2018rhpZeroesNotNecessaryForInverseResponse}.

The name BRL-URF is suitable becuase, for such plants, 
the corresponding RFS has: 
\begin{itemize}
    \item{only {\textit{bounded}}  trajectories - this is guaranteed by the Hurwitz stability of the plant~(Section~\ref{section:ultimateMagnitudeOfstate}), and,}
    \item{only {\textit{restless}} trajectories, 
        meaning that: 
         \begin{itemize}
            \item{the RFS has no equilibrium point, because the DC gain is positive, and,}
            \item{the RFS has no chattering point, because: (i)~the DC gain is positive, and (ii)~the coefficient~$b_{n-1}$ is negative.
                See Figure~\ref{fig:fieldsAtSwitchingLine}, and~\cite{johanssonRantzerAstrom1999fastSwitchesInRFS}.}
         \end{itemize}
        }
\end{itemize}
Strictly speaking, a sliding set exists, where trajectories cannot be considered to  be restless. But at every point on it, non-sliding Filippov solutions also exist. On those points where this non-uniqueness is present, we choose to work with the non-sliding solutions, because with the addition of even the smallest amount of noise, trajectories shall get away from the sliding set, and never return to it. 

The step response of a BRL-URF plant matches the template of Megretski - it starts with an inverse response, is bounded, and has a positive steady state.

\subsection{The state space realization}
We use 
the Observer canonical realization:
\begin{gather}
    A = \begin{bmatrix}
        0 & 0 & \cdots  & \cdots  & 0 & -a_0 \\
        1 & 0 &\cdots   &  \cdots & 0 & -a_1 \\
        0 & 1 & 0   & \cdots      & 0 & -a_2 \\
        \vdots & \vdots & \ddots & \ddots & \vdots  & \vdots \\
        0 & 0 & \cdots    & 1 & 0  & -a_{n-2} \\
        0 & 0 & \cdots    & 0  & 1  & -a_{n-1} 
        \end{bmatrix},
    \quad
    B = \begin{bmatrix}
        b_0 \\ b_1 \\ b_2 \\ \vdots \\ b_{n-2} \\ b_{n-1}
        \end{bmatrix}, 
        \label{eqn:companionA}
        \\
    C = \begin{bmatrix}
        \ 0 \quad  0 \quad   0 \quad   \cdots \quad   0 \quad  1 \ \
        \end{bmatrix}.
        \nonumber
\end{gather}
Then the relay feedback system evolves as per:
\begin{align}
    {\frac{d}{dt}} x & =
    A x  - B \, \, {\text{sign}} \Bigl( C \, x \Bigr) .
    \label{eqn:RFSdynamics}
\end{align}

\subsection{Definitions of exit times and exit maps\label{section:notationAndDefinitions}}
The switching hyperplane~${\mathcal{S}}$ is the set where the relay switches sign. For our chosen realization, it is that hyperplane whose points have zero as the~$x_n$ coordinate: 
\begin{align*}
{} {\mathcal{S}} & = \left\{ x \in  {\mathbb{R}}^n : C x = 0 \right\} .
\end{align*}
\begin{definition}
    The first exit time from positive sign 
    is the nonnegative function:
\begin{gather*}
    {\tau_{+}} \ : \ {\mathbb{R}}^n   \to
    {\mathbb{R}} \ \text{such that} \ \
{\tau_{+}} \left( \xi \right) \ = \  
    \inf \left\{ t > 0 :  C x\left(t\right) < 0 \right\} , \\
    \text{where} \
    {\dot{x}} \ = \ A x - B, \ \text{and} \ x\left( 0 \right) = \xi.
\end{gather*}
\end{definition}
\begin{definition}
    The first exit map from positive sign 
    is defined at every point for which the trajectory starting there  crosses the switching plane in finite time:
\begin{gather*}
    {\psi_{+}} \left( \cdot ; 1 \right) : \ {\mathbb{R}}^n   \to
    {\mathbb{R}}^n \ \text{such that} \ \
    {\psi_{+}} \left( \xi ; 1 \right) \ = \  x\left(\tau_{+}\right), 
    \\ 
    \text{where} \
    {\dot{x}} \ = \ A x - B, \ \text{and} \ x\left( 0 \right) = \xi.
\end{gather*}
\end{definition}
Exactly in analogy to $ {\tau_{+}} \left( \cdot \right), \; {\psi_{+}} \left( \cdot ; 1 \right) $, we can define the first exit time function~$ {\tau_{-}} \left( \cdot \right), $ and first exit map~${\psi_{-}} \left( \cdot  ; 1 \right) $ from negative sign, where the ODE in operation shall be:~$ {\dot{x}} = Ax + B $. Then
the following identities hold: 
\begin{gather*}
{\tau_{-}} \left( x \right) \; = \; {\tau_{+}} \left(  - x \right), \quad \text{and,} \quad
{\psi_{-}} \left( x ; 1 \right) \; = \;
 - {\psi_{+}} \left( - x  ; 1 \right).
\end{gather*}
\begin{definition}
 For positive integers~$k,$   the $k$-th exit map from positive sign
for the RFS~\eqref{eqn:RFSdynamics} is defined recursively as follows.
    For~$k \ge 2,$
\begin{align}
    \psi_{+} \left( x ; k \right) 
    &  = 
\psi_{+} \Bigl( 
\
    {\mathbf{-}} \, \psi_{+} \left( x ; k - 1 \right) 
\
    ; 1 \Bigr) .
    \label{eqn:firstExitMap_k}
\end{align}
\end{definition}
At any point~$x,$ regardless of whether or not a closed orbit passes through it, we shall call any~$k$th exit map~$ \psi_+\left( \cdot ; k \right) ,$ as a Poincar{\'e} map. 

The discussions so far lead to the following lemmas. 
\begin{lemma} 
    The first exit time~$\tau_+\left(  \,\cdot\, \right)$ is finite everywhere, when the plant is stable, and has a positive DC gain.
\end{lemma} 
\begin{lemma} 
    The first exit map~$\psi_+\left(  \,\cdot\, ; 1 \right)$ exists everywhere,
    when the plant is stable and has a positive DC gain.
\end{lemma} 

In this paper, the adjective {\textit{analytic}} means  {\textit{complex analytic,}} which is the same as {\textit{holomorphic.}} 
A function from~$ {\mathbb{C}}^n $  to~$ {\mathbb{C}} $ is called analytic, if it is analytic separately in each each element of its argument vector.
We say that a map 
from the real Euclidean space~$ {\mathbb{R}}^n $ 
to itself 
is analytic at a point~$x,$ if there is an analytic extension to this map at~$x.$ In other words,  when we embed the point~$x$ in the complex vector space~$  {\mathbb{C}}^n ,$ there should be an open ball around~$x$ in the complex domain, where the extended map is complex analytic~(it has complex-derivatives of all orders at~$x,$ and has a power series expansion that converges everywhere within the open ball). 

Where do analytic maps appear in our analysis ? 
The RFS trajectories that we study are not analytic, but are continuous  and  piecewise analytic functions of time. In specific, the trajectory coordinates are analytic functions of time, between any two successive relay switches. 
 Next we establish that  the first exit times and maps are 
  analytic.


\section{The Poincaré map is analytic\label{section:continuityOfPoincare}\label{section:poincareMapIsAnalytic}}
\begin{lemma}
    \label{lemma:usingImplicitFunctionTheorem}
    Suppose that the plant poles are all non-zero, so that the matrix~$A$ is invertible. If at a point~$\xi\in{\mathbb{R}}^n$, 
    \begin{gather}
        C A^{-1} B +  C  e^{A\tau_+\left(\xi\right)} \left( \xi -  A^{-1}B  \right) \ = \ 0, \label{eqn:crossing} \\
        C A e^{A\tau_+\left(\xi\right)} \left( \xi -  A^{-1}B  \right) \ \neq \ 0, \label{eqn:transversality}
    \end{gather}
    then the function~$\tau_+\left(\,\cdot\,\right)$ and the map~$\psi_+\left(\,\cdot\, ; 1 \right)$ are analytic in a small neighbourhood of~$\xi.$  
\end{lemma}
\begin{proof}
Consider the following analytic function of~$ \xi , t $:
\begin{align*}
    {\rm{output}}  \left( \xi, t \right) & \triangleq CA^{-1}B + C e^{At} \left( \xi -  A^{-1}B  \right)
\end{align*}
    Then $\tau_+\left(\xi\right)$  is implicitly defined through the  transcendental equation: $  {\rm{output}}   \left( \xi, \tau_+ \left( \xi \right) \right) = 0. $ The partial derivative 
\begin{align*}
    {\left.   {\frac{\partial}{\partial t}} {\rm{output}}   \left( \,\cdot\, , \,\cdot\, \right) \right\rvert}_{\xi, \tau_+\left(\xi\right) }  & 
      =  C A e^{A\tau_+\left(\xi\right)}  \left( \xi -  A^{-1}B  \right)  . 
\end{align*}
    The Implicit function theorem~\cite{krantz2002implicitFunctionTheorem} states that,
    if the above partial derivative is non-zero, then there is an open
 neighbourhood of~$\xi,$ 
    on which  the implicit function~$ \tau_+\left( \,\cdot \,\right) $ exists, and inherits the continuity, differentiability  and analyticity properties of the function:~$ {\rm{output}} \left( \cdot , \cdot \right) $.
Since the map~$  \psi_+\left( \xi ; 1 \right)   $ is an analytic function
of~$ \tau_+, $ it too inherits those properties. 
\end{proof}
\begin{lemma}
    \label{lemma:crossingHypersurface}
    For a given initial condition let the trajectory of a constant coefficient affine ODE cross an analytic hypersurface, after a finite duration. If the crossing is transversal, then locally around the initial condition, there is an open ball where the crossing time, and the crossing map are both defined, and are analytic.
\end{lemma}
\begin{proof}
   Crossing the hypersurface transversally at any given point, is nothing other than transversally crossing the tangent hyperplane at  that point. The implicit function can be then be used, exactly like for  Lemma~\ref{lemma:usingImplicitFunctionTheorem}.
\end{proof}


An {\textit{inflection point}} for a smooth space curve is a point at which the curve crosses its tangent line. 
Alternatively, it is a point
where there is at least one hyperplane that the curve crosses tangentially.

\begin{figure}
\centering
    \subcaptionbox{Conjugate poles%
\label{fig:secondOrderOscillatory}}
{\includegraphics[width=0.42\linewidth]{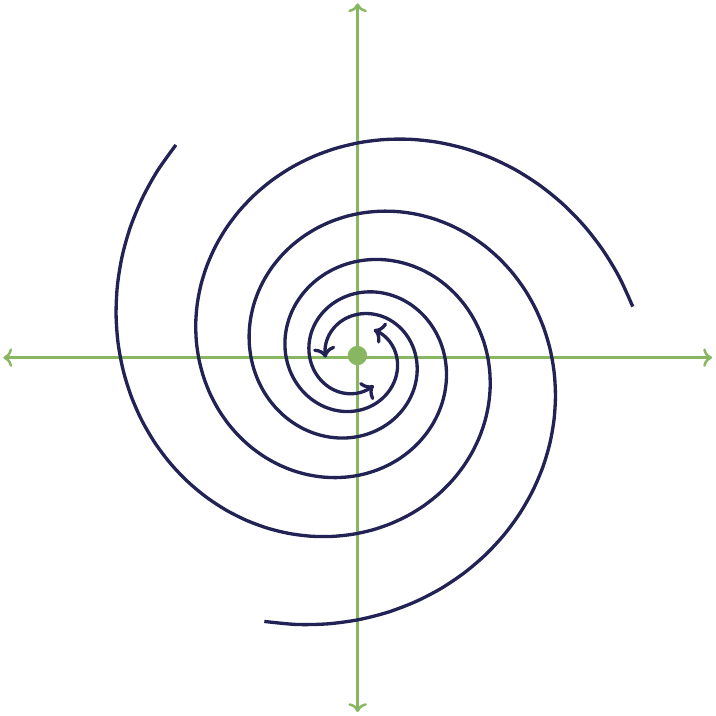}}%
\hfill 
    \subcaptionbox{Distinct real negative poles%
\label{fig:secondOrderDistinctReal}}
{\includegraphics[width=0.56\linewidth]{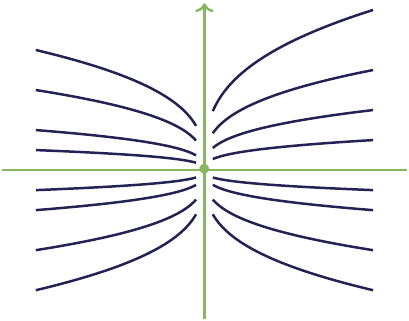}}%
    \caption{Trajectories having no inflection points.
    \label{fig:secondorderNoInflection}}
\end{figure}
If a second order transfer function has only non zero poles, then there is no inflection point in any trajectory of the corresponding affine ODE.
The curvature at points on any trajectory has the same sign throughout the trajectory, for every trajectory of such an ODE. 
For example, Figure~\ref{fig:secondorderNoInflection}
 illustrates the shapes of trajectories for second order plants having stable poles. The following theorem is the equivalent statement for any finite dimension.
\begin{theorem}
    Consider the affine system:
    \begin{gather}
        {\dot{x}} \ = \ A x - B , \quad
        y \ = \ C x ,
        \label{eqn:affineODE}
    \end{gather}
    where the matrix~$ A $ and the vectors are given by~\eqref{eqn:companionA}. 
    If~$ A $  is stable, then every trajectory that crosses the hyperplane
\begin{align*}
    {\mathcal{H}}  
    & 
    \triangleq
    \left\{  \
          x \in {\mathbb{R}}^n  \ : \ 
          C x =  0   
    \right\} ,
\end{align*}
must cross it transversally.
\end{theorem}
\begin{proof}
\begin{figure}
\centering
    \subcaptionbox{The hypothetical trajectory that tangentially crosses the hyperplane~${\mathcal{H}} ; $  note that to the left of the hypothetical  point of tangential crossing, the vector field crosses~${\mathcal{H}}$ downwards, while to the right, the vector field crosses~${\mathcal{H}}$ upwards.
\label{fig:crossingTangentially}}
    {\includegraphics[width=0.96\linewidth,trim={1.3mm 2.2mm 1.5mm 8mm},clip]{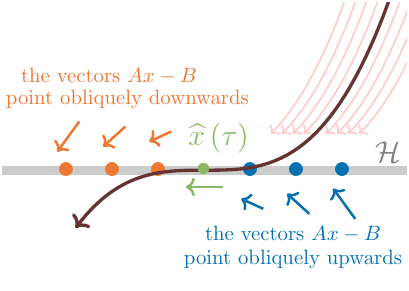}}%
    \\
    \subcaptionbox{The only possible tangential crossing of a smooth hypersurface by a trajectory, for the ODE:~${\tfrac{d}{dt}}{\widehat{x}} = {\mathbf{-}} {\widehat{x}} .$ The vector field crosses downwards the hypersurface~$ 
    {\widetilde{\mathcal{H}}} , $ that is $
    h \left( {\mathcal{H}} \right) $, on either side of the point of tangential crossing.%
\label{fig:topologicalEquivalence}}
    {\includegraphics[width=0.96\linewidth,trim={6mm 5mm 4mm 8mm},clip]{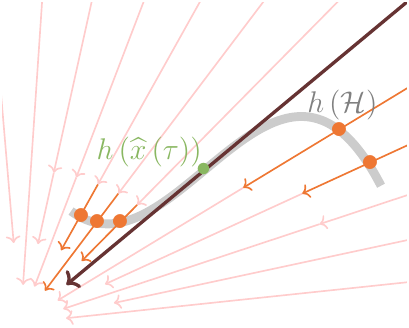}}%
    \caption{Contradiction between two descriptions of the vector field, around the hypothetical tangential crossing.\label{fig:crossingTheSwitchingHyperplane}}
\end{figure}
We shall prove this by contradiction.
Assume for the sake of argument, that a trajectory starting at some~$ x_0  $ crosses the hyperplane tangentially, after a duration of~$ \tau . $  If the crossing is indeed tangential, then the crossing point ought to be 
   within at the following set  intersection:
    \begin{align*}
     {\mathcal{I}}     &
    \triangleq
    \left\{  \
          x \in {\mathbb{R}}^n  \ : \ 
          C x =  0   
    \right\} 
        \cap
    \left\{  \
          x \in {\mathbb{R}}^n  \ : \ 
          C A x  - C B =  0   
    \right\} , 
        \\
      & = \left\{  \
          x \in {\mathbb{R}}^n  \ : \ 
         x_n = 0 , \ \text{and} \  x_{n-1} - a_{n-1} x_{n} -  b_{n-1}  =  0   
    \right\} .
    \end{align*}
    Using a topological argument, we shall prove that a crossing via the above set~$ {\mathcal{I}}   $ is impossible.

    Consider the phase portrait for the observer realization~\eqref{eqn:companionA}, which we have used so far. The two dimensional slice 
    spanned by the coordinates~$ x_{n-1}  $ and $ x_{n}  $ is shown in Figure~\ref{fig:crossingTangentially}. 
    To the left of the hypothetical crossing point, the vector field crosses downwards, in the direction of decreasing~$ x_n . $
    But to the right, 
    the vector field crosses upwards, in the direction of increasing~$ x_n . $

   Arnold~(Lemma~3, in \textsection~22, of Chapter~3
    from~\cite{arnold1992odes}) gives a simpler affine system, that is topologically equivalent~(homeomorphic) to the given affine ODE~\eqref{eqn:affineODE}.
    For the ODE~\eqref{eqn:affineODE}, consider any positive definite, quadratic Lyapunov function.  Level sets of this function are bounded, non-degenerate ellipsoidal shells. Fix such a level set, and call it~$ {\mathcal{L}} .  $  Let~$ \tau \left( x \right) $ denote the finite
    time it takes for the trajectory starting at a given~$ x $ to cross
    the level set~$ {\mathcal{L}} .  $
This crossing must be transversal, since the Lyapunov function is a strictly decreasing function, along any trajectory.
    With this and Lemma~\ref{lemma:crossingHypersurface}, the crossing time~$ \tau\left( x \right) $ is an analytic function, whenever it is finite.
Consider the invertible map:
    \begin{align*}
        h :  x \in   {\mathbb{R}}^n 
        \to 
        \xi \in   {\mathbb{R}}^n 
        & \triangleq
        e^{ {\mathbf{-}}  \tau\left( x \right) } 
          e^{A \tau\left( x \right) } \left( x -  A^{-1} B \right) 
        + A^{-1} B .
    \end{align*}
    Then clearly~$ h $ is an analytic diffeomorphism, whenever~$ x  $ is not equal to the equilibrium point~$  A^{-1} B . $ 
    
    The flow in the~$ \xi $ space is given by the simpler ODE:
    \begin{align*}
        {\frac{d}{dt}} \xi
        &
        = 
        - \left(  \xi - A^{-1} B   \right) .
    \end{align*}
    Figure~\ref{fig:topologicalEquivalence} depicts the topology of every possible two-dimensional slice, around the mapped crossing point.
    On either side of the mapped crossing point, the vector field crosses the hypersurface~$ {\widetilde{\mathcal{H}}}  $
    in the same direction, w.r.t. the normal vectors on the hypersurface.

    Every analytic diffeomorphism preserves orientedness. 
    The hyperplane~$  {\mathcal{H}}  $ separates the two open, half spaces of~$ {\mathbb{R}}^n . $ 
 These must correspond to the two space regions on either side of the analytic hypersurface~$   {\widetilde{\mathcal{H}}} . $

 But we have topologically contradictory pictures of the vector fields, around the supposed crossing point. Hence there can be no such tangential crossing point.
\end{proof}
\begin{corollary}
If a transfer function is rational, proper and stable, then
its step response cannot have an inflection point that is also a stationary point.
\end{corollary}
\begin{corollary}
\label{theorem:analytic} 
If the plant is stable with a positive DC gain, then the corresponding first exit
time~$ 
    \tau_+ \left( \,\cdot\, ; 1 \right)  
$ and first exit map~$
      \psi_+ \left( \,\cdot\, ; 1 \right)
$
are analytic  everywhere on the switching set.
\end{corollary}

\section{A periodic orbit exists\label{section:existence}}
We show that a RFS has a periodic orbit if the plant is BRL-URF, by showing that for sufficiently large number of iterations~$k,$ the 
map~$ {\mathbf{-}} \psi\left( \,\cdot\, ; k \right)$ has a fixed point.
\subsection{Bounds on the magnitude of the state%
\label{section:ultimateMagnitudeOfstate}
}
Because the matrix~$A$ is Hurwitz~\cite{desoerVidyasagar2009}):
\begin{align*}
    { \left\lVert e^{A t}  \right\rVert }_2  & 
    \le  M_{\rm{initial}} \; e^{  - \sigma_{\rm{slowest}}}
    , \ \text{for all nonegative times}~t ,
\end{align*}
where~$ M_{\rm{initial}}, \; \sigma_{\rm{slowest}}  $ are positive constants, that  depend 
on the matrix~$A.$ In specific, 
\begin{align*}
 \sigma_{\rm{slowest}}  
    & \triangleq
    \epsilon 
    + \min_{1 \le i \le n}{ \left\lvert {\text{real part of $i$th eigenvalue of $A$}}  \right\rvert },
\end{align*}
where~$\epsilon$ is a positive number that can be arbitrarily small.

Then it follows that the trajectories of the RFS are bounded. 
In specific, consider
an initial
state~$x_0 \in {\mathbb{R}}^n .$ Then at any time~$t,$ we have
\begin{align}
    { \left\lVert x \left( t \right)  \right\rVert }_2  & 
    = 
    { \left\lVert 
        e^{A t}  x_0 
        +
        \int_0^t{  e^{A\left( t - s \right)} B \;
                   {\text{sign}} \left( C x \left( s \right) \right)
             ds   }
    \right\rVert }_2 ,  
 \nonumber   \\
    &  \le
    { \left\lVert 
        e^{A t} 
        x_0 
    \right\rVert }_2  
        +
        \int_0^t{ 
    { \left\lVert 
        e^{A\left( t - s \right)} B \;
                   {\text{sign}} \left( C x \left( s \right) \right)
    \right\rVert }_2  
             ds   } ,
 \nonumber   \\
   & \le
      M_{\rm{initial}} \; e^{-\sigma_{slowest} t}  
    { \left\lVert 
        x_0 
    \right\rVert }_2  
    \nonumber \\
    & \quad \quad \quad \quad \quad +
      M_{\rm{initial}} \;
      {\frac{\left(  1 - e^{-\sigma_{\rm{slowest}} t} \right)  }
      {\sigma_{\rm{slowest}}}}
    { \left\lVert 
        B
    \right\rVert }_2  .
\label{eqn:stateMangnitudeBound}
\end{align}
This gives the long-term bound:
\begin{align}
    \lim_{t \to \infty}{   
    { \left\lVert x \left( t \right)  \right\rVert }_2
    }
&  \le
      {\frac {M_{\rm{initial}}}
      {\sigma_{\rm{slowest}}}}
    { \left\lVert 
        B
    \right\rVert }_2  .
\label{eqn:limitingStateMangnitudeBound}
\end{align}
To bound the state after a finite number of switches, we need bounds on the state that hold
after a finite duration of time. Hence we shall use Inequality~\eqref{eqn:stateMangnitudeBound}  to derive the looser bounds given in the following two lemmas.
\begin{lemma} 
    \label{lemma:ultimateBound}
Let the matrix~$A$ be Hurwitz stable. Let%
\begin{align}
    {\widetilde{M}}_{\rm{loose}}  & \triangleq
      \quad   {\mathbf{2}} \times
     {\frac{
      M_{\rm{initial}}
     }
    {\sigma_{\rm{slowest}} }
}
    { \left\lVert 
        B
    \right\rVert }_2  , 
\label{eqn:boundR}
    \\
\intertext{and  let the following spherical ball have a  radius~$ {\widetilde{M}}_{\rm{loose}}$:} 
    {\mathcal{B}}_{\rm{ultimate}} & \triangleq
          \left\{
              x \in {\mathbb{R}}^n : {\left\lVert x \right\rVert}_2 \le {\widetilde{M}}_{\rm{loose}}
          \right\} ,
      \label{eqn:sphericalSetR}
\end{align}
and let the initial state~$x \left( 0 \right) $ lie in the ball~$ {\mathcal{B}}_{\rm{ultimate}}  .  $ 
For some of the time, the 
    state of the RFS described by Equation~\eqref{eqn:RFSdynamics} 
    could have excursions outside 
the ball~$ {\mathcal{B}}_{\rm{ultimate}}  . $
But nevertheless, for all times greater than
\begin{align}
    t_{\rm{excursions-over}}  & \triangleq
     {\frac
    { 1 }
    {\sigma_{\rm{slowest}} }
    }
    \log{ 
    \bigl(  2  M_{\rm{initial}}
    \bigr) 
    }
    ,
    \label{eqn:excursionsOver}
\end{align}
the 
    state of the RFS 
    is confined to the ball~$ {\mathcal{B}}_{\rm{ultimate}}  . $ 
\label{lemma:excursion} 
\end{lemma} 
\begin{lemma} 
Let the matrix~$A$ be Hurwitz stable.
    If the initial state~$ x \left( 0 \right) $ lies in the ball~$ {\mathcal{B}}_{\rm{ultimate}}  , $
then the subsequent maximum magnitude of the state has the upper bound:
\begin{align}
     \max_{t \ge 0}{
    {\left\lVert  x \left( t  \right)  \right\rVert}_2
    }
    & \le \ 
    {\widehat{M}}_{\rm{excursion}}  
  \  \triangleq \
     {\frac
    { M_{\rm{initial}}\left( 2 M_{\rm{initial}}+ 1 \right) }
    {\sigma_{\rm{slowest}} }
    }
    {\left\lVert  B  \right\rVert}_2
    .
\label{lemma:maximumMagnitudeOfStateFromR}
\end{align}
\end{lemma} 

%
\subsection{Inter-switch intervals have a positive minimum duration\label{section:minimumSwitchDuration}}
\begin{figure}
\begin{center}
\begin{tikzpicture}
\begin{axis}[
        axis lines = middle, 
        ticks = none, 
        xmin = -3, xmax = 3, 
        ymin = -1.0, ymax = 1.0, 
        xlabel=$x_{n-1}$,  ylabel=$x_n$,
            x post scale = 1.1, 
            x = 1.2cm, y = 1.4cm,
     ]
\fill [red1!40, fill opacity=0.6] (-3,0) rectangle (3,1.7);
\draw[color=red1,->,>=angle 60, very thick]   (axis cs:-0.5, 0) to (axis cs: -0.55,0.2);
\draw[color=red1,->,>=angle 60, very thick]   (axis cs:0, 0) to (axis cs: -0.1,0.4);
\draw[color=red1,->,>=angle 60, very thick]   (axis cs:0.5, 0) to (axis cs: 0.4,0.5);
\draw[color=red1,->,>=angle 60, very thick]   (axis cs:1, 0) to (axis cs: 0.85,0.6);
\draw[color=red1,->,>=angle 60, very thick]   (axis cs:2, 0) to (axis cs: 1.8,0.9);
\draw[color=red1,->,>=angle 60, very thick]   (axis cs:-1.5, 0.2) to (axis cs: -1.55,0);
\draw[color=red1,->,>=angle 60, very thick]   (axis cs:-2.0, 0.4) to (axis cs: -2.1,0);
\draw[color=red1,->,>=angle 60, very thick]   (axis cs:-2.6, 0.6) to (axis cs: -2.65,0);
\draw[color=tolDarkYellow,->,>=angle 60, very thick]   (axis cs:+0.5, 0) to (axis cs: +0.55,-0.2);
\draw[color=tolDarkYellow,->,>=angle 60, very thick]   (axis cs:0, 0) to (axis cs: 0.1,- 0.4);
\draw[color=tolDarkYellow,->,>=angle 60, very thick]   (axis cs:-0.5, 0) to (axis cs: -0.42,- 0.5);
\draw[color=tolDarkYellow,->,>=angle 60, very thick]   (axis cs:-1, 0) to (axis cs: -0.85,- 0.6);
\draw[color=tolDarkYellow,->,>=angle 60, very thick]   (axis cs:-2, 0) to (axis cs: -1.8,-0.9);
\draw[color=tolDarkYellow,->,>=angle 60, very thick]   (axis cs:-2.6, 0) to (axis cs: -2.62,-1.0);
\draw[color=tolDarkYellow,->,>=angle 60, very thick]   (axis cs:+1.5,- 0.2) to (axis cs: +1.55,0);
\draw[color=tolDarkYellow,->,>=angle 60, very thick]   (axis cs:+2.0, -0.4) to (axis cs: +2.1,0);
\draw[color=tolDarkYellow,->,>=angle 60, very thick]   (axis cs:+2.6, -0.6) to (axis cs: +2.65,0);
\draw[color=blueForRed1, dashed, thick]   (axis cs:-1, 0) to (axis cs: -1,-0.9);
\draw[color=blueForRed1, dashed, thick]   (axis cs:1, 0) to (axis cs: 1,-0.9);
\draw[color=blueForRed1, <-> , >=angle 60, dashed, thick]   (axis cs:-1, -0.65) to (axis cs: 0,-0.65);
\node[color=blueForRed1,scale=1,below]  at (axis cs:-0.5, -0.66) {$  \left\lvert b_{n-1} \right\rvert$};
\draw[color=blueForRed1, <-> , >=angle 60, dashed, thick]   (axis cs:+1, -0.55) to (axis cs: 0,-0.55);
\node[color=blueForRed1,scale=1,below]  at (axis cs:+0.5, -0.58) {$\left\lvert b_{n-1} \right\rvert$};
\draw[color=blueForRed1 , ultra thick]   (axis cs:-1, 0.03) to (axis cs: 1,0.03);
\draw[color=blueForRed1 , ultra thick]   (axis cs:-1, -0.03) to (axis cs: 1,-0.03);
\draw[color=blueForRed1 , ultra thick]   (axis cs:-1, 0) to (axis cs: 1,0);
\node[color=blueForRed1,scale=1.5]  at (axis cs:1, 0) {\Large{$\bullet$}};
    \node[color=blueForRed1,scale=1.5]  at (axis cs:-1, 0) {\Large{$\bullet$}};
\end{axis}
\end{tikzpicture}
\end{center}
    \caption{\label{fig:fieldsAtSwitchingLine}A slice of the vector fields at switching hyperplane}
\end{figure}
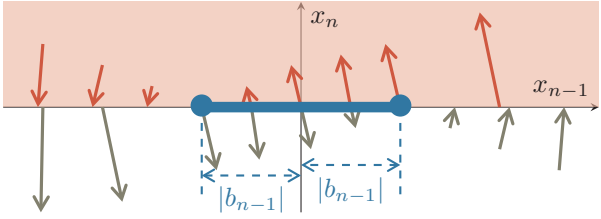
If the plant is BRL-URF, then the coefficient~$b_{n - 1}$ is negative. 
The set we define next is a subset of the switching set, and it gains special significance:
\begin{align}
    {\mathcal{Q}} & \triangleq
    \left\{
          x \in {\mathcal{S}} : 
              \left\lvert x_{n-1} \right\rvert <  
              \left\lvert b_{n-1} \right\rvert
        \right\}
    \label{eqn:evanescentStrip}
\end{align}
We call it the {\textit{evanescent strip,}} because no trajectory from outside the strip lands on it, and the first exit maps take the points on this strip to points on the switching set outside the evanescent strip.
As Figure~\ref{fig:fieldsAtSwitchingLine} shows, the vector field points away from the strip, under both positive and negative state of the relay. In other words:
\begin{itemize}
    \item{no point on the strip~$ {\mathcal{Q}} $ is in the range of the first exit map from positive sign~$ \psi_+\left( \,\cdot\, ; 1 \right) ,$ and,}
    \item{no point on the strip~$  {\mathcal{Q}} $ is in the range of the first exit map from negative sign~$ \psi_-\left( \,\cdot\, ; 1 \right) .$}
\end{itemize}
We can construct the sequence of switching points in the trajectory,
by alternately applying the maps~$ \psi_+\left( \,\cdot\, ; 1 \right) .$ and~$ \psi_-\left( \,\cdot\, ; 1 \right) .$ Then, for any arbitrary trajectory, only the starting point could lie in the strip~$  {\mathcal{Q}} $. No subsequent switching point can lie in the strip.

Therefore, starting from the second interval, over every subsequent inter-switching interval, the trajectory must leap over the strip~$  {\mathcal{Q}} $.   are the times taken by either the affine ODE:
Over such an interval of time, the trajectory must  travel from a point on the switching set lying on one side of the strip~$  {\mathcal{Q}} $, to a point on the switching set lying on the opposite side from the strip.

We can then calculate a positive minimum duration for leaping from one side of the strip to the other side, if we restrict the starting point to be from a bounded set. 

Consider a trajectory that starts at a point which: 
(i)~is in the ball~$ {{\mathcal{B}}_{\rm{ultimate}}} , $  
(i)~also lies on the switching set~$ {\mathcal{S}} , $  
and (iii)~has its~$ \left( n - 1 \right)$th coordinate greater than~$ \left\lvert  b_{n-1}  \right\rvert .$ 
Since~$b_{n-1}$ is negative, the starting and ending points of the trajectory segment have a mutual distance that is bounded below by the shortest of distances between 
points on~$ {\mathcal{S}} $ with~$x_{n-1}$ greater than or equal to~$ \left\lvert  b_{n-1}  \right\rvert ,$ 
and 
 points on~$ {\mathcal{S}} $ with~$x_{n-1}$ lesser than or equal to~$  - \left\lvert  b_{n-1}  \right\rvert . $
 This shortest distance is  bounded below by: 
\begin{gather*}
    2  \left\lvert  b_{n-1}  \right\rvert  .
\end{gather*}
Along any trajectory starting
from inside the ball~$ {{\mathcal{B}}_{\rm{ultimate}}} , $
the magnitude of the velocity 
is bounded above by: 
\begin{gather*}
    {\left\lVert  A  \right\rVert}_2 {\widehat{M}}_{\rm{excursion}}
     +
    {\left\lVert  B  \right\rVert}_2  . 
\end{gather*}
\begin{lemma} 
Consider the RFS described by Equation~\eqref{eqn:RFSdynamics}, where the matrices~$A , B, C$ are given by Equation~\eqref{eqn:companionA} .
Let the plant be BRL-URF.
Let ${\widetilde{M}}_{\rm{loose}} , {{\mathcal{B}}_{\rm{ultimate}}} , {\widehat{M}}_{\rm{excursion}}$ be defined by 
Equations~\eqref{eqn:boundR},~\eqref{eqn:sphericalSetR}, and~\eqref{lemma:maximumMagnitudeOfStateFromR} respectively.
Then for every trajectory of the RFS that starts from the ball~$
  {{\mathcal{B}}_{\rm{ultimate}}}  , $ 
the second and later switching intervals are lower bounded by
    the following:
\begin{align}
    t_{\rm{minimum-inter-switch}}  & \triangleq
    &
  {\frac{
      2  \left\lvert  b_{n-1}  \right\rvert  }
    {  {\left\lVert  A  \right\rVert}_2 {\widehat{M}}_{\rm{excursion}}
       + {\left\lVert  B  \right\rVert}_2  }
  }
  .
    \label{eqn:minimumInterSwitch}
\end{align}
    \label{lemma:minimumInterSwitch}
\end{lemma} 
\subsection{Brouwer's fixed point for the map~$ {\mathbf{-}} \psi_+\left( \,\cdot\, ;  k \right)$}

The Brouwer fixed point theorem 
 guarantees a fixed point for a continuous map that takes a closed, bounded and simply connected subset of an Euclidean space to itself. We apply this theorem to the map~$   {\mathbf{-}} \psi_+\left( \,\cdot\, ;  k \right) .$
\begin{lemma}\label{lemma:fixedPointExists}
Consider the RFS~\eqref{eqn:RFSdynamics} for a plant that is BRL-URF.
There exists a positive integer~$K$ such that for every integer~$k$ that is  greater than or equal to~$ K ,$ the map
    \begin{gather*}
  {\mathbf{-}} \psi_+\left( \,\cdot\, ;  k \right) 
    \end{gather*}
has a fixed point. 
\end{lemma}
\begin{proof}
The set we define next 
lies on the switching set, and is completely 
on one side of the evanescent strip:
    \begin{align}
        {\mathcal{D}} & \triangleq
         \left\{  
         x \in  {\mathcal{S}} 
         :
         \;
           x \in {{\mathcal{B}}_{\rm{ultimate}}} , \ 
           x \notin {{\mathcal{Q}}} , \
           x_{n-1}  \ge 0
         \;
         \right\} .
        \label{eqn:theSetD}
    \end{align}
    This is the set of points on the switching set that: (i)~have lengths less than or equal to~${\widetilde{M}}_{\rm{loose}},$ (ii)~are not on the evanescent strip, and (iii)~have a positive value for the $\left( n - 1 \right)$th coordinate.
In other words, 
if $  x \in  {\mathcal{D}}  ,$ then
    \begin{gather*}
       {\left\lVert  x  \right\rVert}_2  \le {\widetilde{M}}_{\rm{loose}} ,
        \
         x_{n}  = 0 ,  
        \
        \text{and,}
        \ \
         x_{n-1}  \ge
        \left\lvert  b_{n-1}  \right\rvert
        .
    \end{gather*}
For any trajectory that starts from~$ {\mathcal{D}} , $ 
    Lemma~\ref{lemma:excursion} implies that the  magnitude of switching points is less than or equal to~$  {\widetilde{M}}_{\rm{loose}} , $ if the number of switches has been at least:
    \begin{align}
        K & \triangleq
             \left\lceil
               {\frac
               {t_{\rm{excursions-over}}}
               {t_{\rm{minimum-inter-switch}}
               }
               }
             \right\rceil
            ,
            \label{eqn:K}
    \end{align}
    where the times~$  {t_{\rm{excursions-over}}} , \; {t_{\rm{minimum-inter-switch}}} $ are defined by Equations~\eqref{eqn:excursionsOver},~\eqref{eqn:minimumInterSwitch}.

Hence 
    if~$  
    k \ge K , 
     $ then
~\[ 
  x 
     \in
     {\mathcal{D}} 
    \quad
   \implies
    \quad
 - \psi_+ \left( x ; k \right)
     \in
     {\mathcal{D}} .
    \]
%
%
Because~$ {\mathcal{D}} $ is bounded and convex, and because 
the  map~$  \psi_+ \left( \,\cdot\, ; k \right)  $ is continuous, 
the lemma follows from Brouwer's fixed point theorem.
\end{proof}
\begin{theorem}
If the plant in the RFS~\eqref{eqn:RFSdynamics} is BRL-URF, then the RFS possesses at least one  periodic orbit.
\end{theorem}
\begin{proof}
    Pick an integer~$ k $ that is greater than or equal to~$ K .$
Since 
    $
    {\mathbf{-}} \psi_+ \left( \,\cdot\, ; k \right)
    $
has a fixed point in the set~$ {\mathcal{D}} $, 
    the RFS admits a closed orbit that passes through the set~$ {\mathcal{D}} $.

The number of distinct switching points that fall within a full traversal of this closed orbit is less than or equal to~$ 
        2 k . 
        $
    We view the orbit as the trajectory that starts at the Brouwer fixed point in the set~$ {\mathcal{D}} .$  Then by Lemma~\ref{lemma:minimumInterSwitch}, the minimum time interval between two successive switches is at least~$ t_{\rm{minimum-inter-switch}} . $

    Similarly, the maximum interval between two successive switches can be bounded above. 
 For any trajectory that starts from  the set~$ {\mathcal{D}} , $  
Lemma~\ref{lemma:maximumMagnitudeOfStateFromR} guarantees that all  
subsequent
    switching points have magnitudes less than or equal to~$  
{\widehat{M}}_{\rm{excursion}} .
    $

    Because the plant is BRL-URF, the first exit time~$ \tau_+\left( \cdot \right) $ is continuous everywhere on the switching plane~$ {\mathcal{S}} $. 
    Then by the Weierstrass theorem for continuous functions on closed and bounded sets, the function~$ \tau_+ \left(  \,\cdot\, \right) $ attains its maximum over the following closed and bounded set:
    \begin{gather*}
        {\mathcal{S}} \ \cap \
        \left\{  x \in {\mathbb{R}}^n : 
       {\left\lVert  x  \right\rVert}_2  \le {\widehat{M}}_{\rm{excursion}} ,
        \ \text{and} \ x_{n-1}  \ge
        \left\lvert  b_{n-1}  \right\rvert
        \right\} .
    \end{gather*}
Therefore the maximum inter-switch time must be  finite. 

    Since the orbit is traversed in finite time, and since it has a finite number of switches, and also because any two successive switches have a positive duration between them, it follows that this is a regular periodic orbit.
\end{proof}

\subsection{Examples of locally stable Poincaré maps\label{section:numericalExploration}}
\begin{figure}
\centering
    \subcaptionbox*{}%
    {\includegraphics[width=0.96\linewidth,trim={26mm 0mm 10mm 104mm}, clip ]{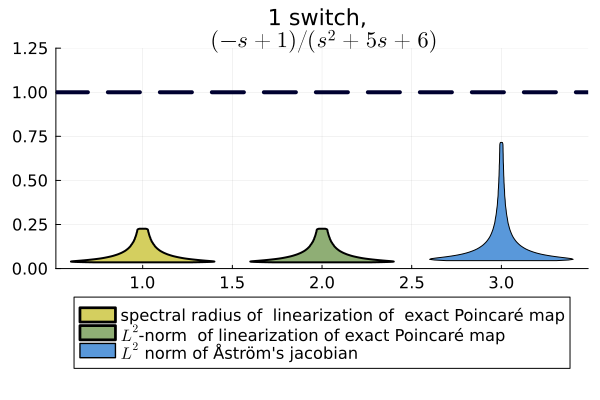}}%
\vspace{-30pt}
\\    
    \subcaptionbox{Both local Schur stability and local contractivity hold%
\label{fig:numericalLocalStabilitySecondOrderPlant}}
    {\includegraphics[width=0.66\linewidth,trim={0mm 44mm 0mm 0mm}, clip ]{figuresAndData/posit_invariant_SecondOrderNumMinus1Plus1Den1Plus5Plus6.png}}%
    \hfill 
    \subcaptionbox{Error is low%
\label{fig:numericalLocalStabilitySecondOrderPlant_bauer_fike}}
    {\includegraphics[width=0.32\linewidth, ]{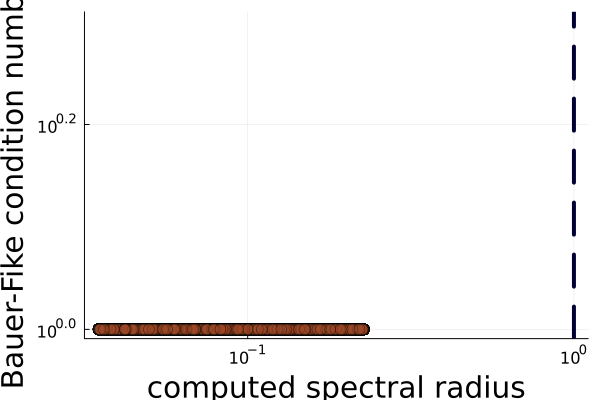}}%
\\    
    \subcaptionbox{Local Schur stability mostly  holds, 
    but fails  at some points
\label{fig:numericalLocalStabilityFourthOrderPlant1}}
    {\includegraphics[width=0.66\linewidth, trim={0mm 44mm 0mm 0mm}, clip ]{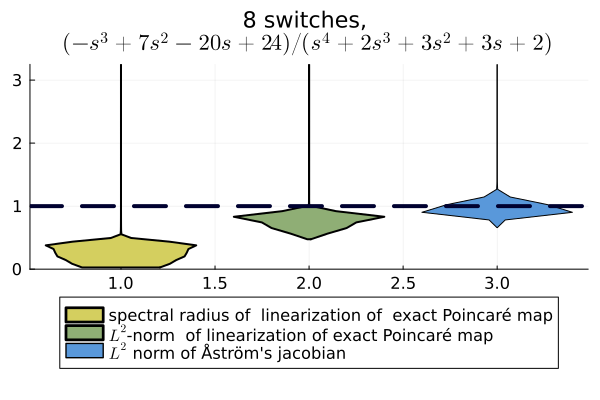}}%
    \hfill 
    \subcaptionbox{Error is low%
\label{fig:numericalLocalStabilityFourthOrderPlant1_bauer_fike}}
    {\includegraphics[width=0.32\linewidth]{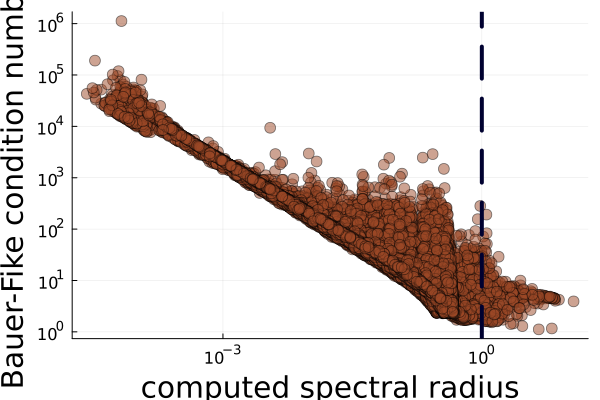}}%
\\    
    \subcaptionbox{Local Schur stability holds%
\label{fig:numericalLocalStabilityThirdOrderPlant1}}
    {\includegraphics[width=0.66\linewidth, trim={0mm 44mm 0mm 0mm}, clip ]{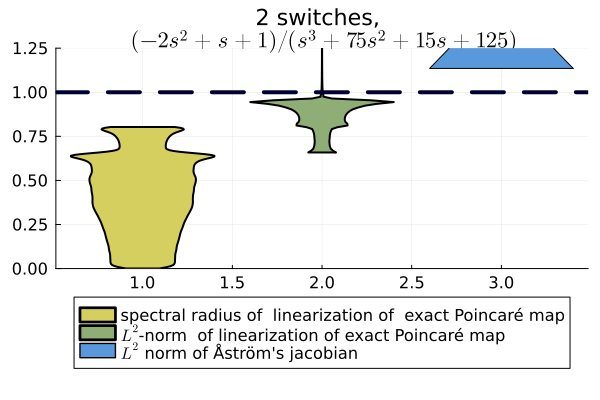}}%
\hfill 
    \subcaptionbox{Error is low
\label{fig:numericalLocalStabilityThirdOrderPlant1_bauer_fike}}
    {\includegraphics[width=0.32\linewidth]{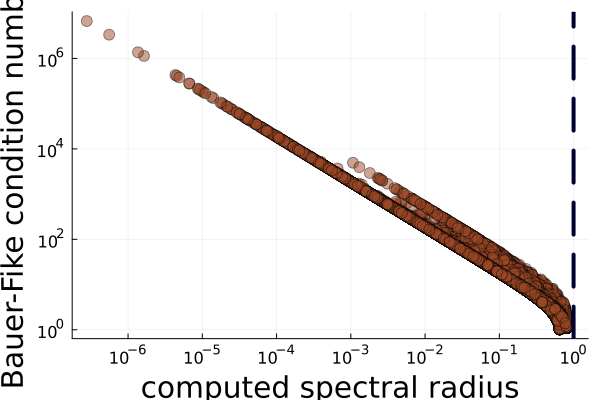}}%
\\    
    \subcaptionbox{Both local Schur stability and local contractivity hold%
\label{fig:numericalLocalStabilityThirdOrderPlant2}}
    {\includegraphics[width=0.66\linewidth, trim={0mm 51mm 0mm 0mm}, clip ]{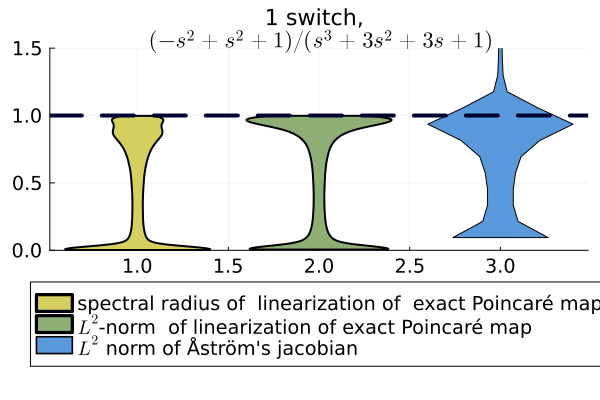}}
\hfill
    \subcaptionbox{Error is low
\label{fig:numericalLocalStabilityThirdOrderPlant2_bauer_fike}}
    {\includegraphics[width=0.32\linewidth]{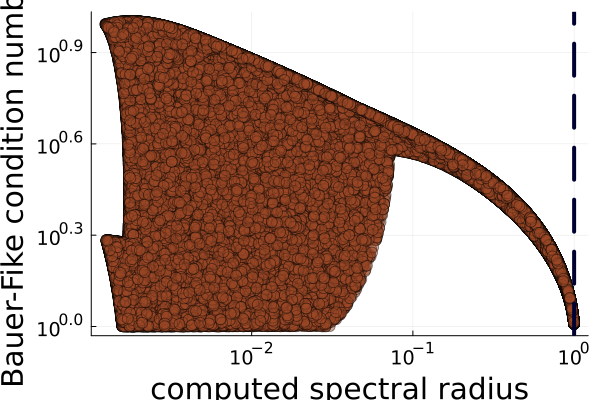}}
%
    \\    
    \subcaptionbox{Local Schur stability holds%
    \label{fig:numericalLocalStabilityThirdOrderPlant2_second}}
    {\includegraphics[width=0.66\linewidth, trim={0mm 51mm 0mm 0mm},   clip, ]{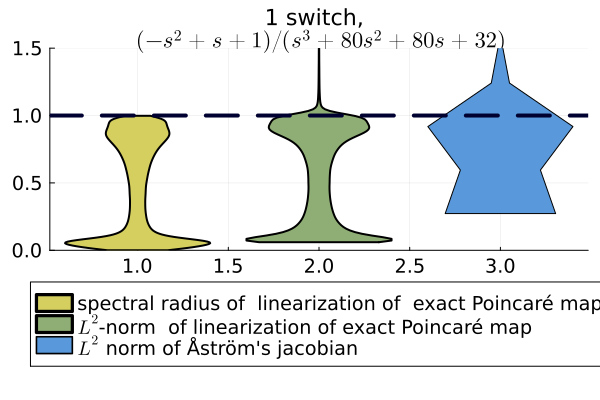}}%
\hfill
    \subcaptionbox{Error is low%
\label{fig:numericalLocalStabilityThirdOrderPlant2_second__bauer_fike}}
    {\includegraphics[width=0.30\linewidth,  ]{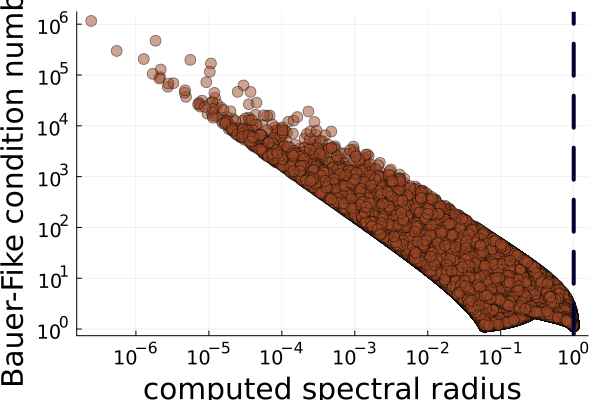}}%
\\    
    \subcaptionbox{Local Schur stability mostly holds%
\label{fig:numericalLocalStabilityThirdOrderPlant3}}
    {\includegraphics[width=0.66\linewidth, trim={0mm 44mm 0mm 0mm}, clip ]{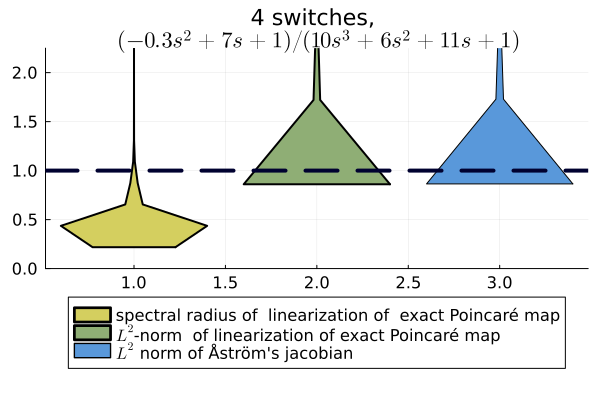}}%
\hfill 
    \subcaptionbox{Error is not low
\label{fig:numericalLocalStabilityThirdOrderPlant3_bauer_fike}}
    {\includegraphics[width=0.32\linewidth]{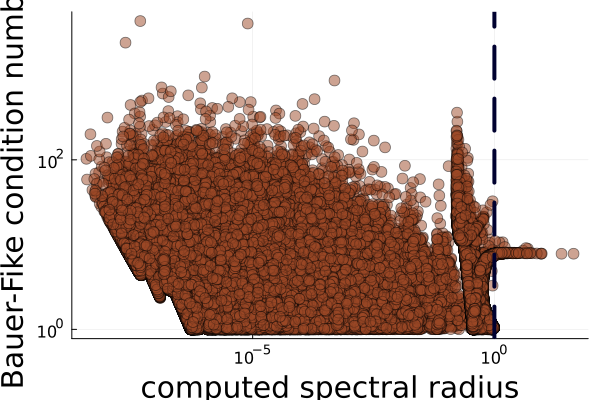}}%
 \caption{Examples of violin plots of the empirical distribution of spectral radius and $L^2$~norm, accompanied by scatter plots of the condition numbers
    - based on a random sample of $10^6$ points from the set~$ {\mathcal{D}} . $
    In these examples, this set happens to be positively invariant for the map~$ \psi\left( \cdot ; k \right) $, for the given small number~$k$ of iterations~(switches).%
\label{fig:violinPlotsSpectralRadius}}
\end{figure}
Simple conditions are available for  the parameters of possible periodic orbits. The half period~$ \tau_+ \left( {\widehat{x}}  \right) ,  $ and the real vector~$ {\widehat{x}} $ correspond to
a symmetric and unimodal periodic orbit, if and only if $ C {\widehat{x}}  = 0 , $
and~\cite{astrom1995oscillationsRelay}:
\begin{gather*}
    \begin{cases}
   C   
        { \left(   e^{ A \tau_+ \left( x \right)  }  + I \right)}^{-1}
        { \left(   e^{ A \tau_+ \left( x \right)   }  - I \right)}
         A^{-1} B 
           = 0  ,  \ \text{and,}&\\
    C  e^{ A t  }  ( {\widehat{x}}   - A^{-1} B )
        + C A^{-1} B   \ge 0, \ \text{if}\ 0 \le t  \le 
        \tau_+ \left( {\widehat{x}} \right)   . &
    \end{cases}
\end{gather*}
For stable second order plants, the first exit map is locally stable, and
even contracting~\cite{rabi2018relay}.
But for higher order plants, we do not know of any guarantees that some Poincaré map is either contracting, or at least 
locally Schur stable everywhere on some positively invariant set.
Numerical experiments indicate that at least for some plants or order three and four, Poincaré maps are locally Schur stable and positively invariant over the set~$  {\mathcal{D}} . $ We first describe the derivatives that we compute numerically.

We have two geometric notions~(see~Appendix~A) of the derivative 
of the Poincaré map.
Both derivatives are calculated by applying the implicit function theorem.

{\AA}str{\"o}m's jacobian~\cite{astrom1995oscillationsRelay} has the expression:
\begin{align}
    J_n \left( x \right) & =  \left( I - {\frac{1}{C u}}u C \right) e^{A\tau_+\left( x \right)}, 
    \label{eqn:astromjacobian}
\end{align}
where the column vector~$ u =
    A e^{ A \tau_+ \left( x \right)} ( x - A^{-1} B )
,$  which is the direction of the vector field at the point of crossing the 
switching hyperplane, and the row vector~$ C $ is 
the normal to the switching hyperplane, and is given by~\eqref{eqn:companionA}.

The exact derivative of the Poincaré map is:
\begin{align}
    J_{n - 1} \left( x \right) & =  
    \left( I - {\frac{1}{C u}}u C \right) 
    e^{A\tau_+\left( x \right)}
    \left( I - {\frac{1}{v^T v}} v v^T \right) ,
    \label{eqn:exactjacobian}
\end{align}
where the column vector~$ v =
    A  ( x - A^{-1} B )
,$  which is the direction of the vector field at the starting point~$ x . $

Both derivatives have the same spectral radius~(see~Appendix~A). But the $ L_2 $~norm of {\AA}str{\"o}m's
jacobian is bigger than or equal to that of the exact jacobian.


For some examples of BRL-URF plants, 
Figure~\ref{fig:violinPlotsSpectralRadius} describes the empirical 
distributions of the spectral radius, and the $ L_2 $~norms of the two   jacobians. 
The sensitivity of a calculated spectral radius is the Bauer-Fike condition number, which is
the condition number of the square matrix created by assembling all the eigenvectors.

In our examples, the Bauer-Fike condition number can be high. 
But for those eigenvalues that are close to or exceed one,
the corresponding Bauer-Fike condition numbers are quite low.
Hence Schur stability or instability can be concluded with reasonable confidence.

\section{Local stability implies Global stability\label{section:globalStability}}

\subsection{Why we embed the dynamics in~${\mathbb{C}}^n$}
We cannot generally guarantee global asymptotic stability of a fixed point, even for a smooth map on a bounded Euclidean set that possesses the following two properties:  (a)~the map meets the requirements of the Brouwer fixed point theorem, and (b)~the map has a Schur stable jacobian everywhere.
The following example from Shih and Wu~\cite{shihWu1998studiaMathematica} illustrates this.
\begin{example}[Schur stable jacobian everywhere does not guarantee global stability of Brouwer's fixed point]
    Consider the domain~$E \subset {\mathbb{R}}^2 ,$ which is defined as the following square shaped set,  centered at the origin:
    \begin{gather*}
        \left\{ \left(x_1 , x_2 \right) \in  {\mathbb{R}}^2 
     \mathbf{:}
         \left\lvert x_1 \right\rvert 
         +
         \left\lvert x_2 \right\rvert 
         \le 
         1
    \right\}   .
\end{gather*}
Consider the map~$ f :  {\mathbb{R}}^2 \to {\mathbb{R}}^2 $ defined as follows:
    \begin{align*}
        f \left(  x_1 , x_2 \right) & = 
        \left(  {\widetilde{\mu}}\left(x_2\right) , {\widetilde{\mu}}\left(x_1\right)  \right) ,
    \end{align*}
where, ${\widetilde{\mu}} : {\mathbb{R}} \to  {\mathbb{R}}$ is any function that satisfies:
    \begin{gather*}
        {\widetilde{\mu}}\left( \alpha \right) = 0, \ {\text{if}} \ \left\lvert \alpha \right\rvert \le 1 / 2 ,
        \quad \quad
        {\widetilde{\mu}}\left( 1 \right) =  {\widetilde{\mu}}\left( -1 \right) = 1, 
 \quad \quad
        \\
 \left\lvert  {\widetilde{\mu}}\left( \alpha \right) \right\rvert \le 1,
        \ {\text{if}} \ \left\lvert \alpha \right\rvert \le 1 .
    \end{gather*}
One such smooth map~${\widetilde{\mu}}$ is the following~$C^4$~function:
\begin{align*}
    {\widetilde{\mu}}\left( \alpha \right) &
    =
    \begin{cases}
        0 ,  & {\text{if}} \ \left\lvert \alpha \right\rvert \le 1 / 2 , \\
        2^5 { \left(  \left\lvert \alpha \right\rvert - 1 / 2 \right)}^5 , & {\text{if}} \ \left\lvert \alpha \right\rvert \ge 1 / 2 .
    \end{cases}
\end{align*}
Other alternatives exist for~${\widetilde{\mu}}$ that are arbitrarily smooth. Regardless of the specific form of~${\widetilde{\mu}} , $ we have the following: 
\begin{itemize}
    \item{Brouwer's fixed point exists, because $f\left( E \right) \subset E,  $}
    \item{the origin is the only fixed point of the function~$f , $}
    \item{the jacobian at any point~$ \left( x_1 ,  x_2 \right)$ from $E$ is given by:
        \begin{align*}
            f_x  \left( x_1 , x_2 \right)  & = 
            \begin{pmatrix}
                0   
                & {\widetilde{\mu}}_{x_2}\left( x_2 \right)   \\
                {\widetilde{\mu}}_{x_1}\left( x_1 \right)  
                & 0
            \end{pmatrix}
            .
        \end{align*}
        Both eigenvalues of this matrix are zero~!
        }
    \item{But we have at least one cycle of period~1, because: $ 
        f \left(  0 , 1 \right) =  \left(  1 , 0 \right)  $
        and $
        f \left(  1 , 0 \right) =  \left(  0 , 1 \right)  . $
        }
\end{itemize}
Therefore the origin does not attract all points on~$E.$ 
Clearly differentiability and Schur stability everywhere on the domain of the map is not sufficient to guarantee the global stability of Brouwer's fixed point.
 In specific, for the particular~${\widetilde{\mu}}$ function mentioned above, the origin can only attract those points in the smaller square:~$
        \left\{ \left( x_1 , x_2 \right) \in
        {\mathbb{R}}^2 
     \mathbf{:}
         \left\lvert x_1 \right\rvert 
         \le 1 / 2 , \ {\text{and}} \
         \left\lvert x_2 \right\rvert 
         \le  1 / 2
    \right\}   .
$
\end{example}
We can guarantee the global stability of Brouwer's fixed point, if we have in addition, the holomorphy of the mapping, whose fixed point we are interested in. 
\subsection{Shih-Wu fixed point theorem for bounded, analytic maps}
    Next comes a theorem that guarantees the global stability of a locally stable fixed point, when Brouwer's fixed point theorem is applied on analytic map acting on a bounded, convex domain.   
    Because we depend on this theorem to study the stability of our RFS limit cycles, we should understand why this result holds. Therefore we repeat below both the theorem and the proof that is given in~\cite{shihWu1998studiaMathematica}. 
\begin{theorem}[the finite-dimensional special case of Theorem~1 from Shih and Wu~\cite{shihWu1998studiaMathematica}]
    \label{theroem:shihWu}
    Let the domain~$\Delta \subset {\mathbb{C}}^n , $ and assume that:
    \begin{itemize}
        \item{$\Delta$ is non-empty, closed and bounded,}
        \item{the interior of $\Delta$ is an open subset of  ${\mathbb{C}}^n$, and}
        \item{$\Delta$ is convex.}
    \end{itemize}
%
%
Let the map~$ f :  {\mathbb{C}}^n  \to {\mathbb{C}}^n ,$ 
and assume that:
    \begin{itemize}
        \item{$f$ is analytic on~$\Delta,$}
        \item{the sequence of iterates of~$f$ is a uniformly bounded sequence of maps, on the domain~$\Delta$, and,}
        \item{there is a fixed point, where the jacobian 
              matrix~$f_z\left( \,\cdot\, \right)$ is Schur stable.}
    \end{itemize}
    Then
    the said fixed point 
    is globally asymptotically stable. 
\end{theorem}


\begin{proof} 
Essentially the proof consists of showing that 
    at least a subsequence of the iterated maps~$f^k$ must converge to 
a constant map - in specific, the constant map that maps the whole of the domain~$\Delta$ to a single fixed point.
We have laid out the proof in the following three steps:
\begin{description}
    \item[Step 1:]{~Showing that a locally stable fixed point exists,}
    \item[Step 2:]{~Showing that 
        the sequence of iterated  maps:
        $ \left\{ f \left( \,\cdot\, \right) ,    f^2 \left( \,\cdot\, \right) ,  \ldots  \right\}  $ has an infinite subsequence that converges to a analytic map, and,}
    \item[Step 3:]{~Showing that the sequence of iterated maps:
        $ \left\{ f \left( \,\cdot\, \right) ,    f^2 \left( \,\cdot\, \right) ,  \ldots  \right\}  $ converges to a constant map, which maps the whole of~$\Delta$ to the fixed point.}
\end{description}
    \noindent
    {\textcolor{cyanTeal}{\textit{{\colorbox{tolPaleYellow}{Step 1}}: Showing that  a locally stable fixed point exists}}}


    Denote this fixed point by~${\widehat{z}}.$
    Since the jacobian matrix is Schur stable, it follows that there is a small open set~$U$ around the fixed point~${\widehat{z}}, $ such that for any point in~$U$, repeated iteration of~$f$ converges to
    the fixed point~${\widehat{z}}.$

    \noindent
    {\textcolor{cyanTeal}{\textit{{\colorbox{tolPaleYellow}{Step 2}}: 
    Showing that 
        the sequence of iterated  maps:
        $ \left\{ f \left( \,\cdot\, \right) ,    f^2 \left( \,\cdot\, \right) ,  \ldots  \right\}  $ 
        has an infinite subsequence that converges to a analytic map}}}

In this step, we use results from Functional analysis that concern the convergence of sequences of maps that are both uniformly bounded and analytic.

Since every map in the sequence
        $ \left\{ f \left( \,\cdot\, \right) ,    f^2 \left( \,\cdot\, \right) ,  \ldots  \right\}  $ 
 takes values inside the bounded set~$\Delta,$ it follows that this sequence is uniformly bounded.

    For any bounded domain in~${\mathbb{C}}^n ,$ the space of continuous mappings from it to itself, is a complete metric space under the supremum norm~$ {\lVert \,\cdot\, \rVert}_{\infty} $.    
    In any complete metric space, every bounded sequence 
    has a convergent subsequence. 
  
   Hence  it follows that our sequence of iterated maps:
    $ \left\{ f \left( \,\cdot\, \right) ,    f^2 \left( \,\cdot\, \right) ,  \ldots  \right\}  $ 
    has a convergent subsequence.
In fact such a subsequence must converge to a analytic limit, as we shall see next.

    Any uniformly bounded sequence of analytic maps is uniformly Lipshitz, and therefore {\textit{uniformly equicontinuous}}\footnote{A sequence of mappings is said to be uniformly equicontinuous if for any arbitrary positive number~$\epsilon,$ it is possible to choose a corresponding positive number~$\delta,$ such that every map in the sequence satisfies the following closeness condition: every pair of points that are within a distance of $\delta$ from each other, is mapped to a pair of image points that are with a distance of~$\epsilon$ from each other.}  (see Theorem~6.1 and Corollary~6.1 of~\cite{khatskevichShoiykhet1994differentiableOperatorsAndNonlinearEquations}).

By the Arzel{\`a}-Ascoli theorem, 
if a sequence of mappings is uniformly equicontinuous, then it has a subsequence that converges uniformly on every closed and bounded subset of the domain. 
For us this means that there is an infinite sequence of increasing positive integers:
$
 \left\{  
   k_1, k_2, \ldots
    \right\}
$
such that the sequence of iterated mappings:
    \begin{gather*} 
    \left\{ f^{k_1} \left( \,\cdot\, \right) ,    f^{k_2} \left( \,\cdot\, \right) ,  \ldots  \right\}  
    \end{gather*} 
converges uniformly over the set~$\Delta.$

    On an $n$-dimensional complex domain\footnote{The equivalent of this statement for Real Euclidean domains is false. The Weierstrass nondifferentiable function~\cite{wikipediaWeierstrassNonDifferentiableFunction} is an example of a Fourier series converging uniformly on a finite interval, but to a limit that is continuous but not differentiable.},
 if a  sequence of analytic maps is uniformly convergent, then its
 limit must also be 
 analytic (see~\cite{krantz2013limitsOfHolomorphicSequences}).
 Thus the limit: 
    \begin{align*}
        f_{\infty} \left( \,\cdot\, \right) & \triangleq  \lim_{k_i \to\infty} f^{k_i} \left( \,\cdot\, \right) ,
    \end{align*}
    must be analytic on  the domain~$\Delta.$

    \noindent
    {\textcolor{cyanTeal}{\textit{{\colorbox{tolPaleYellow}{Step 3}}: 
    Showing that 
        the sequence of iterated  maps:
        $ \left\{ f \left( \,\cdot\, \right) ,    f^2 \left( \,\cdot\, \right) ,  \ldots  \right\}  $ converges to the constant map that maps every point in~$\Delta$ to the fixed point.}}}

        In this step, we shall use the Identity theorem for analytic mappings, to deduce that the limiting map~$ f_{\infty} $ is constant over the whole of the domain~$\Delta.$

Recall from Step~1 that the open set~$U$ around the fixed point~$ {\widehat{z}}, $ is a basin of attraction for~${\widehat{z}}.$ Therefore 
\begin{align*}
    f_{\infty} \left( z \right) & \equiv {\widehat{z}} \
      \text{for every point} \ z \ \text{in the open set} \  U.
\end{align*}

Analytic functions behave locally like polynomials.
If a polynomial takes a fixed value over a finite interval,
then all of its coefficients for positive powers must be necessarily zero.
In other words, such a polynomial function must be a constant function.
Analytic maps possess a similar property.  
The identity theorem
for several complex variables states that over a $n$-dimensional complex domain, any analytic function that takes a fixed value over an open subset must be a constant function over the whole domain. 
This is essentially because constancy over an open subset  means that derivatives of all orders are zero, and this 
reduces the map's power series to the constant term. 

Hence it follows that the limit~$f_{\infty}$ is the constant map:
\begin{align*}
    f_{\infty} \left( z \right) & \equiv {\widehat{z}} , \ {\text{for all}} \ z \in \Delta .
\end{align*}
We shall now show that every infinite subsequence of the iterates of~$f$ converges to the above limit.
Consider a number~$k_i$ in the sequence~$
    \left\{ {k_1} ,    {k_2} ,  \ldots  \right\}  
.
$
The larger~$k_i$ is, the more closely the
map~$f^{k_i} \left( \,\cdot\, \right)$ approximates
the 
map~$f_{\infty} \left( \,\cdot\, \right).$
In fact as~$k_i$  goes to infinity, the approximation error, namely
$  f^{k_i} \left( \,\cdot\, \right) -  {\widehat{z}} , $
converges to zero, uniformly over the domain~$\Delta.$
Hence for all large enough~$ k_i $
\begin{align*}
    f^{k_i}\left( z \right) & \in U , \ {\text{for all}} \ z \in \Delta , \ \text{and,} 
\\
\implies    f^{k_i + k }\left( z \right) & \in U , \ {\text{for every}} \ k > 0 ,
\end{align*}
which is to say that starting from large enough~$k_i$
all further iterates of~$f$ must also stay inside the open set~$U,$ and converge to the fixed point~$ {\widehat{z}} . $
%
\end{proof}


\subsection{The Poincaré map has a globally convergent fixed point}

\begin{theorem}[global convergence to fixed point of Poincaré map]
\label{theorem:globallyStableFixedPoint}
    Let the closed, bounded and convex  set~$ {\mathcal{D}} \subset {\mathcal{S}} \subset {\mathbb{R}}^n $ be described by Equation~\eqref{eqn:theSetD}. 
     Let the integer~$ K  $ be as per Equation~\eqref{eqn:K}, 
    so that the Poincaré map~$  \psi_+\left( \xi ; K \right) $ takes the set~$  {\mathcal{D}} $ to itself.
    Suppose that the Poincaré map~$  \psi_+\left( \xi ; K \right)     $  is locally Schur stable, everywhere on an open set that has the set~$ {\mathcal{D}}  $ in its interior. By this we are merely assuming that the said Poincaré map is locally Schur stable at every point on~$ {\mathcal{D}} , $ and also at those points outside~$ {\mathcal{D}}    $ that are within a short prescribed distance from its boundary.

    Then every point of~$ {\mathcal{D}}  $ is taken by repeated iterates of the Poincaré map~$  \psi_+\left( \xi ; K \right)     $  to the map's  Brouwer fixed point.
\end{theorem}
\begin{proof}
    We apply the Shih-Wu fixed point theorem for: 
    \begin{align}
        f \left( z \right) & = {\mathbf{-}} \psi\left( z ; K \right), 
    \end{align}
    on a suitable domain~$ \Delta  $ in~$ {\mathbb{C}}^n . $ For our choice of domain, it shall be easy to prove that the map~$ f $ is analytic on the domain. But the bulk of this proof shall be devoted  to carefully showing that the map~$ f $ takes points on the domain, back inside the domain. By showing that, we shall show that the iterates of~$ f $ are uniformly bounded.

    {\textit{Notational conveniences:}}
    (i)~given a complex vector~$ z \in {\mathbb{C}}^n , $ and a
    set of complex vectors~$ X \subset {\mathbb{C}}^n , $ we denote by
    $ z + X $  that set of complex vectors in  $ {\mathbb{C}}^n , $ formed by adding~$ z $ to members of the set~$ X, $ and similarly 
    (ii)~given a positive real number~$ \alpha , $ and a
    set of complex vectors~$ X \subset {\mathbb{C}}^n , $ we denote by
    $ \alpha \cdot X $  that set of complex vectors in  $ {\mathbb{C}}^n , $ formed by scaling~$ z $ members of the set~$ X, $ by the factor alpha,
    (iii)~given a complex vector~$ z \in {\mathbb{C}}^n  $ we denote by~$ 
       {\left\lvert z \right\rvert}_{\infty}     
    $ the infinity norm of the vector, which is the largest of the magnitudes of its~$ n $ component complex numbers,
    (iv)~given a complex vector~$ z \in {\mathbb{C}}^n  , $ and a
    set of complex vectors~$ X \subset {\mathbb{C}}^n , $ we denote by
    $ \left\lvert z -  X \right\rvert $  the infimum of all infinity norms of complex vectors  formed by subtracting from~$ z , $ members of the set~$ X, $ ; similarly, $ {\left\lvert z -  X \right\rvert}_{\infty} $  denotes the infimum of norms of the differences between~$ z $ and members of~$ X  , $ and
    (v)~for a complex vector~$ z \in {\mathbb{C}}^n  , $ we denote by~$\Re{ z } $ and $ \Im{ z } $
    the two real vectors of of the real parts and imaginary parts 
    of~$ z . $

    Our choice of~$ \Delta   $ shall 
    contain the real set~$ {\mathcal{D}} $ in its interior.  Let~$ \epsilon , \eta $
    be  two positive real numbers of small magnitude. Consider the bounded domain 
    formed by taking the bounded real set~$ {\mathcal{D}} $ and slightly thickening it along all the real and imaginary axes:
    \begin{align*}
        {\widetilde{\mathcal{D}}}  & \triangleq 
         \left\{
             x \in  {\mathbb{R}}^n :
             \left\lvert \Re{ x }  - {\mathcal{D}} \right\rvert  \le
             \epsilon 
             \right\} , \\
        \Delta  & \triangleq 
         \left\{
             z \in  {\mathbb{C}}^n :
        { \bigl\lvert \Im{  z  } - {\widetilde{\mathcal{D}}} \bigr\rvert }_{\infty}  \le  \eta
             \right\} .
    \end{align*}
    In directions parallel to the imaginary axes, we have thickened up to an  infinity norm constraint.
    Hence the boundary~$ \partial\Delta $ has two flat faces perpendicular to each imaginary axis. This is used in Steps~2, 3 below.

The rest of the proof consists of showing that:
\begin{description}
    \item[Step 1:]{~$ f $ is analytic on~$  \Delta  , $}
    \item[Step 2:]{~$ f $ maps small Schur ellipsoids around $ x \in  {\mathcal{D}}   $ to smaller Schur ellipsoids around~$ f (x)  , $ and,}
    \item[Step 3:]{~$ f $ maps even the boundary points of~$ \Delta $ into~$ \Delta . $}
\end{description}
    \noindent
    {\textcolor{cyanTeal}{\textit{{\colorbox{tolPaleYellow}{Step 1}}: $ f $ is analytic everywhere on~$ \Delta $}}}
  The domain~$ \Delta $ inherits from the real set~$ {\mathcal{D}}   $   the Brouwer fixed point of~$ f  . $ 
   We can choose~$ \epsilon , \eta $ small enough so that~$ f $ is analytic on the set of $ \epsilon - $ and $ \eta -$ sized open balls,  around points of the real set~$ {\mathcal{D}}  . $ The reason is that the set of singularities of~$ f $ is at a non zero distance from the closed set~$  {\mathcal{D}}  . $ 

    We can choose~$ \epsilon $ small enough so that the Schur stability of the jacobians~$ f_z ( x ) $ on the set~$ {\mathcal{D}} $ also holds on the set~$   {\widetilde{\mathcal{D}}}   . $   

    \noindent
    {\textcolor{cyanTeal}{\textit{{\colorbox{tolPaleYellow}{Step 2}}  
        {~$ f $ maps small Schur ellipsoids around $ x \in  {\mathcal{D}}   $ to smaller Schur ellipsoids around~$ f (x)  $}%
      }}}
    At every~$ x \in
   {\mathcal{D}}  
    , $ Schur stability of~$ f_z \left(  x \right)  $ implies 
    that 
    a symmetric, positive definite matrix~$ V\left( x \right) $ exists satisfying the Schur equation:
    \begin{align*}
        {\left[ f_z \left( x \right) \right]}^T \cdot V\left( x \right)  \cdot f_z \left( x \right)  - V \left( x \right) & = - I .
    \end{align*}
    The sublevel sets of the positive definite quadratic form~$ 
      z^H V \left( x \right) z
    $ are ellipsoids. We shall consider the action of the map~$ f $ on such ellipsoids, around points~$ x \in {\mathcal{D}}  . $ We set the sizes of these ellipsoids so that they graze two parallel flat faces of the boundary~$ \partial\Delta . $ One of those faces is the following set of all points in~$ \Delta, $ such that the last coordinate has the maximum possible imaginary part: 
    \begin{align*}
        {\mathcal{F}} 
        &
        \triangleq 
         \left\{
             z \in  \Delta :
        { \Im{ C z  }  }  = + \eta
             \right\} .
    \end{align*}
    At each point~$ x \in  {\mathcal{D}}  $ define the Schur ellipsoid~$ {\mathcal{E}}_x   , $ as:
    \begin{align*}
       {\mathcal{E}}_x   &  \triangleq
        {\left\{
             z \in  {\mathbb{C}}^n :
             z^H V \left( x \right)  z \le r\left( x \right), 
            \right\}}
    \end{align*}
    where the positive real number~$ r \left( x \right)  $ is chosen such that 
    $ \max_{ z \in {\mathcal{E}}  } \Im{ C z  } =  + \eta . $ 

   For the complex vector~$ \delta z  $ from the Schur ellipsoid~$ {\mathcal{E}}_x   , $  we have the first order Taylor expansion:
       \begin{align*}
         f\left( x + \delta z \right)   & = 
         f\left( x \right)   + 
           f_z\left( z  \right) \delta z + o\left( \eta \right) \delta z  .
       \end{align*}
   For small enough~$ \eta $, the Schur stability of~$ f_z (x) $ implies: 
   \begin{align*}
       f \left(  x + {\mathcal{E}}_x    \right) \
       & \ \subset \
        f \left( x \right) + 
          {\frac{\sqrt{ {\left\{\sigma_{\rm{max}}\left( V(x) 
                                     \right)\right\}}^2 - 1 }}
                {\sigma_{\rm{max}}\left( V(x) \right) 
                }
          } 
        \cdot \  {\mathcal{E}}_x   ,
   \end{align*}
   where~$  \sigma_{\rm{max}}\left( V(x) \right) $ denotes the largest eigenvalue of~$ V ( x ) . $
  Thus~$ f $ maps the Schur ellipsoid centered at~$ x  $  to within a shrunken version of the ellipsoid, centered at~$ f(x) . $
\begin{figure}
    \begin{center}
    {\includegraphics[width=0.99\linewidth]{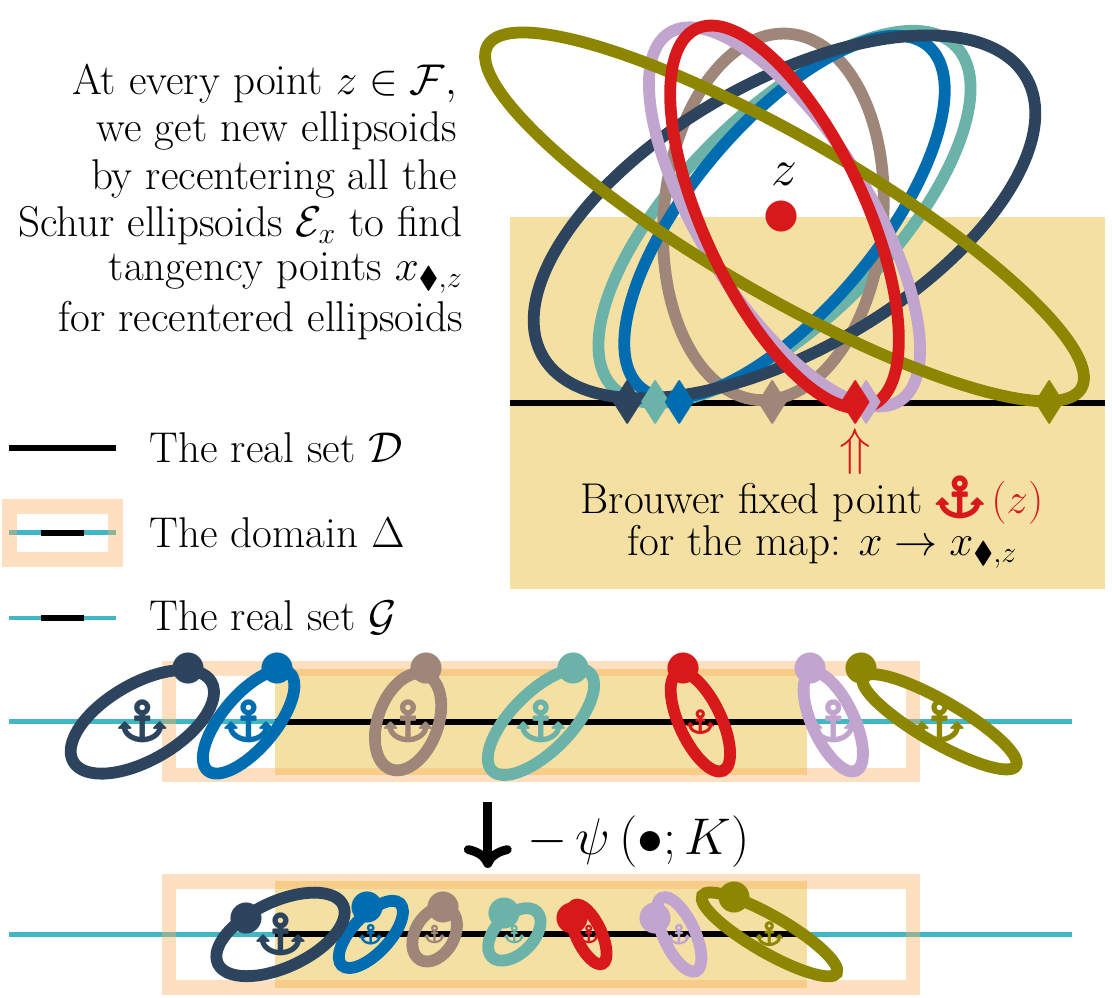}}%
    \end{center}
    \caption{Schur ellipsoids.\label{fig:schurEllipses}}
\end{figure}

    \noindent
    {\textcolor{cyanTeal}{\textit{{\colorbox{tolPaleYellow}{Step 3}}: $ f $ maps even the boundary points of~$ \Delta $ back into~$ \Delta  $}}}
  The domain~$ \Delta $ can be viewed as the set of an infinite collection of sheets that are parallel to the face~$ {\mathcal{F}} . $   Each such sheet is no farther from the set~$  {\mathcal{D}}  $  than~$ {\mathcal{F}}   $ is. 

  If a point~$ z $ on the face~$  {\mathcal{F}}   $ can be included in some Schur ellipsoid centered at some point~$ x $ in the real set~$  {\mathcal{D}}  , $ then we say that~$ z $ is anchored at $ x , $ and we denote that by:
\begin{align*}
    x
    & =
    {\text{\faAnchor}} \left( z  \right) .
\end{align*}
Since the different Schur ellipsoids have different shapes and orientations, it is not clear whether every point~$ z $ on the face~$ {\mathcal{F}}  $ can be anchored at a suitable point in~$  {\mathcal{D}} . $ 

We need a slightly bigger set, to guarantee such anchoring.
Consider the following bounded, convex  set
\begin{align*}
    {\mathcal{G}}  
    & 
        \triangleq 
         \left\{
             x 
             \in {\mathbb{R}}^n :
             \left\lvert x - {\widetilde{\mathcal{D}}} \right\rvert  \le
             \epsilon  
             \right\} .
\end{align*}
If $\epsilon $ is small enough, then the jacobian~$ f_z \left(  x \right) $
is Schur stable, everywhere on~$ {\mathcal{G}} .  $

Now consider any point~$  z $ on the face~$ {\mathcal{F}} .   $ 
For each point~$ x \in  {\mathcal{G}}   , $ construct a corresponding Schur ellipsoid~$ {\mathcal{E}}_x  , $  but  centered at~$ z . $
This ellipsoid has a unique point of tangency with the real space~$  {\mathbb{R}}^n . $ Call this point~$ x_{\blacklozenge , z} .$


The map from a point~$ x $ to the corresponding tangency point~$ x_{\blacklozenge , z} $
is continuous, because~$ V \left( x \right) $ and $ r\left( x \right) $ depend 
continuously on~$ x . $  If~$ \eta $ is small enough, then the Schur
ellipsoids~$ {\mathcal{E}}_x   $ are small, so that  every tangency point~$ x_{\blacklozenge , z} $  falls within the real set~$ {\mathcal{G}} .  $ 

Applying the Brouwer fixed point theorem for the map~$ x \to x_{\blacklozenge , z}  , $  on the bounded, convex set~$  {\mathcal{G}} , $ 
we get at least one fixed point. This fixed point is the anchor for~$ z .  $

Since~$ f $ maps all points in~$  {\mathcal{G}}  $
to within the set~$  {\mathcal{D}} , $ it follows that~$ f $ maps all points
 on~$  {\mathcal{F}}   $ to within~$ \Delta . $
 
 For points on the sheets that are parallel to the face~$ {\mathcal{F}} ,  $ 
 and are closer than~$ {\mathcal{F}}   $ to the real set~$ {\mathcal{D}}  , $
 we can apply the same geometric constructions and arguments as in Steps~2,~3.
 Then we can conclude that~$ f $ maps all points on those sheets to within
 the domain~$ \Delta . $
%
%
%
%
\end{proof}
\begin{theorem}\label{theorem:symmetryAndUnimodality} 
If the requirements of Theorem~\ref{theorem:globallyStableFixedPoint} are met, then the RFS~\eqref{eqn:RFSdynamics} 
has a unique limit cycle that is symmetric, unimodal and globally asymptotically stable.
\end{theorem}
\begin{proof}
    Let~$ {\widehat{x}} \in {\mathcal{D}} $ be the Brouwer fixed point for the map~$ {\mathbf{-}} \, \psi_+\left( \cdot ; K \right)  $ promised by     
  Theorem~\ref{theorem:globallyStableFixedPoint}.  
  The  sequence:
    \begin{gather*}
        \bigl\{
       {\widehat{x}} ,  
       \ {\mathbf{-}} \, \psi_+\left( {\widehat{x}} ; 1 \right)  ,
        \ldots ,
      \  {\mathbf{-}} \, \psi_+\left( {\widehat{x}} ; K - 1 \right) 
        \bigr\}
    \end{gather*}
    is a closed cycle for the map~$ {\mathbf{-}} \, \psi_+\left( \cdot ;  1 \right)  . $  
    Every element in this sequence is a fixed point for the map~$ {\mathbf{-}} \, \psi_+\left( \cdot ;  K \right)  . $
   Therefore, all elements of this sequence have to be the same. Otherwise, we would have a violation Theorem~\ref{theorem:globallyStableFixedPoint}.  

    Hence~$ {\widehat{x}} $ is a fixed point for the 
    map~$ {\mathbf{-}} \, \psi_+\left( \cdot ;  1 \right)  . $    Hence the RFS has a symmetric, unimodal limit cycle. Every trajectory is ultimately confined to the ball~$ {{\mathcal{B}}_{\rm{ultimate}}}  , $  and  all trajectories through it must cross the set~$  {\mathcal{D}} . $ Hence 
the said limit cycle is globally asymptotically stable.
\end{proof}

\section{Smooth approximation of relay element\label{section:smoothApproximation}}
In this section, we do two things: (i)~we study a wider class of plant transfer functions than the BRL-URF class, and (ii)~we do this by analyzing a smooth approximation of the RFS, by approximating the discontinuous relay operator by the~$\tanh$ sigmoidal function. We shall vary the approximation parameter for the sigmoid, and as we do that, we look for structural stability of the entire ODE flow, and  Hopf bifurcations in which  the equilibrium point at the origin
produces a limit cycle.

\subsection{The $\tanh$ approximation of the relay element}
\begin{figure}
\begin{center}
 \begin{tikzpicture}
     \begin{axis}[
             width = \linewidth, height = 0.4\linewidth,
          xmin=-4.5, xmax=4.5,
          ymin=-1.3, ymax=1.3,
          axis lines=center,
          axis on top=true,
          domain=-4.5:4.5,
          ylabel=$u$,
          xlabel=$y$,
          ticks = none,
          y label style={at={(0.45,0.95)}},
      ]
     \addplot [mark=none,draw=redSeries1,ultra thick, samples = 40] {tanh(0.5*\x)};
     \addplot [mark=none,draw=redSeries2,ultra thick, samples = 80] {tanh(\x)};
     \addplot [mark=none,draw=redSeries3,ultra thick, samples = 100] {tanh(2*\x)};
     \addplot [mark=none,draw=redSeries4,ultra thick, samples = 200] {tanh(10*\x)};
\end{axis}
\end{tikzpicture}
\end{center}
     \caption{Ideal relay characteristic asymptotically modelled by the family of smooth functions:~$u = \tanh{\left( \gamma \, y\right)}$.\label{fig:tanh}}
\end{figure}
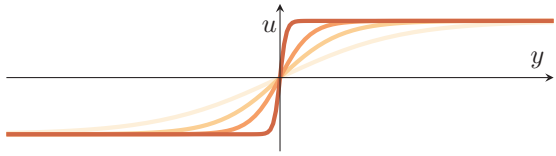

We shall use the smooth approximation:
\begin{align*}
    {\text{sign}}  \bigl(  y \bigr) & 
    \approx \tanh{ \bigl(  \gamma \cdot  y \bigr)} ,
    \ \ y \in {\mathbb{R}} , 
\end{align*}
for large enough values of the parameter~$ \gamma .$
As shown in Figure~\ref{fig:tanh}, the approximation gets better and better for larger and larger choices of the parameter~$ \gamma . $
The closed loop is then the Smooth feedback system~(SFS): 
\begin{align}
    {\frac{d}{dt}} \zeta & =
    A \zeta  - B \, \, \tanh \bigl( \gamma \,C \, \zeta \bigr) .
    \label{eqn:tanhDynamics}
\end{align}
where the companion matrix~$A$ and the vectors~$B , C$ are given by Equation~\eqref{eqn:companionA}.
The origin is an equilibrium point for this closed loop, 
at every value of~$ \gamma . $
Around this equilibrium point, we get the following linearized flow:
\begin{align}
    {\frac{d}{dt}} \delta& =
    \bigl( A   - \gamma B C \bigr) \delta.
    \label{eqn:linearizationAtOrigin}
\end{align}
Because the row vector~$C$ has only zero elements except for its last entry, the $ n~\times~n$~matrix~$ B C$ has non-zero elements only on its last column. 
Therefore the matrix~$A - \gamma B C$ is a companion matrix, whose characteristic polynomial is:
\begin{gather*}
    \sum_{i= 0}^{n-1} {\left( a_{i} + \gamma b_{i} \right) } \lambda^i
 \ \  \ + \ \ \ \lambda^n  .
\end{gather*}

\subsection{Local, global versions of the Hopf bifurcation theorem}
By  {\textit{local bifurcation}} we mean a bifurcation~(change in the qualitative behaviour of orbits) that holds over an infinitesimally small interval of the bifurcation parameter. By a {\textit{global bifurcation}} we 
mean one that holds over some nontrivial finite interval, or an infinite interval.
\begin{theorem}[Local version of the Hopf bifurcation theorem in terms of the Describing function~\cite{allwright1977harmonicBalanceAndHopfBifurcation,mees1981dynamicsOfFeedbackSystems}]
\label{theorem:describingFunctionHopf}
    Consider the nonlinear feedback system in Lure form:
\begin{align}
    {\frac{d}{dt}} \chi & =
    A \chi  - B \, \cdot \, Y \bigl(  C \chi , \gamma \bigr) , \quad \text{for} \
    \chi \in {\mathbb{R}}^n , 
    \label{eqn:lureFeedbackSystem}
\end{align}
where~$ \gamma  $ is a real-valued parameter. 
Assume that the scalar-valued function~$ Y $ is four times continuously differentiable w.r.t. to the state~$ \chi , $ and continuously differentiable w.r.t. the parameter~$ \gamma . $ 
    Assume that the origin is an equilibrium point for this feedback system, for all the considered values of~$ \gamma .$
   Let~$ G \left( s \right) $ denote the transfer function of the linear plant in the feedback system, meaning that~$ 
     G \left(  s \right)  = C { \left(  s I - A \right) }^{-1} B . 
    $
    Let~$ {\widetilde{G}}\left( s \right) $ denote the open loop transfer
function of the combination  of plant and {\textit{linearization}} of the nonlinear element, meaning that
    \begin{align*}
        {\widetilde{G}} \left(  s \right)  & = 
        {{G}} \left(  s \right) \, \cdot \,
        {\left.
        {\frac{\partial}{\partial {\left\{C \chi\right\}}}} Y \left( C \chi   , \gamma \right) \right\rvert}_{ C \chi = 0 }  .
    \end{align*}
Suppose that as~$ \gamma $ passes the value~$ \gamma_0 ,  $ 
the Nyquist locus of~$ 
{\widetilde{G}} \left(  j \omega \right)  
    $
    passes the Nyquist critical point~$ {\mathbf{ - 1 }} , $ at a unique positive  frequency~$ \omega_0 . $ Suppose that the partial derivatives
    \begin{gather*}
        {\left.
        {\frac{\partial}{\partial {\gamma }}} 
        {\widetilde{G}} \left(  j \omega \right)  
        \right\rvert}_{
            \substack{ \gamma = \gamma_0 \\ \omega = \omega_0  }
        }
        ,
        \quad 
        {\left.
        {\frac{\partial}{\partial {\omega }}} 
        {\widetilde{G}} \left(  j \omega \right)  
        \right\rvert}_{
            \substack{ \gamma = \gamma_0 \\ \omega = \omega_0  }
        }
    \end{gather*}
    are nonzero, and as complex numbers are not parallel. This would mean that the Nyquist locus of~$ {\widetilde{G}} \left(  j \omega \right)   $
  has a generic  and non-pathological crossing of the Nyquist critical point. 
\newline
    Let~$ L \left( \theta , \omega \right)  $ 
    denote the 
Describing function of the nonlinear element~$ Y , $ for harmonics up to the second one. 
    Suppose that for~$ \gamma = \gamma_0 , $ and as we vary~$ \theta $ the locus of the Describing function~$ L \left( \theta , \omega = \omega_0 \right)  $ crosses the Nyquist locus of~$ {\widetilde{G}} \left( j \omega \right)  $ at the Nyquist critical point. If this crossing is transversal, then a local Hopf bifurcation is guaranteed. In specific,  
there is  a small interval of parameter values that lie on one side of~$ \gamma_0 $, such that a limit cycle exists around the equilibrium point.
    Close to the origin, 
this limit cycle is unique.
\newline
Suppose that the equilibrium point at the origin is stable, before  the Hopf bifurcation. If the locus of~$ L \left( \theta , \omega 
    = \omega_0 \right)  $  points in a direction outward from the Nyquist locus of~$  {\widetilde{G}} \left( j \omega \right)  , $  then the limit cycle is 
    stable.
\end{theorem}
We now complement the above local version, with  the following global version.
The  version below requires the same smoothness on the nonlinear element, as the local version above required. The global version predicts the effect of increasing the bifurcation parameter beyond the critical value.  It   permits only four possible qualitative trends in the 
 subsequent  evolution of  limit cycle. 
\begin{theorem}[Global Hopf bifurcation theorem\cite{alexanderYorke1978globalBifurcationsOfPeriodicOrbits}\label{theorem:globalHopf}] 
Consider the nonlinear feedback system~\eqref{eqn:lureFeedbackSystem}. Assume  that  the nonlinear element in its RHS satisfies the continuous differentiability properties demanded by Theorem~\ref{theorem:describingFunctionHopf}. 
    Suppose that when~$ \gamma =  \gamma_0 , $ the linearized dynamics at the origin:
\begin{align*}
    {\frac{d}{dt}} \delta & =
    \left( A   -  B C \, 
      \cdot
      {\left.
     {\frac{\partial}{\partial {\left\{C \chi\right\}}}} Y \left( C \chi   , \gamma \right) \right\rvert}_{ C \chi = 0 }  
    \,\right) \delta .
\end{align*}
    has a pair of complex conjugate eigenvalues~$ \pm j \omega_0 . $ 
    Assume that this eigenvalue has an odd multiplicity.
\newline
    Assume also that for $ \gamma $
   near but not equalling $ \gamma_0 , $ the linearization has no eigenvalue with real part zero. Then as as $ \gamma $ is increased from the value~$ \gamma_0 , $
   there is at least one periodic orbit that persists in the following sense. 
   The periodic orbit evolves with some non-contradictory combination of 
    the following trends:
    (T1)~a periodic orbit orbit exists for every value of~$ \gamma $ that is larger than~$ \gamma_0 , $  such that the sizes and time periods  are both uniformly bounded, and/or
    (T2)~a reverse Hopf bifurcation occurs subsequently, and the family of periodic orbits  terminates at an equilibrium point, at some value of~$ \gamma $ that is bigger than~$ \gamma_0  , $
    and/or
    (T3)~as the value of~$ \gamma $ increases, so does the period of the orbit, 
    and/or
    (T4)~as the value of~$ \gamma $ increases, so does the size of the orbit.
\end{theorem}
We of course care about the first of the above four possible evolution trends for the periodic orbit. Unfortunately, we have found no arguments that could guarantee this trend, and rule out the others.

We shall now apply the above theorems to our SFS~\eqref{eqn:tanhDynamics}.


\subsection{Explanatory power of the plant's root locus diagram}
\begin{figure}
\begin{center}
    \includegraphics[width=\linewidth]{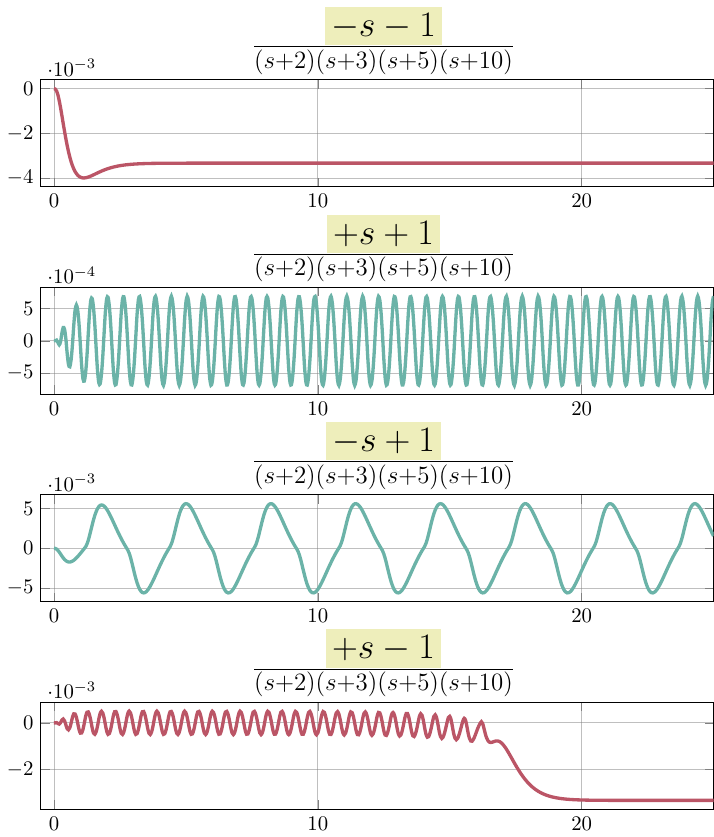}
\end{center}
\caption{The sign of plant zero and of the DC gain affect the closed loop behaviour of the RFS.\label{fig:effectOfzeroesignAndDCgainSign}}
\end{figure}
\begin{figure}
\centering
%
    \subcaptionbox{1st plant in Fig.~\ref{fig:effectOfzeroesignAndDCgainSign}. Pitchfork bifurcation.%
\label{fig:rootLocusp1}}
    {\includegraphics[width=0.49\linewidth]{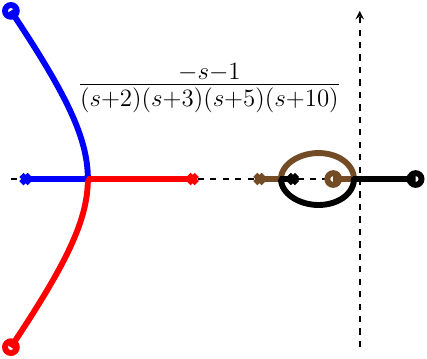}}%
\hfill
    \subcaptionbox{2nd plant in Fig.~\ref{fig:effectOfzeroesignAndDCgainSign}. A supercrtical Hopf bifurcation.%
\label{fig:rootLocusp2}}
    {\includegraphics[width=0.49\linewidth]{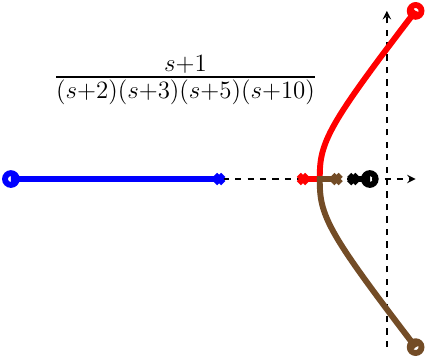}}%
    \\    
    \subcaptionbox{3rd plant in Fig.~\ref{fig:effectOfzeroesignAndDCgainSign}. A supercritical Hopf bifurcation.%
\label{fig:rootLocusp3}}
    {\includegraphics[width=0.49\linewidth]{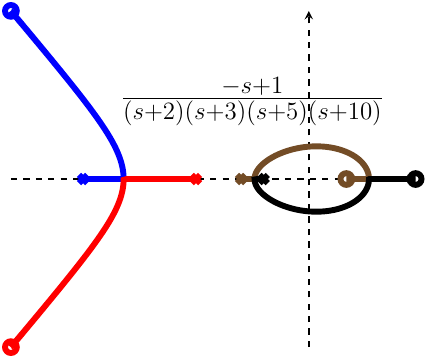}}%
\hfill
    \subcaptionbox{4th plant in Fig.~\ref{fig:effectOfzeroesignAndDCgainSign}%
    \label{fig:rootLocusp4}. Pitchfork bifurcation, followed by a subcritical Hopf bifurcation.}
    {\includegraphics[width=0.49\linewidth]{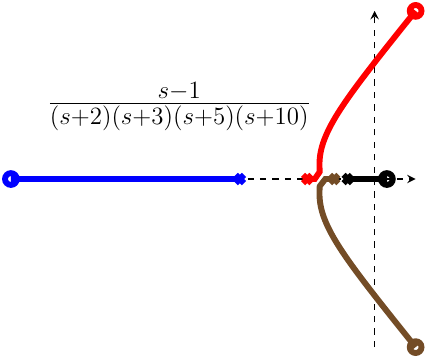}}%
 \caption{Explanation for behaviours shown in Figure~\ref{fig:effectOfzeroesignAndDCgainSign}, via the bifurcations predicted by the root loci of the plants.%
\label{fig:exampleRootLocusPlotsForFigure2}}
\end{figure}
The root locus of the plant shall track the bifurcations undergone by the equilibrium point that the SFS has at the origin.
With the examples below, we illustrate the predictive power of the root locus.

Figure~\ref{fig:effectOfzeroesignAndDCgainSign} shows RFS outputs, for trajectories that start close to the origin. All four plants are stable, and have the same set of poles. But they differ in the sign of the plant zero, and hence in the sign of the DC gain.
Explanations of their RFS behaviours are offered by Figure~\ref{fig:exampleRootLocusPlotsForFigure2}.

    In the first plant,  
   a Pitchfork bifurcation produces two stable equilibrium points, and renders the origin unstable. %
   In the second plant, 
    a supercritical Hopf bifurcation produces a stable limit cycle, having a high frequency and a relatively small amplitude. %
    In the third plant, 
    a supercritical Hopf bifurcation produces a stable limit cycle of relatively small frequency. %
 And in the fourth plant, 
  a Pitchfork bifurcation produces two stable equilibrium points, and renders the origin unstable. A subsequent Hopf bifurcation is subcritical. 

\begin{figure}
\centering
%
    \subcaptionbox{Example~1\phantom{1} in~Johansson~et~al.~\cite{johanssonRantzerAstrom1999fastSwitchesInRFS}
\label{fig:rootLocusPlant1JohanssonEtAl}}
    {\includegraphics[width=0.44\linewidth]{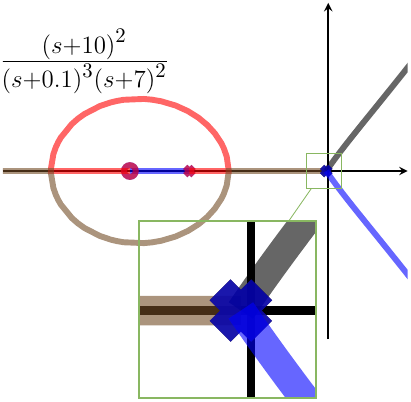}}%
%
    \hfill
    \subcaptionbox{Minimum phase, relative degree~2%
\label{fig:rootLocusFourthOrderRelDeg2}}
    {\includegraphics[width=0.49\linewidth]{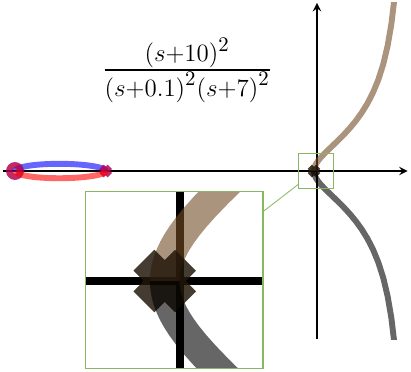}}%
%
    \\    
    \subcaptionbox{Minimum phase, relative degree~2
\label{fig:rootLocusPlantThirdOrderRelDeg2}
    }
    {\includegraphics[width=0.49\linewidth]{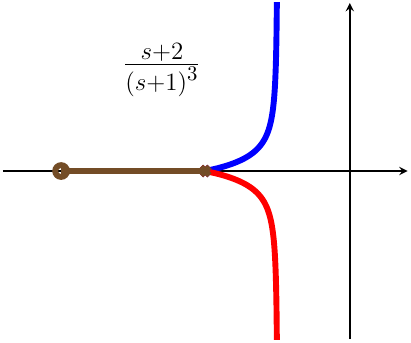}}%
%
    \hfill
    \subcaptionbox{Minimum phase, with relative degree~1%
\label{fig:rootLocusThirdOrderMinimumPhase}}
    {\includegraphics[width=0.49\linewidth]{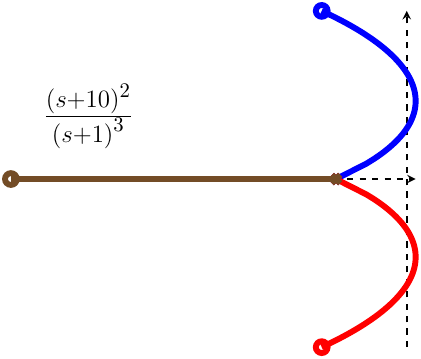}}%
%
    \\  
    \vskip 1mm
    \subcaptionbox{BRL-URF plant%
\label{fig:rootLocusBRL_URF}}
    {\includegraphics[width=0.49\linewidth]{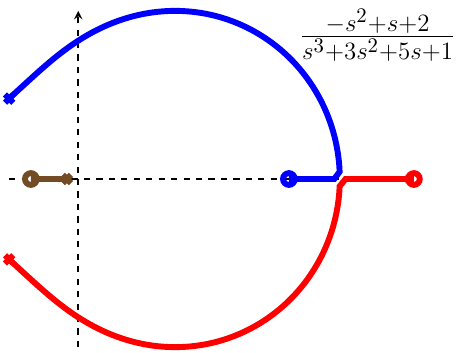}}%
    \hfill
    \subcaptionbox{Minimum phase plant, with relative degree~3%
\label{fig:rootLocusSixthOrderMinimumPhase}}
    {\includegraphics[width=0.49\linewidth]{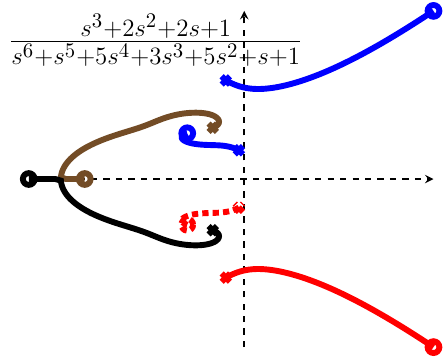}}%
%
    \caption{Example root locus plots for a variety of plants.\label{fig:exampleRootLocusPlotsExtra}}
\end{figure}
%
%
%
%
%
%
\subsection{RHP zeros, and positive DC gain imply Hopf bifurcations}
The SFS shall have Hopf bifurcations, whenever the plant root locus has
generic crossings of the imaginary axis, at complex conjugate values. 
The latter happens, if the plant has a positive DC gain, and RHP zeros.

Two types of zeros, can be associated with any
transfer function, as is well known.  
The transmission zeros come from those finite points in the complex plane where the transfer function vanishes. The `zeros at infinity' come from those 
close loop poles that grow unboundedly, as the scalar linear feedback gain is
 increased to infinity. 

 Suppose all poles of the open loop plant are stable. If any of either kind of a plant's zeros are in the RHP, then conditions are ripe for Hopf bifurcations. Every branch of the plant root locus has to start from the LHP. And some branches have to cross over to the RHP; as many branches as there are RHP zeros.

 Suppose the plant's DC gain is positive. Then 
zero cannot be a closed loop pole, for any
 value of the negative feedback gain.
 Hence every cross-over by a root-locus branch must cross
 the imaginary axis as complex conjugate values.
These conditions are sufficient for a local Hopf bifurcation,
because of Theorem~\ref{theorem:describingFunctionHopf}.

The global version due to Theorem~\ref{theorem:globalHopf} also comes into force. We can calculate the Describing function of up to the second harmonic, (see Page~175,~\cite{mees1981dynamicsOfFeedbackSystems}) to be:
 \begin{align*}
     L \left( \theta , \omega \right)
     & =
     - 1 + {\theta}^2 \cdot 
            {\frac{\gamma^2}{4}}
     \cdot  G \left( j \omega \right)   . 
 \end{align*}
 At the critical frequency~$ \omega_0  , $ the locus of~$ L_1 $ is parallel to the negative real semi-axis. Hence, this locus shall be transversal to, and pointing outwards of the Nyquist locus of~$ {\widetilde{G}} \left( j \omega \right)   , $ unless the latter passes by the point~$ - 1 $ tangentially.
 Hence 
 generically we can expect that for a stable plant with RHP zeros: (i)~as many Hopf bifurcations happen for the SFS, as there are plant  root locus crossings of the imaginary axis, and (ii)~the first such crossing must give a supercritical Hopf bifurcation.

For a RFS to have a periodic orbit, it is not sufficient that some
Describing function analysis predicts it. This can be seen from
the example of the root locus  diagram in 
Figure~\ref{fig:rootLocusThirdOrderMinimumPhase}.
 From this root locus, and the corresponding Nyquist diagram,
 we see that one periodic orbits is predicted by Describing function analysis. 
 In specific there are two separate transversal intersections of the plant's Nyquist locus, with the locus of the  Describing function that includes up to the second harmonic.
The SFS for this plant follows the bifurcation trend~(T2), and has a Hopf bifurcation, followed by a reverse Hopf bifurcation.
\subsection{SFS for which all trajectories decay towards the origin}
For some  minimum phase plants having relative degree one or two, we  show below that the SFS phase portrait has the same qualitative features, regardless of whether the amplifier gain~$\gamma$ is small or large. 

The notion of structural stability concerns the preservation of qualitative features of the behaviour of a dynamical system, under
small  changes in the governing dynamics.

\begin{definition}[Structural stability]
    Consider a dynamical system governed by an ODE. 
Consider the set of all ODEs obtained by adding to the right hand side of the ODE, perturbations that are both small and continuously differentiable. 
    If for every sufficiently small perturbation, trajectories of the
    perturbed ODE 
    are topologically equivalent~(homeomorphic) to those  of the unperturbed ODE, then 
    the dynamical system is said to be structurally stable.
\end{definition}

For any ODE, a  {\textit{hyperbolic equilibrium point}} is one  whose linearization has no eigenvalues with real part zero. Locally around
any such point, the ODE is structurally stable. This follows from the following two ideas: (i)~locally around the hyperbolic equilibrium point, the trajectories of the ODE are homeomorphic to those of the linearization, as per the Hartman-Grobman theorem, and (ii)~any small change to the right hand side of the ODE induces continuous changes in the linearization, which in turn induces continuous changes in the eigenvalues.

Similarly, if the jacobian of the right hand side of an ODE is hyperbolic everywhere in a closed and bounded state space, then the ODE is structurally stable. This idea is captured by the notion of Anosov diffeomorphisms and Anosov's guarantee of their structural stability.

\begin{definition}[Anosov diffeomorphism, for Euclidean spaces]
 A differential equation of the form:
    \begin{align*}
        {\frac{dx}{dt}}  & =
        g \left( x \right) , \ \ \text{for} \ x \in X \subset {\mathbb{R}}^n ,
    \end{align*}
    is called an Anosov diffeomorphism over the set~$X$, if~$X$ is positively invariant under the flow of the ODE, and at every point~$x$ from the set~$X,$ the jacobian matrix
    \begin{gather*}
        {\frac{dg{\left(x\right)}}{dx}}  
    \end{gather*}
    has no eigenvalues with real part zero.
\end{definition}
\begin{theorem}[Theorem~7.9 of~\cite{irwin2001smoothDynamicalSystems}]%
{Anosov diffeomorphisms are structurally stable.
}
\end{theorem}
We now apply these notions to the SFS, for 
those minimum plants for which no bifurcations can happen, as the parameter~$ \gamma $ ranges from~$ 0  $ to $ + \infty .$
\begin{lemma}
{If the plant root locus neither crosses nor touches the imaginary axis, then for every nonegative value of
~$ \gamma , $ the SFS of Equation~\eqref{eqn:tanhDynamics} is structurally stable.
}
\end{lemma}
\begin{proof}
    We shall show that the jacobian matrix of the RHS of the SFS ODE is hyperbolic, for every nonegative value of the parameter~$ \gamma , $ at every point in the state space.

    At any given point~$ x ,$ the jacobian of the RHS is:  
    \begin{gather*}
        A 
        - {\frac{ \gamma }{ {\left\{ \cosh{ \left( \gamma \cdot C x\right)} \right\}}^2 }} 
        B C .
    \end{gather*}
    We compare this jacobian with that at the origin~(Equation~\eqref{eqn:linearizationAtOrigin}), by comparing the functional behaviour of
    the respective  coefficients in front of~$ B C .$ 
    For~$ \gamma \ge 0 , $  we have:
    \begin{align}
      {\frac{ \gamma }{ {\left\{ \cosh{ \left( \gamma \cdot C x\right)} \right\}}^2 }} 
        &  \le \gamma, \nonumber \\
      {\frac{ \gamma }{ {\left\{ \cosh{ \left( \gamma \cdot C x\right)} \right\}}^2 }} 
        & \le {\frac{0.4478}{\left\lvert C x \right\rvert}} , \quad \text{and,} \nonumber \\
        \lim_{\gamma\to + \infty}  {\frac{ \gamma }{ {\left\{ \cosh{ \left( \gamma \cdot C x\right)} \right\}}^2 }} 
        & =  
        \begin{cases}
            + \infty, & \text{if} \ C x = 0, \ \text{and,} \\
            0, & \text{if} \ C x \neq 0.
        \end{cases}
        \label{eqn:limitOfGammaLikeCoeff}
    \end{align}
    By assumption, as~$ \gamma $ ranges from~$ 0  $ to $ + \infty , $ the matrix~$  A - \gamma B C $ remains Hurwitz. 
    Hence as~$ 
{{ \gamma } / { {\left\{ \cosh{ \left( \gamma \cdot C x\right)} \right\}}^2 }} 
    $ ranges from~$ 0  $ to $ 0.4478 / \left\lvert C x \right\rvert , $ the matrix~$ 
        A 
        - {\left(  { \gamma } / { {\left\{ \cosh{ \left( \gamma \cdot C x\right)} \right\}}^2 } \right)} 
        B C $ remains Hurwitz. Hence the jacobian is hyperbolic everywhere, for nonegative values of~$ \gamma . $

        A sufficiently large Lyapunov sublevel set for the linear ODE:~$
    {\dot{\chi}} = A \chi, 
        $ is positively invariant for the SFS~\eqref{eqn:tanhDynamics}, because the extra term involving the~$ \tanh $ function is bounded. Applying Anosov's theorem on this positively invariant set yields the structural stability. 
\end{proof}

The  asymptotic behaviour the SFS suggests but does not confirm the asymptotic behaviour of the RFS.

For plants that meet the requirements of the previous lemma,
regardless of the nonegative value of~$ \gamma ,  $ 
the SFS trajectories decays asymptotically to the origin.
If we believe that the behaviour of the RFS is practically indistinguishable from that of the SFS for very high values of~$ \gamma , $ then we can also believe that the RFS trajectories move asymptotically in the direction of the origin. If the plant's relative degree is 2, then this suggests that the trajectories decay to the origin. But if the plant's relative degree is 1, then the SFS can only predict the smoothed result of solutions averaged over the sliding portions of the
RFS trajectories.
Figure~\ref{fig:rootLocusPlantThirdOrderRelDeg2} gives an example, where the SFS  and RFS behaviours seem to coincide empirically.
\subsection{Liouville's theorem and local volume shrinkage\label{section:liouville}}

For a broad class of periodic orbits that can arise in the SFS~\eqref{eqn:tanhDynamics}, we shall now show local volume shrinkage around those orbits, as the gain parameter~$ \gamma $ tends to infinity.

Let the SFS~\eqref{eqn:tanhDynamics} possess a periodic orbit of finite period: 
 \begin{gather}  
 \left\{  
    p\left( t ; \gamma \right) : 
     0 \le t \le T_{\rm{period}}
 \right\} .
    \label{eqn:periodicOrbitForSFS}
 \end{gather}  
Then the linearized flow around this orbit is given by:
\begin{align}
    {\frac{d}{dt}} \delta & =
    \left[ A   -
    {\frac { \gamma } { {\left\{ \cosh{ \left( \gamma \cdot C p\left( t ; \gamma \right) \right)} \right\}}^2 } } 
    B C \right] \delta .
    \label{eqn:linearizationAroundPeriodicOrbit}
\end{align}
Let~$ \Phi \left( t , 0 \right) $ denote the fundamental matrix for the above time-periodic linear system. 
Differentiating both sides of~\eqref{eqn:tanhDynamics} w.r.t. time, and using
the periodicity: 
\begin{align*}
    {\frac{d}{dt}}p \left( t + T_{\rm{period}} \right) 
    & = 
    {\frac{d}{dt}}p\left( t ; \gamma \right) ,
    \quad \text{for all} \ \
     t \ge 0, 
\end{align*}
we get~(Equation~2.6, Chapter~13 of \cite{coddingtonLevinson1955theoryOfODEs}):
\begin{align*}
    \Phi \left( T_{\rm{period}} , 0 \right)  
    {\frac{d}{dt}} p\left( t ; \gamma \right)
    & =
    {\frac{d}{dt}} p\left( t ; \gamma \right) .
\end{align*}
This means that~$ 1 $ is an eigenvalue for~$ 
    \Phi \left( T_{\rm{period}} , 0 \right)  . $

The product of all eigenvalues can be calculated, by following
Willems~\cite{jacquesWillems1968stabilityOfNonlinearOscillations}, to give us the following theorem.
\begin{theorem}[local volume shrinkage]%
    \label{theorem:localVolumeShrinkage}
Let the periodic orbit~\eqref{eqn:periodicOrbitForSFS}
have a finite period. Let the orbit never merely graze without crossing, the hyperplane:~$
    \left\{ 
     x \in {\mathbb{R}}^n : C x = 0
    \right\} .
    $
    Let the possible crossings be at isolated time instants, transversal, and 
    finite in number. 
In other words, if~$ Z_p $ denotes the zero set of the 
signal~$ C p\left( t ; \gamma \right) , $ then this set must be finite,
and transversality means that 
    \begin{gather*}
        C {\frac{d}{dt}} p \left(  t_0 ; \gamma \right)   \ \neq \ 0 ,
        \quad \text{if} \ \ t_0 \in Z_p 
    \end{gather*}
    for very large but finite~$ \gamma , $ and also in the limit when $ \gamma $ is taken to infinity. Assume that for every time~$ t_0 \in Z_p , $ as $ \gamma $~goes to infinity, the signal~$ C {\tfrac{d}{dt}} p $  has finite left and right limits of the same sign, at all zero crossing times~$ t_0 . $

    Then 
    as the parameter~$\gamma$ goes to infinity, 
    the determinant 
    of the monodromy matrix~$   \Phi \left( T_{\rm{period}} , 0 \right)  $
 goes to the limit: 
    \begin{gather*} 
        \exp\left(%
        { - \; a_{n-1} \, T_{\rm{period}}  \; - \; b_{n-1} 
        \sum_{t_0\in Z_p} 
            { \frac
                 { 
                 \left\lvert  
                    \rho^{\mathbf\huge{-}} \left( t_0 \right) 
                 \right\rvert  
                 +
                 \left\lvert  
                    \rho^{\mathbf\huge{+}} \left( t_0 \right) 
                 \right\rvert  
                 }
                {
             \left\lvert  
                    \rho^{\mathbf\huge{-}} \left( t_0 \right) 
                    \rho^{\mathbf\huge{+}} \left( t_0 \right) 
             \right\rvert
             } 
             }
        }
        \right) ,
    \end{gather*} 
    where the left and right limits~$
    \rho^{\mathbf\huge{-}} \left( t_0 \right) , 
     \ 
   \rho^{\mathbf\huge{+}} \left( t_0 \right) ,    
     $ 
     are: 
     \begin{align*} 
    \rho^{\mathbf\huge{-}} \left( t_0 \right) 
         &  \triangleq 
    \lim_{
         \gamma \to + \infty 
     }
    \lim_{
         \nu \to 0 
     }
     {
                 C {\frac{d}{dt}} p \left( t_0 - \left\lvert \nu \right\rvert ; \gamma \right) 
    } 
         ,  \ \text{and} \\
    \rho^{\mathbf\huge{+}} \left( t_0 \right) 
         &  \triangleq
    \lim_{
         \gamma \to + \infty 
     }
    \lim_{
         \nu \to 0 
     }
     {
                 C {\frac{d}{dt}} p \left( t_0 + \left\lvert \nu \right\rvert ; \gamma \right) 
    } .
     \end{align*}  
\end{theorem}
\begin{proof}
By Liouville's theorem, the 
    determinant of the state transition matrix, 
 namely~$
{\rm{det}} \left(  \Phi \left( t , 0 \right)  \right) , 
$ equals
\begin{gather*}
         \exp\left(%
         {%
             \int_0^t 
         {\rm{trace}}
\left(  A - 
    {\frac { \gamma } { {\left\{ \cosh{ \left( \gamma \cdot C p\left( {\widehat{t}} ; \gamma \right) \right)} \right\}}^2 } } 
    B C 
        \right)  d{\widehat{t}}
        } 
       \right) .
\end{gather*}
Because the two matrices~$ A, C $  have special structure~\eqref{eqn:companionA}, 
    \begin{multline}
      \int_0^t   {\rm{trace}}
\left(  A - 
    {\frac { \gamma } { {\left\{ \cosh{ \left( \gamma \cdot C p\left(  {\widehat{t}}  \right) \right)} \right\}}^2 } } 
    B C 
        \right) d{\widehat{t}}
        \\ = 
        - a_{n-1} \, t
        - b_{n-1}
         \int_0^t
    {\frac { \gamma } { {\left\{ \cosh{ \left( \gamma \cdot C p\left( {\widehat{t}}  \right) \right)} \right\}}^2 } } 
         d{\widehat{t}} .
         \label{eqn:expressionForTrace}
    \end{multline}
    In the RHS above, the integral multiplying~$ b_{n-1} $ can be 
    broken into the following two types of  integrals. 
The first type covers short durations around fast changes of the sigmoidal relay element - these integrals are over 
    infinitesimally small intervals 
    of time, around those 
    instants when the scalar signal~$ C p\left( t ; \gamma \right) $ crosses zero.
    And the second type of integrals covers longer durations when there are slight or no changes of the sigmoidal relay element.
    Thus the integral appearing on the RHS of~\eqref{eqn:expressionForTrace} equals: 
    \begin{gather*}
        \sum_{
            \left\{  t_0 \in Z_p
            \right\}
        } 
        \int_{t_0 - \epsilon}^{t_0 + \epsilon}
    {\frac { \gamma } { {\left\{ \cosh{ \left( \gamma \cdot C p\left( {\widehat{t}}  \right) \right)} \right\}}^2 } } 
         d{\widehat{t}} 
       \\
        +
        \int_{\rm{rest~of~the~horizon}}
    {\frac { \gamma } { {\left\{ \cosh{ \left( \gamma \cdot C p\left( {\widehat{t}}  \right) \right)} \right\}}^2 } } 
         d{\widehat{t}} .
    \end{gather*}
As~$ \gamma $ goes to infinity,  the second integral above goes to zero, in light of~\eqref{eqn:limitOfGammaLikeCoeff}. In the remaining sum of integrals, each term can be dealt as follows. For finite~$ \gamma $ we have:
\begin{multline*}
    \mu_{t_0} \left( \epsilon ; \gamma \right) 
    \ \triangleq \
        \int_{t_0 - \epsilon}^{t_0 + \epsilon}
    {\frac { \gamma } { {\left\{ \cosh{ \left( \gamma \cdot C p\left( {\widehat{t}} ; \gamma \right) \right)} \right\}}^2 } } 
         d{\widehat{t}} , 
         \\
           =
        \int_{t_0 - \epsilon}^{t_0 + \epsilon}
        \left( 
        {\frac { d }
               { dt}
        }
{\tanh{\left(\gamma \cdot C p \left( {\widehat{t}} ; \gamma \right) \right)}}  
        \right)
        {\frac { 1 }
               {  C {\frac{d}{dt}}p\left( {\widehat{t}}  \right)  } 
        }
         d{\widehat{t}} ,
         \\
           =
        \int^{t_0}_{t_0 - \epsilon}
        \left( 
        {\frac { d }
               { dt}
        }
{\tanh{\left(\gamma \cdot C p \left( {\widehat{t}} ; \gamma \right) \right)}}  
        \right)
        {\frac { 1 }
               {  C {\frac{d}{dt}}p\left( {\widehat{t}}  \right)  } 
        }
         d{\widehat{t}} , \\
           + 
        \int^{t_0 + \epsilon}_{t_0}
        \left( 
        {\frac { d }
               { dt}
        }
{\tanh{\left(\gamma \cdot C p \left( {\widehat{t}} ; \gamma \right) \right)}}  
        \right)
        {\frac { 1 }
               {  C {\frac{d}{dt}}p\left( {\widehat{t}}  \right)  } 
        }
         d{\widehat{t}} .
    \end{multline*}
  We have split the interval of integration into two, because as ~$ \gamma $ is taken to infinity, in the limit, the signal $ C {\frac{d}{dt}} p\left( t \right) $ is not an analytic function of~$ t , $ at~$ t = t_0 . $
  But, even in the limit of~$ \gamma $ tending to infinity, this signal is analytic 
  for~$ t $ on the semi-open intervals:~$
\left[ t_0 - \epsilon , t_0 \right) , $
and~$
\left( t_0 , t_0 + \epsilon  \right] .
  $

  As~$ \gamma $ is taken to infinity, the above RHS tends to:
\begin{align*}
    \lim_{\epsilon\to 0} \lim_{\gamma\to\infty} \mu_{t_0} \left( \epsilon ; \gamma \right)  
    & = 
    {\frac{1}{
        \left\lvert \rho^{\mathbf\huge{-}}\left( t_0  \right) \right\rvert 
    }}
  +
    {\frac{1}{
       \left\lvert   \rho^{\mathbf\huge{+}}\left( t_0  \right) \right\rvert 
    }} .
\end{align*}
Setting~$ t  = T_{\rm{period}} $ in~\eqref{eqn:expressionForTrace}
yields the result sought.
\end{proof}
With this theorem, we cannot conclude outright, the stability of a periodic orbit, except in the case of second order plants. However, we can conclude instability of a periodic orbit, if the determinant of the monodromy matrix is bigger than one. For stable plants, the determinant can exceed one, only if the coefficient~$ b_{n-1} $ is negative. 

For BRL-URF plants, the coefficient~$ b_{n - 1} $ is negative. 
Consider a Hopf bifurcation for the SFS with a BRL-URF plant.
 Suppose that 
the critical value of the bifurcation parameter~$ \gamma $ is relatively large. Then for values of~$ \gamma $ that are only slightly past the critical value, the Hopf limit cycle is still newly born.
  Hence its amplitude  is bound to be infinitesimally small, at these values of~$ \gamma  . $ Therefore the corresponding output speeds~$ 
\rho^{\mathbf\huge{+}}\left( t_0  \right) 
, 
\rho^{\mathbf\huge{-}}\left( t_0  \right) 
  $, being proportional to
  the amplitude, are also bound to be small. On the other hand the period is well approximated by~$ 2 \pi / \omega_0 , $ since~$ \gamma $ is 
  close to the critical value.
Hence, the determinant of the monodromy matrix exceeds one, and the corresponding Hopf limit cycle is unstable at birth.

\begin{corollary}
Let the plant transfer function be proper.
Let the sum of its poles be negative,
    so that the coefficient~$ a_{n-1} $ is positive. 
Let the first nonzero Markov parameter be positive,
    so that the coefficient~$ b_{n-1} $ is either positive, or zero. 
Now suppose that the SFS~\eqref{eqn:tanhDynamics}  possess a 
 periodic orbit~\eqref{eqn:periodicOrbitForSFS}. Suppose also that this periodic orbit 
has a finite period, has a finite number of zero crossings, and that such crossings are transversal.
%
Then, as the gain~$ \gamma $ is taken to infinity, around this periodic orbit, the linearized flow shrinks volume in phase space.
\end{corollary}
\subsection{The monodromy matrix from the orbit's parameters\label{section:monodromyMatrix}}
Continuing the previous subsection, we illustrate below
how to obtain the monodromy matrix~$ 
 \Phi \left( T_{\rm{period}} , 0 \right)  . 
$
We shall arrive at this, from knowing  the SFS output's  speeds at zero crossings, and the orbit's time period. If the plant has a relative degree of two or more, then the SFS output is well approximated by a pure sinusoid, and we shall only need to know the amplitude and  time period. 
%
%
We shall
calculate the monodromy matrix, by integrating
\begin{align}
    {\frac{d}{dt}} 
      \left[ e^{- A t} \Phi \left( t , 0 \right)  \right]
      & =
    {\frac { {\mathbf{-}} \, \gamma } { {\left\{ \cosh{ \left( \gamma \cdot C p\left( t ; \gamma \right) \right)} \right\}}^2 } } 
    B C 
      \left[ e^{- A t} \Phi \left( t , 0 \right)  \right]
    ,
    \label{eqn:fundamentalMatrixODE}
\end{align}
separately 
over shorter durations when it changes very fast, and 
over longer durations when the sigmoidal relay output changes very little,  like in the proof of Theorem~\ref{theorem:localVolumeShrinkage}. 

Over the second sort of interval, as~$  \gamma  $ tends to infinity, the RHS of~\eqref{eqn:fundamentalMatrixODE} tends to zero, and
\begin{align*}
      \Phi \left( t_2 , t_1 \right)  
      & \rightarrow
      e^{A \left( t_2 - t_1 \right) } .
\end{align*}

The first sort of interval is infinitesimally small, and so the factor~$
 e^{-A t}
$ in the ODE~\eqref{eqn:fundamentalMatrixODE} has no effect over the short duration. 
However there is an effect due to the  matrix~$ B C .$ This effect describes the change between the tangents to the trajectory, before and after switching.

In the time-varying linear ODE~\eqref{eqn:fundamentalMatrixODE}, the    coefficient matrix at any time commutes with its value at any other time. Hence around any switching time~$ t_0  $, we get:
\begin{multline*}
      \Phi \left( t_0 + \epsilon , t_0 - \epsilon \right)  
       = \\
      \exp\left(
      { - B C 
          \int_{ t_0 - \epsilon }^{ t_0 + \epsilon } {
                 {\frac { \gamma } 
                      { {\left\{ \cosh{ \left(
                             \gamma \cdot C p\left( {\widehat{t}} ; \gamma \right) 
                              \right)} 
                        \right\}}^2 } 
                 } 
              d{\widehat{t}} } }
          \right) , \ \ \text{and} 
\end{multline*}
since~$ {\left( B C \right)}^k = b_{n-1}^{ k - 1 } B C , $ for~$ k \ge 1 ,$  the above RHS equals
\begin{gather*}
      I -
      B C 
       \sum_{k = 1}^{\infty}{
           {\frac{ {\left( - b_{n-1} \right)}^{k-1} }{k!}}
       {\left(
          \int_{ t_0 - \epsilon }^{ t_0 + \epsilon } {
                 {\frac { \gamma } 
                      { {\left\{ \cosh{ \left(
                             \gamma \cdot C p\left( {\widehat{t}} ; \gamma \right) 
                              \right)} 
                        \right\}}^2 } 
                 } 
              d{\widehat{t}} } 
          \right) }^k
          } .
\end{gather*}
If~$ b_{n-1} $ is not zero, then the RHS equals:
\begin{gather*}
      I -
      B C \, 
       {\frac{ 1 - 
       \exp  
       { 
       {\mathbf{-}} \, b_{n-1} \,
          \int_{ t_0 - \epsilon }^{ t_0 + \epsilon } {
                 {\frac { \gamma } 
                      { {\left\{ \cosh{ \left(
                             \gamma \cdot C p\left( {\widehat{t}} ; \gamma \right) 
                              \right)} 
                        \right\}}^2 } 
                 } 
              d{\widehat{t}} } 
          } 
      }
      {b_{n-1}}
      } ,
\end{gather*}
and if~$ b_{n-1} $ equals zero, then the RHS equals
\begin{gather*}
      I -
      B C \, 
          \int_{ t_0 - \epsilon }^{ t_0 + \epsilon } {
                 {\frac { \gamma } 
                      { {\left\{ \cosh{ \left(
                             \gamma \cdot C p\left( {\widehat{t}} ; \gamma \right) 
                              \right)} 
                        \right\}}^2 } 
                 } 
              d{\widehat{t}} } 
       .
\end{gather*}
If $ \gamma $ tends to infinity, then~$ 
{\Phi \left( t_0 + \epsilon , t_0 - \epsilon \right)  } $
tends to
\begin{gather*}
\begin{cases}
      I - B C \, 
      {\frac{ 1 -   
      e^{ 
           {\mathbf{-}} \, b_{n-1} \, 
        \left( 
            { \frac
                {1}  {
                     \left\lvert
                       \rho^{\mathbf\huge{-}}\left( t_0  \right)
                     \right\rvert 
                      }
            }
           +
           {  \frac
                {1}   {
                     \left\lvert
                       \rho^{\mathbf\huge{+}}\left( t_0  \right) 
                     \right\rvert  
                       }
           } 
        \right)
       }}
      {b_{n-1}}
      } 
             ,
             & \text{if} \, b_{n-1} \neq 0, 
             \\
      I - B C \, 
      { 
        \left( 
            { \frac
                {1}  {
                     \left\lvert
                       \rho^{\mathbf\huge{-}}\left( t_0  \right)
                     \right\rvert 
                      }
            }
           +
           {  \frac
                {1}   {
                     \left\lvert
                       \rho^{\mathbf\huge{+}}\left( t_0  \right) 
                     \right\rvert  
                       }
           } 
        \right)
       }
             ,
             & \text{if} \, b_{n-1} =  0 . 
         \end{cases}
\end{gather*}
The above expression is different from that in Equations~28,~34 in~\cite{majhiAtherton1998stabilityOfLimitCyclesInRelay}. Those equations were obtained with an arbitrary correction factor. Our derivation is cleaner.

If the plant has a high relative degree, and the periodic orbit in question is symmetric and unimodal, then the output of the SFS as that of the RFS is well approximated by a pure sinusoid. Then the 
limit:~$   
    {
        {1} / {
        \left\lvert \rho^{\mathbf\huge{-}}\left( t_0  \right) \right\rvert 
    }}
  +
    { 
        {1} /  {
       \left\lvert   \rho^{\mathbf\huge{+}}\left( t_0  \right) \right\rvert  
   }
   }
$
is well approximated by~$
    T_{\rm{period}}  / \left( \pi  M_{\rm{peak}} \right) ,
$
where~$
   M_{\rm{peak}} 
$
denotes the peak amplitude of the output of the SFS/RFS, for the said periodic orbit. Hence, if the periodic orbit is symmetric and unimodal, then we can use: 
\begin{align*}
      \Phi \left( T_{\rm{period}} , 0 \right)  
      & \approx
   {\left(
    \left[
        I
          - B C 
                 {\frac { T_{\rm{period}}   } 
                      { \pi M_{\rm{peak}} 
                      } 
                 } 
          \right]
    \,
    e^{A T_{\rm{period}} / 2}
    \right)}^2 .
\end{align*}

\appendix[A. Two jacobians for the  Poincaré map\label{section:appendixjacobian}]
Consider a switched system, where the state space equals~${\mathbb{R}}^n ,$ where the dynamics between any two successive switches is affine and time-invariant, and where the switching happens at the crossing of hyperplanes. Consider two successive switch times. 
Then 
we call as the {\textit{launching hyperplane}}~(${\mathcal{L}}$) the hyperplane whose crossing was marked by the first time, 
and we call as the {\textit{hitting hyperplane}}~(${\mathcal{H}}$)  the hyperplane whose crossing is marked by second time.
In the Relay feedback system that we have studied, these two hyperplanes are one and the same. But in general they could be different, and do not even need to be parallel.

Consider a trajectory that starts at~$x \in {\mathcal{L}} $ and after time~$\tau,$ crosses~$ {\mathcal{H}} $ at the point~$\xi .$
There are two jacobians:
\begin{itemize}
    \item{a matrix of size $ n \times n$ 
        obtained by linearizing the first crossing map for crossing~$ {\mathcal{H}} , $ from an $n$-dimensional open ball around~$x$ in~${\mathbb{R}}^n$ to an $n$-dimensional open ball around~$\xi.$ We call this the unrestricted jacobian, and denote it by the symbol~$ J_{n}. $}
    \item{and a matrix of size $ n - 1 \times n - 1 $  obtained by linearizing the first crossing map across~${\mathcal{H}},$ from an $n - 1$-dimensional open ball around~$x$ in~${\mathcal{L}}$ to an $n-1$-dimensional open ball around~$\xi \in {\mathcal{H}}.$ We call this the restricted jacobian, and denote  it  by the symbol~$ J_{n-1} . $ Some extra care is needed for calculating this jacobian, because we need to use the same ${n-1}$-dimensional coordinate system for expressing the vectors in the launching and hitting hyperplanes.}
\end{itemize}

The unrestricted jacobian of the first crossing map from any point~$x$ in the ${\mathbb{R}}^n$-valued state space to the hitting hyperplane can be written~({\AA}str{\"o}m~\cite{astrom1995oscillationsRelay})
as a composition of 
the following two linear transformations: (i)~the fundamental matrix of the linearized flow,
from time zero to the time of first crossing of the hitting hyperplane, and (ii)~the oblique projection onto the hitting hyperplane, 
along the direction of the vector field at the point of crossing.

In specific, {\AA}str{\"o}m's jacobian is given by Equation~\eqref{eqn:astromjacobian}. 

\subsection*{A.1 When launching and hitting hyperplanes are not parallel}
 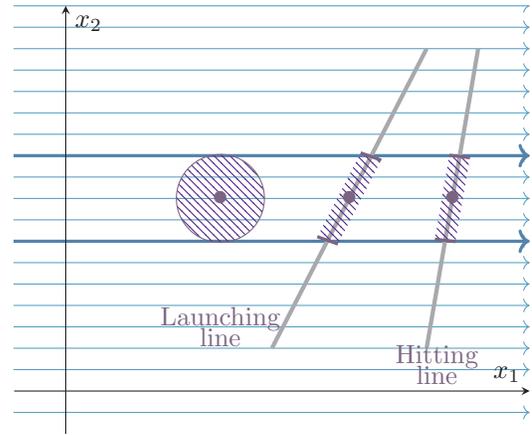
\begin{figure}
\begin{center}
 \begin{tikzpicture}
     \begin{axis}[
          xmin=-0.5, xmax=4.5,
          ymin=-0.5, ymax=4.5,
          axis lines=center,
          axis on top=true,
          domain=-4.5:4.5,
          ylabel=$x_2$,
          xlabel=$x_1$,
          ticks = none,
      ]
  \foreach \n in {1,...,18} {
      \addplot [>->,mark=none,,draw= blue3, samples = 5] {0.25*\n};
  }
      \addplot [>->,mark=none,,draw=blue3, samples = 5] {-0.25};
  \foreach \n in {7,11} {
      \addplot [>->,mark=none,very thick,draw=blue4, samples = 5] {0.25*\n};
  }
  \addplot [mark=none,draw=black!40,ultra thick, samples = 5] coordinates 
  {(2,0.5)
   (3.5,4)};
  \addplot [mark=none,draw=black!40,ultra thick, samples = 5] coordinates 
  {(3.5,0.5)
   (4,4)};
   \draw[c5, pattern = north west lines , pattern color = c5 ] (1.5,2.25) circle (5.85mm);
   \draw [draw=none,pattern = north west lines , pattern color = c5] 
     (2.75-0.21+0.1, 2.25-7*0.21/3-3*0.1/7)
   --(2.75-0.21-0.1, 2.25-7*0.21/3+3*0.1/7)
   --(2.75+0.21-0.1, 2.25+7*0.21/3+3*0.1/7)
   --(2.75+0.21+0.1, 2.25+7*0.21/3-3*0.1/7)
   --(2.75-0.21+0.1, 2.25-7*0.21/3-3*0.1/7) ;
   \draw [draw=c5, very thick]
     (2.75-0.21+0.1, 2.25-7*0.21/3-3*0.1/7)
   --(2.75-0.21-0.1, 2.25-7*0.21/3+3*0.1/7);
   \draw [draw=c5, very thick]
     (2.75+0.21-0.1, 2.25+7*0.21/3+3*0.1/7)
   --(2.75+0.21+0.1, 2.25+7*0.21/3-3*0.1/7);
   \draw [draw=none,pattern = north west lines , pattern color = c5] 
     (3.75-0.07+0.1, 2.25-7*0.07-0.1/7)
   --(3.75-0.07-0.1, 2.25-7*0.07+0.1/7)
   --(3.75+0.07-0.1, 2.25+7*0.07+0.1/7)
   --(3.75+0.07+0.1, 2.25+7*0.07-0.1/7)
   --(3.75-0.07+0.1, 2.25-7*0.07-0.1/7) ;
   \draw [draw=c5, very thick]
     (3.75-0.07+0.1, 2.25-7*0.07-0.1/7)
   --(3.75-0.07-0.1, 2.25-7*0.07+0.1/7);
  \node[] at (1.5,0.85) {{{\textcolor{c5}{Launching}}}};
  \node[below] at (1.5,0.85) {{{\textcolor{c5}{line}}}};
  \node[] at (3.6,0.40) {{{\textcolor{c5}{Hitting}}}};
  \node[below] at (3.6,0.40) {{{\textcolor{c5}{line}}}};
   \draw [draw=c5, very thick]
     (3.75+0.07-0.1, 2.25+7*0.07+0.1/7)
   --(3.75+0.07+0.1, 2.25+7*0.07-0.1/7);
   \node at (1.5,2.25) {{\large{\textcolor{c5}{${\mathbf{\bullet}}$}}}};
   \node at (2.75,2.25) {{\large{\textcolor{c5}{${\mathbf{\bullet}}$}}}};
   \node at (3.75,2.25) {{\large{\textcolor{c5}{${\mathbf{\bullet}}$}}}};
\end{axis}
\end{tikzpicture}
\end{center}
\caption{Example of a Poincaré map whose local stability is not predicted by the spectral radius of $J_{n} . $ 
\label{fig:jacobianExample}%
}
\end{figure}
\begin{example}
    Consider the affine flow induced by the  following ODE having a constant right hand side:
    \begin{align*}
        {\frac{d}{dt}}
        \begin{pmatrix}
            x_1 \\
            x_2
        \end{pmatrix}
        & =  
        \begin{pmatrix}
            1 \\
            0
        \end{pmatrix} .
    \end{align*}
    As shown in Figure~\ref{fig:jacobianExample}, the trajectories in phase space are parallel straight lines. 
    Since the coefficient matrix~$A$ is the zero matrix, the exponential~$ e^{A t}  $ is the identity matrix.

    Thus, regardless of the orientation of the hitting line, the unrestricted jacobian~$J_n$ is simply equal to the projection along the $x_1$-axis onto it. Since the eigenvalues of any  projection matrix equal either~$ 0 $ or~$ 1 ,$ it follows that the spectral radius of~$J_n$ is exactly~$1.$

The restricted jacobian is less than~$1,$ because the Poincaré map is a contraction. The Poincaré map takes line segments in~$ {\mathcal{L}} $  to shorter line segments on~$ {\mathcal{H}} .$

\end{example}
\subsection*{A.2 When launching and hitting hyperplanes are parallel}
\fboxsep=0pt
\begin{lemma}
For the $n \times n$ square matrix~$P$ that represents a linear transformation from~${\mathbb{R}}^n$ to itself, assume that:
\begin{itemize}
    \item{The rank 
        of~$P$ is less than or equal to~$ n - 1 .$}
    \item{%
            $P$~maps
            the whole of~${\mathbb{R}}^n$  to 
            either the whole of, or a subspace of the following linear space:
        \begin{align*}
            {\mathcal{S}}  & \triangleq  \left\{  x \in  {\mathbb{R}}^n 
                \left\lvert  C x = 0     \right. \right\} ,
        \end{align*}
         where $C$~is the following row vector:
        \begin{align*}
            C & = 
        \begin{bmatrix} 
               0 & \cdots & 0 & 1 
            \end{bmatrix}
            .
        \end{align*}
        }
\end{itemize}
Let us decompose the matrix~$P$ into the following blocks:
        \begin{align*}
            P & = 
           \underbracket[.4pt]{%
            \begin{bmatrix}
                 \underbracket[.4pt]{%
        P_{%
            {\textcolor{tolRedForContrast}{\colorbox{tolBlueForContrast}{
            1:n-1, 1:n-1
            }}} 
        }
                  }_{n - 1}
                &
        P_{%
            {\textcolor{tolDarkRed}{\colorbox{tolPaleYellow}{
            1:n-1, n
            }}} 
        }
               \\
        P_{%
            {\textcolor{tolDarkRed}{\colorbox{tolPaleYellow}{
            n , 1:n-1
            }}} 
        }
                &
        P_{%
            n , n
        }
            \end{bmatrix}
            }_{n}  .
        \end{align*}
    Since $C P x = 0$ for every~$x \in {\mathbb{R}}^n, $ 
    it must follow that
    \begin{align}
        P x & = 
        \begin{pmatrix}
        P_{%
            {\textcolor{tolRedForContrast}{\colorbox{tolBlueForContrast}{
            1:n-1, 1:n-1
            }}} 
        }
        x_{%
            {\textcolor{tolDarkRed}{\colorbox{tolPaleYellow}{
            1:n-1
            }}} 
        }
            + 
        P_{%
            {\textcolor{tolDarkRed}{\colorbox{tolPaleYellow}{
            1:n-1, n
            }}} 
        }
        x_{%
            n
        }
            \\
             0
        \end{pmatrix}
        \label{eqn:lastElementZero}
        .
    \end{align}
In other words,  both the 
row vector~$ 
        P_{%
            {\textcolor{tolDarkRed}{\colorbox{tolPaleYellow}{
            n, 1:n-1
            }}} 
        }
    $ 
and the real number~$ 
        P_{%
            n, n
        }
    $ must be zeros.

Clearly, the square matrix~$  
        P_{%
            {\textcolor{tolRedForContrast}{\colorbox{tolBlueForContrast}{
            1:n-1, 1:n-1
            }}} 
        }
    $ represents a linear transformation from the
    hyperplane~$ {\mathcal{S}} $ to itself.

    Then the matrices~$P$ and $
        P_{%
            {\textcolor{tolRedForContrast}{\colorbox{tolBlueForContrast}{
            1:n-1, 1:n-1
            }}} 
        }
    $ have the same non-zero eigenvalues, with the same multiplicities.
\label{lemma:whenLaunchingAndHittingHyperplanesAreParallel}
\end{lemma}
\begin{proof}
    First we prove that every eigenvalue of~{\textit{%
$ 
        P_{%
            {\textcolor{tolRedForContrast}{\colorbox{tolBlueForContrast}{
             1:n-1, 1:n-1 
            }}} 
        }
    $}}
    is an eigenvalue of~$P . $
    Suppose~$\lambda$ is an eigenvalue of~{\textit{%
$ 
        P_{%
            {\textcolor{tolRedForContrast}{\colorbox{tolBlueForContrast}{
             1:n-1, 1:n-1 
            }}} 
        }
    $}}
of multiplicity~$k$. Then there must be exactly~$k$ generalized eigenvectors
    \begin{gather*}
        v_1 , \ldots , v_k \in {\mathbb{C}}^{n-1}
    \end{gather*}
    with possibly complex entries such that for $ 0 \le i \le k $:
    \begin{align*}
        \text{~{\textit{%
            $
        P_{%
            {\textcolor{tolRedForContrast}{\colorbox{tolBlueForContrast}{
            1:n-1, 1:n-1
            }}} 
        } $  }} }
        \,  v_i 
        & = \lambda v_i ,
    \end{align*}
    Thus from Equation~\ref{eqn:lastElementZero}, it follows that
    \begin{align*}
        P 
        \begin{pmatrix}
            v_i
            \\
            0
        \end{pmatrix}
            = &
            \lambda 
        \begin{pmatrix}
            v_i
            \\
            0
        \end{pmatrix}
            , 
        \ \text{for} \ 1 \le i \le k .
    \end{align*}
    Hence~$\lambda$ is an eigenvalue of~$P$ with a multiplicity of at least~$k .$
    Next we prove 
    that $P$ cannot have a non-zero eigenvalue that is not 
    an eigenvalue of~{\textit{%
$ 
        P_{%
            {\textcolor{tolRedForContrast}{\colorbox{tolBlueForContrast}{
            1:n-1, 1:n-1
            }}} 
        }
    . $ }}

    Suppose for argument's sake that~$\lambda_{\text{extra}}, v_{\text{extra}} \in {\mathbb{C}}^n $ is an eigenvalue, generalized eigenvector pair for~$P ,$ but not
    for~{\textit{%
$ 
        P_{%
            {\textcolor{tolRedForContrast}{\colorbox{tolBlueForContrast}{
            1:n-1, 1:n-1
            }}} 
        }
    . $ }}
    This implies that~$ C v_{\text{extra}} \neq 0 . $

    Were~$  \lambda_{\text{extra}}  $ to be nonzero, then we could write:
\begin{gather*}
C v_{\text{extra}} \ = \  
 {\frac{1}{ \lambda_{\text{extra}} }}
C P v_{\text{extra}}
    . 
\end{gather*}
    But since~$C P x  \equiv 0$ for all~$x \in {\mathbb{C}}^n  , $ it follows that the vector~$v_{\text{extra}}$ could be a generalized eigenvector of~$P,$ only if the corresponding eigenvalue~$  \lambda_{\text{extra}}  $ is zero.

    Since~$P$ can have only one more eigenvalue than~{\textit{%
$ 
        P_{%
            {\textcolor{tolRedForContrast}{\colorbox{tolBlueForContrast}{
            1:n-1, 1:n-1
            }}} 
        }
    ,
    $}} it follows that the additional one can only be the zero eigenvalue. Hence the non-zero eigenvalues are in common to both the matrices, and for these eigenvalues the multiplicities are the same.
\end{proof}

\begin{corollary}
For a Poincaré map arising in a switched linear system, assume that the launching and hitting hyperplanes are parallel to each other. Then the restricted and unrestricted jacobians have the same non-zero eigenvalues with the same multiplicities. Therefore both the jacobians have the same spectral radius.
\label{corollary:whenLaunchingAndHittingHyperplanesAreParallel}
\end{corollary}


\appendix[B.~RFS output behaviours for some plants\label{appendix:tableOfRFSoutputs}]
The following is  a random sample of RFS behaviours. Each entry gives the plant transfer function, its step response, and a short, partial list of RFS behaviours observed.

There is no clearly visible marker in the shape of the plant step reponse, that could help predict RFS oscillations. What we do notice is that: (1)~there seem to be many minimum phase plants whose step responses have an undershoot, and (2)~some plants give rise to very spiky trajectories of the RFS.
\tablecaption{RFS output for plants that are not BRL-URF} \label{table:timeline}
\tablelasttail{\bottomrule}

\tablefirsthead{
\toprule
\multicolumn{1}{c}{\textit{Plant}}  & \multicolumn{1}{c}{\textit{RFS behaviour}}  \\
\midrule
}
\tablehead{
\multicolumn{1}{c}{\textit{Plant}}  & \multicolumn{1}{c}{\textit{RFS behaviour}}  \\
\midrule
}

\begin{xtabular}{p{3.5cm} p{4.5cm}}
    \multicolumn{2}{c}{\textcolor{cyanTeal}{\textit{Known from the previous literature:}}} 
\\[6pt]
\centerline{%
\( \displaystyle
    {\frac{1}{ s^2 + \alpha s + \beta}} , 
\)%
    }
with 
$ \alpha , \beta  > 0 $
  & 
   Origin is globally asymptotically stable~\cite{mabenRabi2021secondOrderRelayOscillations,rabi2018relay}.
\begin{center}
  \begin{tikzpicture}
\begin{axis}[
              yscale = 1.3,
              yticklabel pos = right,
              width=5cm,
              height=2.3cm,
              ticks=none,
              xmin = -0.3,
              xmax = 2.5,
              ymin = -0.4,
              ymax = 1.2,
           ]
\addplot[thick, tolRedMain,
grid=both, every major grid/.style={gray, opacity=0.5}, 
    domain = -0.5:1.5,
    ]
    file {figuresAndData/relayOutputSecondOrderNum+0+1.csv};
\end{axis}
\end{tikzpicture}
\end{center}
\\
    \centerline{%
\( \displaystyle
        {\frac{1}{ {\left( s + 1 \right)}^3 }} , 
\)
    }
  & 
   Has a unique symmetric, and unimodal limit cycle - this is locally stable~\cite{astrom1995oscillationsRelay}. Trajectories converge to this limit cycle from the origin and from $ {\left[ 1 , 1 , 1 \right]}^T$  
\begin{center}
  \begin{tikzpicture}
\begin{axis}[
              yscale = 1.3,
              yticklabel pos = right,
              width=5cm,
              height=2.3cm,
              ticks=none,
           ]
\addplot[thick, cyanTeal,
grid=both, every major grid/.style={gray, opacity=0.5}, smooth,
    ]
    file {figuresAndData/relayOutputThirdOrderNumerator+1_triplePole.csv};
\end{axis}
\end{tikzpicture}
\end{center}
   \\
    \multicolumn{2}{c}{\textcolor{cyanTeal}{\textit{Simulations of 3rd order plants:}}} 
\\[6pt]
    \centerline{%
\( \displaystyle
        {\frac{1}{  s^3 + 3 s^2 + 5 s +  6 }}  
\)
    }
   & 
        Has a symmetric, unimodal limit cycle. Trajectories converge to this limit cycle from the origin, from $ {\left[ 0.1 , 0.1 , 0.1 \right]}^T$  and from $ {\left[ 1 , 1 , 1 \right]}^T.$ 
\begin{center}
  \begin{tikzpicture}
\begin{axis}[
              yscale = 1.3,
              yticklabel pos = right,
              width=5cm,
              height=2.3cm,
              ticks=none,
           ]
\addplot[thick, cyanTeal,
grid=both, every major grid/.style={gray, opacity=0.5}, smooth,
    ]
    file {figuresAndData/relayOutputThirdOrderNumerator000+1_startingAtOrigin.csv};
\end{axis}
\end{tikzpicture}
\end{center}
\\
    \centerline{%
\( \displaystyle
        {\frac{s + 1}{  s^3 + 3 s^2 + 5 s +  6 }}  
\)
    }
   & Trajectories from $ {\left[ 0.1 , 0.1 , 0.1 \right]}^T$  and from $ {\left[ 1 , 1 , 1 \right]}^T$ converge to the origin.  Trajectory from the origin chatters at the origin.  
\begin{center}
  \begin{tikzpicture}
\begin{axis}[
              yscale = 1.3,
              yticklabel pos = right,
              width=5cm,
              height=2.3cm,
              ticks=none,
           ]
\addplot[thick, tolRedMain,
grid=both, every major grid/.style={gray, opacity=0.5}, 
    ]
    file {figuresAndData/relayOutputThirdOrderNumerator+1+1.csv};
\end{axis}
\end{tikzpicture}
\end{center}
\\
    \centerline{%
\( \displaystyle
        {\frac{ 4  s + 1   }
              {   s^3 + 3 s^2 + 5 s +  6 }
        }  
\)
    }
    Plant has a zero at $- 0.25 $.
    \vspace*{4pt}
    This plant has two repeated zeroes at: $ - 0.5 .$ 
    It is minimum phase. But yet, its step response shows an inverse response, like some minimum phase transfer functions discovered by Maciejowski~\cite{janMaciejowski2018rhpZeroesNotNecessaryForInverseResponse}.
\begin{center}
  \begin{tikzpicture}
\begin{axis}[%
              yscale = 1.3,
              yticklabel pos = right,
              width=5cm,
              height=2.3cm,
              ticks=none,
           ]
            \begin{scope}
            \clip (axis cs:0,0) rectangle (axis cs:140,8);
        \addplot [ tolDarkYellow,] 
    file {figuresAndData/stepResponseWithInverseResponse_3rdthOrder_numerator+4+1_scientific.csv};
            \end{scope}
            \begin{scope}
            \clip (axis cs:0,0) rectangle (axis cs:140,-5);
        \addplot [ tolPaleRed, ultra thick] 
    file {figuresAndData/stepResponseWithInverseResponse_3rdthOrder_numerator+4+1_scientific.csv};
            \end{scope}
\end{axis}
\end{tikzpicture}
\end{center}
   & 
        Trajectory starting from the origin chatters there.
        Trajectories 
        from~$ {\left[ 0.01 , 0.01 , 0.01 , 0.01 \right]}^T , $ 
        from~$ {\left[ 0.1 , 0.1 , 0.1 , 0.1 \right]}^T , $ 
        from~$ {\left[ 1 , 1 , 1 , 1 \right]}^T , $ 
        and from~$ {\left[ 10 , 10 , 10 , 10 \right]}^T  $ 
        settle into the same symmetric, unimodal limit cycle.
\begin{center}
  \begin{tikzpicture}
\begin{axis}[%
              yscale = 1.3,
              yticklabel pos = right,
              width=5cm,
              height=2.3cm,
              ticks=none,
           ]
        \addplot [tolRedMain,  thick]
    file {figuresAndData/asymptoticDecay_3rdthOrder_numerator+4+1_startingAtAllTens.csv};
\end{axis}
\end{tikzpicture}
\end{center}
\\
    \centerline{%
\( \displaystyle
        {\frac{s^2 + s + 1}{  s^3 + 3 s^2 + 5 s +  6 }}  
\)
    }
   & 
        Trajectories starting from the origin, and from~$ {\left[ 1 , 1 , 1 \right]}^T$ chatter on a neighbourhood that is around the origin. 
\\
    \centerline{%
\( \displaystyle
        {\frac{s^2 - 2 s + 2}{  s^3 + 3 s^2 + 5 s +  6 }}  
\)
    }
    \vspace*{2pt}
    Plant has a complex conjugate pair of zeroes on the open RHP.
   & 
        Trajectories starting from the origin, and from~$ {\left[ 1 , 1 , 1 \right]}^T$ chatter on a neighbourhood that is around the origin. 
\\
    \vspace*{1pt}
    \centerline{%
\( \displaystyle
        {\frac{s^2 - 2 s + 1}{  s^3 + 3 s^2 + 5 s +  6 }}  
\)
    }
    \vspace*{2pt}
    Plant has a real RHP zero at +1, with multiplicity 2.
   & 
        Trajectories starting from the origin, and from~$ {\left[ 1 , 1 , 1 \right]}^T$ chatter on a neighbourhood that is around the origin. 
\\
    \centerline{%
\( \displaystyle
        {\frac{s}{  s^3 + 3 s^2 + 5 s +  6 }}  
\)
    }
   & 
        Trajectories starting from the origin, and from~$ {\left[ 1 , 1 , 1 \right]}^T$ chatter on a neighbourhood that is around the origin. 
\\
    \centerline{%
\( \displaystyle
        {\frac{s^2 - 2 s + 2}{  s^3 + 3 s^2 + 5 s +  6 }}  
\)
    }
    Plant has a complex conjugate pair of zeroes on the open RHP.
   & 
        Trajectories starting from the origin, and from~$ {\left[ 1 , 1 , 1 \right]}^T$ chatter on a neighbourhood that is around the origin. 
\\
    \centerline{%
\( \displaystyle
        {\frac{s}{  s^3 + 3 s^2 + 5 s +  6 }}  
\)
    }
    Plant has a zero at 0.
   & 
        Trajectories starting from the origin, and from~$ {\left[ 1 , 1 , 1 \right]}^T$ settle into the same limit cycle.
\begin{center}
  \begin{tikzpicture}
\begin{axis}[
              yscale = 1.3,
              yticklabel pos = right,
              width=5cm,
              height=2.3cm,
              ticks=none,
           ]
\addplot[thick, cyanTeal,
grid=both, every major grid/.style={gray, opacity=0.5}, 
    ]
    file {figuresAndData/relayOutputThirdOrderNumeratorS.csv};
\end{axis}
\end{tikzpicture}
\end{center}
\\
\midrule
    \multicolumn{2}{c}{\textcolor{cyanTeal}{\textit{Simulations of 4th order plants:}}} 
\\[6pt]
    \centerline{%
\( \displaystyle
        {\frac{s}{  s^3 + 3 s^2 + 5 s +  6 }}  
\)
    }
   & 
        Trajectories starting from the origin, and from~$ {\left[ 1 , 1 , 1 \right]}^T$ chatter on a neighbourhood that is around the origin. 
\\
    \centerline{%
\( \displaystyle
        {\frac{s^2 + 2 s + 1}{ s^4 + 2 s^3 + 3 s^2 + 3 s +  2 }}  
\)
    }
    Plant has a complex conjugate pair of zeroes on the open RHP.
    \vspace*{4pt}
    This plant has two repeated zeroes at: $ - 1 .$   It is minimum phase. But yet, its step response shows an inverse response, like some minimum phase transfer functions discovered by Maciejowski~\cite{janMaciejowski2018rhpZeroesNotNecessaryForInverseResponse}.
\begin{center}
  \begin{tikzpicture}
\begin{axis}[%
              yscale = 1.3,
              yticklabel pos = right,
              width=5cm,
              height=2.3cm,
              ticks=none,
           ]
            \begin{scope}
            \clip (axis cs:0,0) rectangle (axis cs:140,40);
        \addplot [ tolDarkYellow,] 
    file {figuresAndData/stepResponseWithInverseResponse_4thOrder_numeratorDoublePlusOne.csv};
            \end{scope}
            \begin{scope}
            \clip (axis cs:0,0) rectangle (axis cs:140,-40);
        \addplot [ tolPaleRed, ultra thick] 
    file {figuresAndData/stepResponseWithInverseResponse_4thOrder_numeratorDoublePlusOne.csv};
            \end{scope}
\end{axis}
\end{tikzpicture}
\end{center}
   & 
        Trajectory starting from the origin chatters there.
        Trajectories 
        from~$ {\left[ 0.01 , 0.01 , 0.01 , 0.01 \right]}^T , $ 
        from~$ {\left[ 0.1 , 0.1 , 0.1 , 0.1 \right]}^T , $ 
        from~$ {\left[ 1 , 1 , 1 , 1 \right]}^T , $ 
        and from~$ {\left[ 10 , 10 , 10 , 10 \right]}^T  $ 
        asymptotically reach the origin.
\begin{center}
  \begin{tikzpicture}
\begin{axis}[%
              yscale = 1.3,
              yticklabel pos = right,
              width=5cm,
              height=2.3cm,
              ticks=none,
           ]
        \addplot [tolRedMain,  thick]
    file {figuresAndData/asymptoticDecay_4thOrder_numeratorDoublePlusOne_startingAtAllTens.csv};
\end{axis}
\end{tikzpicture}
\end{center}
\\
    \centerline{%
\( \displaystyle
        {\frac{ {\left( 2 s + 1 \right)}^2 }
              { s^4 + 2 s^3 + 3 s^2 + 3 s +  2 }
        }  
\)
    }
    Plant has a zero at $- 1 / 2$, with multiplicity two.
    \vspace*{4pt}
    This plant has two repeated zeroes at: $ - 0.5 .$   It is minimum phase. But yet, its step response shows an inverse response, like some minimum phase transfer functions discovered by Maciejowski~\cite{janMaciejowski2018rhpZeroesNotNecessaryForInverseResponse}.
\begin{center}
  \begin{tikzpicture}
\begin{axis}[%
              yscale = 1.3,
              yticklabel pos = right,
              width=5cm,
              height=2.3cm,
              ticks=none,
           ]
            \begin{scope}
            \clip (axis cs:0,0) rectangle (axis cs:140,40);
        \addplot [ tolDarkYellow,] 
    file {figuresAndData/stepResponseWithInverseResponse_4thOrder_numerator+4+4+1.csv};
            \end{scope}
            \begin{scope}
            \clip (axis cs:0,0) rectangle (axis cs:140,-40);
        \addplot [ tolPaleRed, ultra thick] 
    file {figuresAndData/stepResponseWithInverseResponse_4thOrder_numerator+4+4+1.csv};
            \end{scope}
\end{axis}
\end{tikzpicture}
\end{center}
   & 
        Trajectory starting from the origin chatters there.
        Trajectories 
        from~$ {\left[ 0.01 , 0.01 , 0.01 , 0.01 \right]}^T , $ 
        from~$ {\left[ 0.1 , 0.1 , 0.1 , 0.1 \right]}^T , $ 
        from~$ {\left[ 1 , 1 , 1 , 1 \right]}^T , $ 
        and from~$ {\left[ 10 , 10 , 10 , 10 \right]}^T  $ 
        asymptotically reach the origin.
\begin{center}
  \begin{tikzpicture}
\begin{axis}[%
              yscale = 1.3,
              yticklabel pos = right,
              width=5cm,
              height=2.3cm,
              ticks=none,
           ]
        \addplot [tolRedMain,  thick]
    file {figuresAndData/asymptoticDecay_4thOrder_numerator+4+4+1_startingAtAllOnes.csv};
\end{axis}
\end{tikzpicture}
\end{center}
\\
    \centerline{%
\( \displaystyle
        {\frac{ {\left(  s + 1 / 4 \right)}^2 }
              { s^4 + 2 s^3 + 3 s^2 + 3 s +  2 }
        }  
\)
    }
    Plant has a zero at $- 4 $, with multiplicity two.
    \vspace*{4pt}
    This plant has two repeated zeroes at: $ - 0.5 .$ 
    It is minimum phase. But yet, its step response shows an inverse response, like some minimum phase transfer functions discovered by Maciejowski~\cite{janMaciejowski2018rhpZeroesNotNecessaryForInverseResponse}.
\begin{center}
  \begin{tikzpicture}
\begin{axis}[%
              yscale = 1.3,
              yticklabel pos = right,
              width=5cm,
              height=2.3cm,
              ticks=none,
           ]
            \begin{scope}
            \clip (axis cs:0,0) rectangle (axis cs:140,4);
        \addplot [ tolDarkYellow,] 
    file {figuresAndData/stepResponse_4thOrder_numerator+1+8over16+1over16_scientific.csv};
            \end{scope}
            \begin{scope}
            \clip (axis cs:0,0) rectangle (axis cs:140,-4);
        \addplot [ tolPaleRed, ultra thick] 
    file {figuresAndData/stepResponse_4thOrder_numerator+1+8over16+1over16_scientific.csv};
            \end{scope}
\end{axis}
\end{tikzpicture}
\end{center}
   & 
        Trajectory starting from the origin chatters there.
        Trajectories 
        from~$ {\left[ 0.01 , 0.01 , 0.01 , 0.01 \right]}^T , $ 
        from~$ {\left[ 0.1 , 0.1 , 0.1 , 0.1 \right]}^T , $ 
        from~$ {\left[ 1 , 1 , 1 , 1 \right]}^T , $ 
        and from~$ {\left[ 10 , 10 , 10 , 10 \right]}^T  $ 
        settle into the same symmetric, unimodal limit cycle.
\begin{center}
  \begin{tikzpicture}
\begin{axis}[%
              yscale = 1.3,
              yticklabel pos = right,
              width=5cm,
              height=2.3cm,
              ticks=none,
           ]
        \addplot [tolRedMain,  thick]
    file {figuresAndData/asymptoticDecay_4thOrder_numerator+1+2over4+1over16_startingAtAllOnes.csv};
\end{axis}
\end{tikzpicture}
\end{center}
\\
    \centerline{%
\( \displaystyle
        {\frac{s^2 + 5 s + 6}{ s^4 + 2 s^3 + 3 s^2 + 3 s +  2 }}  
\)
    }
    Plant has a complex conjugate pair of zeroes on the open RHP.
    \vspace*{4pt}
    This plant has zeroes at: $ - 2 , -3 .$   It is minimum phase. But its step response is very oscillatory.
\begin{center}
  \begin{tikzpicture}
\begin{axis}[%
              yscale = 1.3,
              yticklabel pos = right,
              width=5cm,
              height=2.3cm,
              ticks=none,
           ]
            \begin{scope}
            \clip (axis cs:0,0) rectangle (axis cs:140,40);
        \addplot [ tolDarkYellow,] 
    file {figuresAndData/stepResponseWithInverseResponse_4thOrder_numerator+1+5+6.csv};
            \end{scope}
            \begin{scope}
            \clip (axis cs:0,0) rectangle (axis cs:140,-40);
        \addplot [ tolPaleRed, ultra thick] 
    file {figuresAndData/stepResponseWithInverseResponse_4thOrder_numerator+1+5+6.csv};
            \end{scope}
\end{axis}
\end{tikzpicture}
\end{center}
   & 
        Trajectory starting from the origin chatters there.
        Trajectories 
        from~$ {\left[ 0.01 , 0.01 , 0.01 , 0.01 \right]}^T , $ 
        from~$ {\left[ 0.1 , 0.1 , 0.1 , 0.1 \right]}^T , $ 
        from~$ {\left[ 1 , 1 , 1 , 1 \right]}^T , $ 
        and from~$ {\left[ 10 , 10 , 10 , 10 \right]}^T  $ 
        reach the same symmetric, unimodal limit cycle.
\begin{center}
  \begin{tikzpicture}
\begin{axis}[%
              yscale = 1.3,
              yticklabel pos = right,
              width=5cm,
              height=2.3cm,
              ticks=none,
           ]
        \addplot [cyanTeal,  thick]
    file {figuresAndData/limitCycle_4thOrder_numerator+1+5+6_startingAtAllTens.csv};
%
\end{axis}
\end{tikzpicture}
\end{center}
\\
    \centerline{%
\( \displaystyle
        {\frac{s}{  s^3 + 3 s^2 + 5 s +  6 }}  
\)
    }
    Plant has a zero at 0.
   & 
        Trajectories starting from the origin, and from~$ {\left[ 1 , 1 , 1 \right]}^T$ settle into the same limit cycle.
\begin{center}
  \begin{tikzpicture}
\begin{axis}[
              yscale = 1.3,
              yticklabel pos = right,
              width=5cm,
              height=2.3cm,
              ticks=none,
           ]
\addplot[thick, cyanTeal ,
grid=both, every major grid/.style={gray, opacity=0.5}, 
    ]
    file {figuresAndData/relayOutputFourthOrderNumerator0+1-2+1_startingAtAll0pt1s.csv};
\end{axis}
\end{tikzpicture}
\end{center}
\\
    \vspace*{1pt}
    \centerline{%
\( \displaystyle
        {\frac{s^2 - 0.02 s + 0.0001}{ s^4 + 2 s^3 + 3 s^2 + 3 s +  2 }}  
\)
    }
   & 
        Trajectory starting from~$ {\left[ 1 , 1 , 1 \right]}^T$ outputs a seris of spikes that decay to zero.
\begin{center}
    \includegraphics[width=1.1\linewidth]{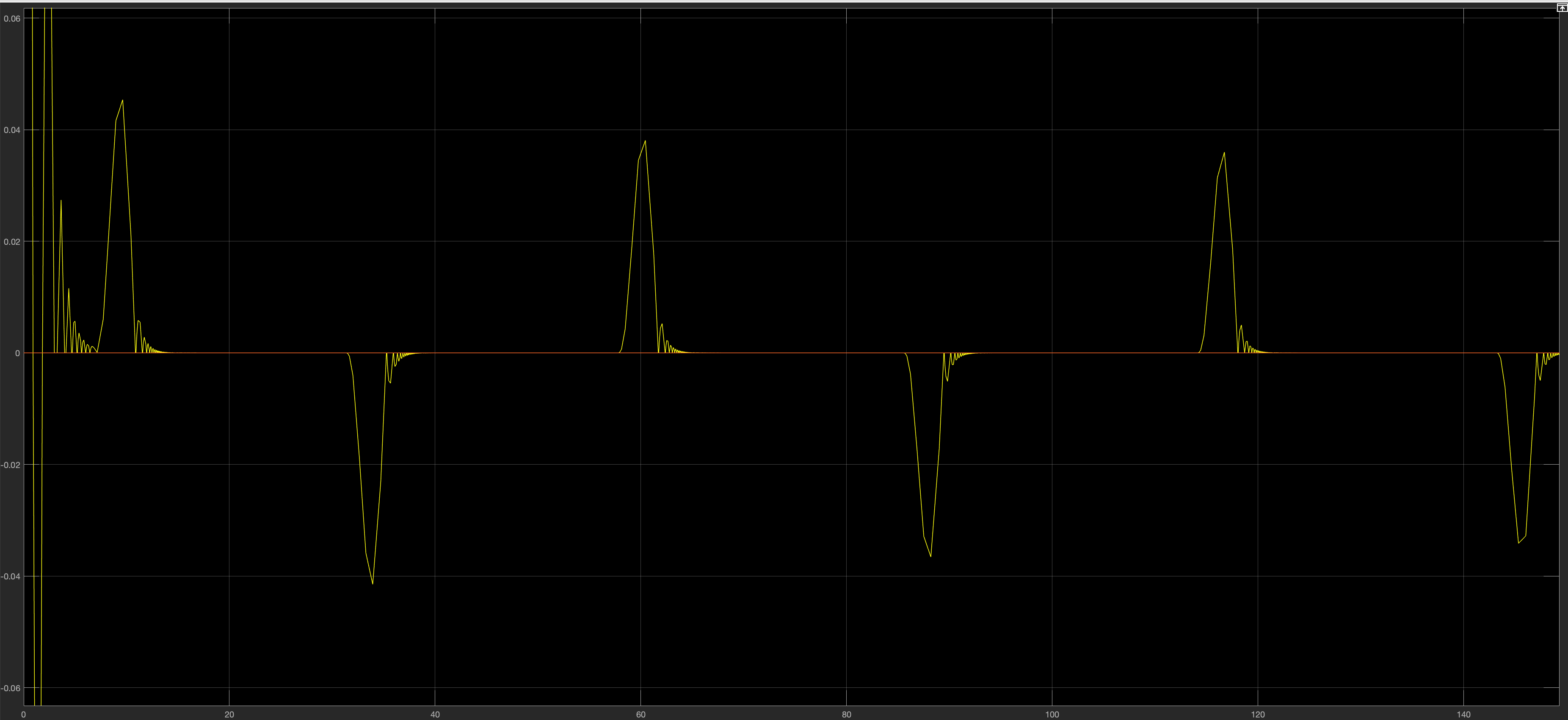}
\end{center}
\\
\midrule
    \multicolumn{2}{c}{\textcolor{cyanTeal}{\textit{Simulations of 5th order plants:}}} 
\\
    \centerline{%
\(  \displaystyle
        {\frac{{\left(s+1\right)}^2}
        { \begin{Bmatrix}
            {\left(s + 0.7\right)}^3 
            {\left(s + 7\right)}^2 
        \end{Bmatrix}
        }
        }  
\)
    }
\begin{center}
  \begin{tikzpicture}
\begin{axis}[%
              yscale = 1.3,
              yticklabel pos = right,
              width=5cm,
              height=2.3cm,
              ticks=none,
           ]
        \addplot [ tolDarkYellow,] 
    file {figuresAndData/stepResponseFifthOrderNumSplus1WholeSqDenSplus0pt7wholeCuSplus7wholeSq.csv};
\end{axis}
\end{tikzpicture}
\end{center}
   &
   There is a stable symmetric and unimodal limit cycle, that passes through the point: 
  $ {\left[
 0.2305 ,  0.4728 , 0.2558 , 0.0138 , 0 
       \right]}^T . $ 
  Trajectory from $ {\left[ 10 , 10 , 10, 10, 10, 10 \right]}^T $ 
   converges to the origin.
\begin{center}
  \begin{tikzpicture}
\begin{axis}[%
              yscale = 1.3,
              yticklabel pos = right,
              width=5cm,
              height=2.3cm,
              ticks=none,
           ]
\addplot[thick, cyanTeal,
grid=both, every major grid/.style={gray, opacity=0.5}, 
    ]
    file {figuresAndData/relayOutputFifthOrderNumSplus1WholeSqDenSplus0pt7wholeCuSplus7wholeSq_limitCycle.csv};
\end{axis}
\end{tikzpicture}
\end{center}
%
\begin{center}
  \begin{tikzpicture}
\begin{axis}[%
              yscale = 1.3,
              yticklabel pos = right,
              width=5cm,
              height=2.3cm,
              ticks=none,
           ]
\addplot[thick, tolRedMain,
grid=both, every major grid/.style={gray, opacity=0.5}, 
    ]
    file {figuresAndData/relayOutputFifthOrderNumSplus1WholeSqDenSplus0pt7wholeCuSplus7wholeSq_decay.csv};
\end{axis}
\end{tikzpicture}
\end{center}
\\
    \centerline{%
\(  \displaystyle
        {\frac{s^2 + 5 s + 6}
        { 
            s^5 + 2 s^4 + 4 s^3 
            + 4 s^2 + 3 s +  1
        }
        }  
\)
    }
   & Trajectories from the origin, 
   from $ {\left[ 0.1 , 0.1 , 0.1 , 0.1 , 0.1 , 0.1 \right]}^T$ 
  and from $ {\left[ 1 , 1 , 1 , 1 , 1 , 1 \right]}^T$ 
   converge to a symmetric unimodal limit cycle. 
   But Trajectories from 
   $ {\left[ 10 , 10 , 10, 10, 10, 10 \right]}^T$ 
   and from
   $ {\left[ 50 , 50 , 50, 50, 50, 50 \right]}^T$ 
   converge to a different symmetric unimodal limit cycle. 
\begin{center}
  \begin{tikzpicture}
\begin{axis}[%
              yscale = 1.3,
              yticklabel pos = right,
              width=5cm,
              height=2.3cm,
              ticks=none,
           ]
\addplot[thick, cyanTeal,
grid=both, every major grid/.style={gray, opacity=0.5}, 
    ]
    file {figuresAndData/limitCycle_5thOrder_numerator+1+5+6_startingAtAllZeroes.csv};
%
%
\end{axis}
\end{tikzpicture}
\end{center}
%
\\
    \centerline{%
\(  \displaystyle
        {\frac{ { \left( s + 1 \right)}^3 }
        { 
            s^5 + 2 s^4 + 4 s^3 
            + 4 s^2 + 3 s +  1
        }
        }  
\)
    }
   & Trajectories from the origin, 
   from $ {\left[ 0.1 , 0.1 , 0.1 , 0.1 , 0.1 , 0.1 \right]}^T$ 
  and from $ {\left[ 1 , 1 , 1 , 1 , 1 , 1 \right]}^T$ 
   converge to a symmetric unimodal limit cycle. 
   But Trajectories from 
   $ {\left[ 10 , 10 , 10, 10, 10, 10 \right]}^T$ 
   and from
   $ {\left[ 50 , 50 , 50, 50, 50, 50 \right]}^T$ 
   converge to a different symmetric unimodal limit cycle. 
\begin{center}
  \begin{tikzpicture}
\begin{axis}[%
              yscale = 1.3,
              yticklabel pos = right,
              width=5cm,
              height=2.3cm,
              ticks=none,
           ]
\addplot[thick, cyanTeal,
grid=both, every major grid/.style={gray, opacity=0.5}, 
    ]
    file {figuresAndData/limitCycle_5thOrder_numeratorTripleZeroPlusOne_startingAtAllOnes.csv};
\end{axis}
\end{tikzpicture}
\end{center}
\\
    \centerline{%
\(  \displaystyle
        {\frac{  {\left(2 s +  1 \right)}^3 }
        { 
            s^5 + 2 s^4 + 4 s^3 
            + 4 s^2 + 3 s +  1
        }
        }  
\)
    }
    \vspace*{4pt}
    This plant has three repeated zeroes at: $ - 1 .$   It is minimum phase. But yet, its step response shows an inverse response, like some minimum phase transfer functions discovered by Maciejowski~\cite{janMaciejowski2018rhpZeroesNotNecessaryForInverseResponse}.
\begin{center}
  \begin{tikzpicture}
\begin{axis}[%
              yscale = 1.3,
              yticklabel pos = right,
              width=5cm,
              height=2.3cm,
              ticks=none,
           ]
            \begin{scope}
            \clip (axis cs:0,0) rectangle (axis cs:140,40);
        \addplot [ tolDarkYellow,] 
    file {figuresAndData/stepResponseWithInverseResponse_5thOrder_numerator+8+12+6+1.csv};
            \end{scope}
            \begin{scope}
            \clip (axis cs:0,0) rectangle (axis cs:140,-40);
        \addplot [ tolPaleRed, ultra thick] 
    file {figuresAndData/stepResponseWithInverseResponse_5thOrder_numerator+8+12+6+1.csv};
            \end{scope}
\end{axis}
\end{tikzpicture}
\end{center}
   & 
   Trajectory from the origin chatters there. 
   Trajectories 
   from $ {\left[ 0.1 , 0.1 , 0.1 , 0.1 , 0.1 , 0.1 \right]}^T , $ 
  from $ {\left[ 1 , 1 , 1 , 1 , 1 , 1 \right]}^T , $ 
   nd from $ {\left[ 10 , 10 , 10, 10, 10, 10 \right]}^T$ 
   asymptotically reach the origin.
\begin{center}
  \begin{tikzpicture}
\begin{axis}[%
              yscale = 1.3,
              yticklabel pos = right,
              width=5cm,
              height=2.3cm,
              ticks=none,
           ]
\addplot[thick, tolRedMain,
grid=both, every major grid/.style={gray, opacity=0.5}, 
    ]
    file {figuresAndData/asymptoticDecay_5thOrder_numerator+8+12+6+1_startingAtAllOnes.csv};
\end{axis}
\end{tikzpicture}
\end{center}
\\
\midrule
    \multicolumn{2}{c}{\textcolor{cyanTeal}{\textit{Simulations of 6th order plants:}}} 
\\[6pt]
    \centerline{%
\(  \displaystyle
        {\frac{s^2}
        { \begin{Bmatrix}
            s^6 + s^5 + 5 s^4 + 3 s^3 \\
            + 5 s^2 +  s +  1
        \end{Bmatrix}
        }
        }  
\)
    }
   & Trajectories from the origin, 
   from $ {\left[ 0.1 , 0.1 , 0.1 , 0.1 , 0.1 , 0.1 \right]}^T$ 
  from $ {\left[ 1 , 1 , 1 , 1 , 1 , 1 \right]}^T$ 
  and from $ {\left[ 10 , 10 , 10, 10, 10, 10 \right]}^T$ 
   converge to the same symmetric unimodal limit cycle. 
   \\
    \centerline{%
\(  \displaystyle
        {\frac{s^3}
        { \begin{Bmatrix}
            s^6 + s^5 + 5 s^4 + 3 s^3 \\
            + 5 s^2 +  s +  1
        \end{Bmatrix}
        }
        }  
\)
    }
   & Trajectories from the origin, 
   from $ {\left[ 0.1 , 0.1 , 0.1 , 0.1 , 0.1 , 0.1 \right]}^T$ 
  and from $ {\left[ 1 , 1 , 1 , 1 , 1 , 1 \right]}^T$ 
   converge to a symmetric unimodal limit cycle. 
   But Trajectories from 
   $ {\left[ 10 , 10 , 10, 10, 10, 10 \right]}^T$ 
   and from
   $ {\left[ 50 , 50 , 50, 50, 50, 50 \right]}^T$ 
   converge to a different symmetric unimodal limit cycle. 
\begin{center}
  \begin{tikzpicture}
\begin{axis}[%
              yscale = 1.3,
              yticklabel pos = right,
              width=5cm,
              height=2.3cm,
              ticks=none,
           ]
\addplot[thick, cyanTeal,
grid=both, every major grid/.style={gray, opacity=0.5}, 
    ]
    file {figuresAndData/limitCycle_sixthOrderNumeratorS2_startingAtAll1s.csv};
\end{axis}
\end{tikzpicture}
\end{center}
\begin{center}
  \begin{tikzpicture}
\begin{axis}[%
              yscale = 1.3,
              yticklabel pos = right,
              width=5cm,
              height=2.3cm,
              ticks=none,
           ]
\addplot[thick, cyanTeal,
grid=both, every major grid/.style={gray, opacity=0.5}, 
    ]
    file {figuresAndData/limitCycle_sixthOrderNumeratorS2_startingAtAll10s.csv};
\end{axis}
\end{tikzpicture}
\end{center}
\\
    \centerline{%
\(  \displaystyle
        {\frac{s^4}
        { \begin{Bmatrix}
            s^6 + s^5 + 5 s^4 + 3 s^3 \\
            + 5 s^2 +  s +  1
        \end{Bmatrix}
        }
        }  
\)
    }
   & 
   Trajectory from the origin chatters there.
   Trajectories 
   from $ {\left[ 0.1 , 0.1 , 0.1 , 0.1 , 0.1 , 0.1 \right]}^T$ 
   goes to a symmetric, unimodal limit cycle 
   Trajectories from
   $ {\left[ 1 , 1 , 1 , 1 , 1 , 1 \right]}^T$ 
   and 
   $ {\left[ 10 , 10 , 10, 10, 10, 10 \right]}^T$ 
   converge to a symmetric limit cycle, which 
   has {\textit{six switches}} per period.
   And trajectories from
   $ {\left[ 30 , 30 , 30, 30, 30, 30 \right]}^T$ 
   and from $ {\left[ 50 , 50 , 50, 50, 50, 50 \right]}^T$ 
   converge to yet another limit cycle, which is symmetric and
   unimodal.
\begin{center}
  \begin{tikzpicture}
\begin{axis}[%
              yscale = 1.3,
              yticklabel pos = right,
              width=5cm,
              height=2.3cm,
              ticks=none,
           ]
\addplot[thick, cyanTeal,
grid=both, every major grid/.style={gray, opacity=0.5}, 
    ]
    file {figuresAndData/limitCycle_sixthOrderNumeratorS4_startingAtAll1s.csv};
\end{axis}
\end{tikzpicture}
\end{center}
\begin{center}
  \begin{tikzpicture}
\begin{axis}[%
              yscale = 1.3,
              yticklabel pos = right,
              width=5cm,
              height=2.3cm,
              ticks=none,
           ]
\addplot[thick, cyanTeal,
grid=both, every major grid/.style={gray, opacity=0.5}, 
    ]
    file {figuresAndData/limitCycle_sixthOrderNumeratorS4_startingAtAll0pt1s.csv};
\end{axis}
\end{tikzpicture}
\end{center}
\\
    \centerline{%
\(  \displaystyle
        {\frac{- s^5 + s^4 - s^3 - s^2 -s}
        { \begin{Bmatrix}
            s^6 + s^5 + 5 s^4 + 3 s^3 \\
            + 5 s^2 +  s +  1
        \end{Bmatrix}
        }
        }  
\)
    }
    \vspace*{4pt}
    Plant has an odd number of positive real zeroes at +1, with multiplicity 2.
   & Trajectories from the origin, from $ {\left[ 0.1 , 0.1 , 0.1 \right]}^T$  and from $ {\left[ 1 , 1 , 1 \right]}^T$ converge to a symmetric unimodal limit cycle. The transient phase is very long, and quasi-periodic.
\begin{center}
  \begin{tikzpicture}
\begin{axis}[%
              yscale = 1.3,
              yticklabel pos = right,
              width=5cm,
              height=2.3cm,
              ticks=none,
           ]
\addplot[thick, cyanTeal,
grid=both, every major grid/.style={gray, opacity=0.5}, 
    ]
    file {figuresAndData/relayOutputSixthOrder_longTransient.csv};
\end{axis}
\end{tikzpicture}
\end{center}
\\
    \vspace*{1pt}
    \centerline{%
\(  \displaystyle
        {\frac{s^4 + s^2}
        { \begin{Bmatrix}
            s^6 + s^5 + 5 s^4 + 3 s^3 \\
            + 5 s^2 +  s +  1
        \end{Bmatrix}
        }
        }  
\)
    }
   & 
   Trajectory
   from $ {\left[ 1 , 1 , 1 \right]}^T$ 
   converges to the origin as a series of decaying spikes
   The transient phase is very long.
\begin{center}
    \includegraphics[width=1.1\linewidth]{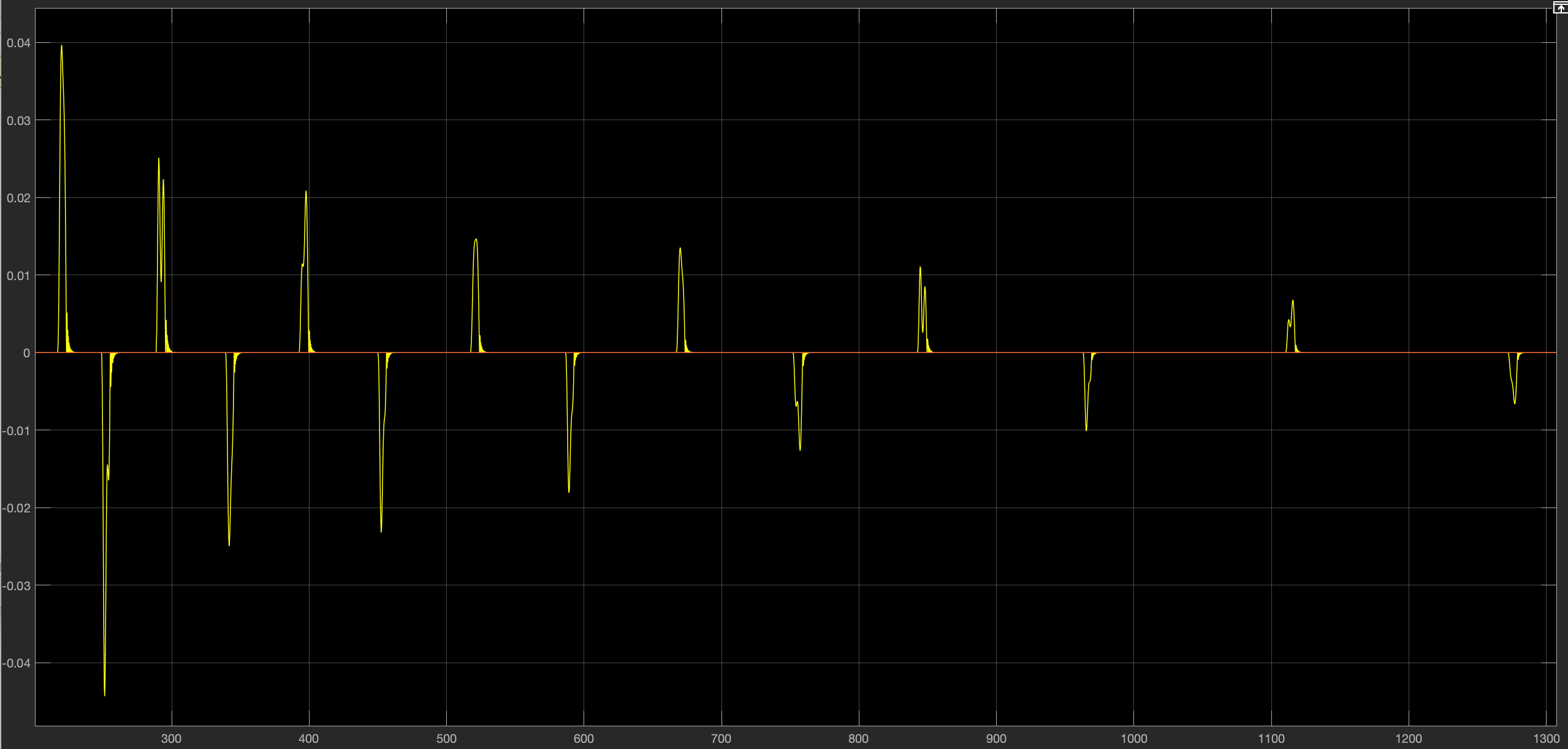}
\end{center}
\\
    \vspace*{1pt}
    \centerline{%
\(  \displaystyle
        {\frac{s^4 + s^3 + s^2}
        { \begin{Bmatrix}
            s^6 + s^5 + 5 s^4 + 3 s^3 \\
            + 5 s^2 +  s +  1
        \end{Bmatrix}
        }
        }  
\)
    }
   & 
   Trajectory
   from $ {\left[ 1 , 1 , 1 \right]}^T$ 
   converges to the origin as a series of decaying spikes
   The transient phase is very long.
\begin{center}
    \includegraphics[width=1.1\linewidth]{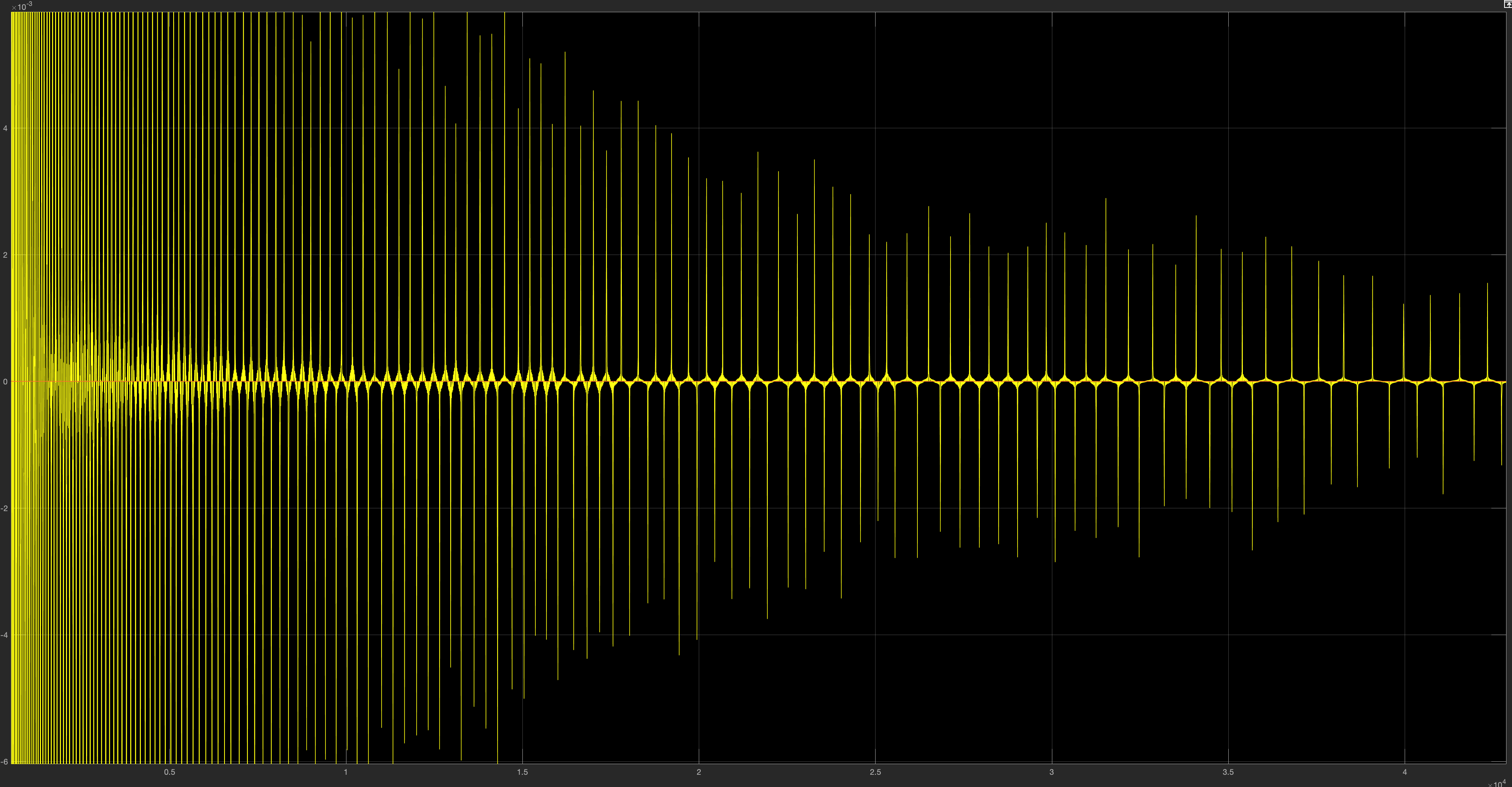}
\end{center}
\\
    \centerline{%
\(  \displaystyle
        {\frac{ s^4 - s^3 + s^2}
        { \begin{Bmatrix}
            s^6 + s^5 + 5 s^4 + 3 s^3 \\
            + 5 s^2 +  s +  1
        \end{Bmatrix}
        }
        }  
\)
    }
   & 
   Trajectory starting at the origin chatters there.
   Tarjectories from $ {\left[ 0.1 , 0.1 , 0.1 \right]}^T$  and from $ {\left[ 1 , 1 , 1 \right]}^T$ converge to a symmetric unimodal limit cycle. 
\begin{center}
  \begin{tikzpicture}
\begin{axis}[%
              yscale = 1.3,
              yticklabel pos = right,
              width=5cm,
              height=2.3cm,
              ticks=none,
           ]
\addplot[thick, cyanTeal,
         grid=both, every major grid/.style={gray, opacity=0.5}, 
    ]
    file {figuresAndData/relayOscillationSixthOrderNumS2minusS3plusS4_startingAtAll0pt1.csv};
\end{axis}
\end{tikzpicture}
\end{center}
\\
    \centerline{%
\(  \displaystyle
        {\frac{- s^4 -  s^2}
        { \begin{Bmatrix}
            s^6 + s^5 + 5 s^4 + 3 s^3 \\
            + 5 s^2 +  s +  1
        \end{Bmatrix}
        }
        }  
\)
    }
   & 
   Trajectories from the origin, and
   from $ {\left[ 1 , 1 , 1 , 1 , 1 , 1 \right]}^T$ 
   converge to a symmetric unimodal limit cycle. 
   Trajectories from 
   from $ {\left[ 10 , 10 , 10, 10, 10, 10 \right]}^T$ 
   converge to a different symmetric unimodal limit cycle. 
\begin{center}
  \begin{tikzpicture}
\begin{axis}[%
              yscale = 1.3,
              yticklabel pos = right,
              width=5cm,
              height=2.3cm,
              ticks=none,
           ]
\addplot[ thick, 
    cyanTeal,
grid=both, every major grid/.style={gray, opacity=0.5}, 
    ]
    file {figuresAndData/second_limitCycle_sixthOrderRelay_WithNumeratorMinusSpower4plusMinusSpower2_startingAtAllOnes.csv};
\end{axis}
\end{tikzpicture}
\end{center}
\begin{center}
  \begin{tikzpicture}
\begin{axis}[%
              yscale = 1.3,
              yticklabel pos = right,
              width=5cm,
              height=2.3cm,
              ticks=none,
           ]
\addplot[thick, 
    cyanTeal,
    grid=both, every major grid/.style={gray, opacity=0.5}, 
    ]
    file {figuresAndData/third_limitCycle_sixthOrderRelay.csv};
\end{axis}
\end{tikzpicture}
\end{center}
\\
    \centerline{%
\(  \displaystyle
        {\frac{s^5 + s^4 + s^3  +  s^2 +  s +  1}
        { \begin{Bmatrix}
            s^6 + s^5 + 5 s^4 + 3 s^3 \\
            + 5 s^2 +  s +  1
        \end{Bmatrix}
        }
        }  
\)
    }
    \vspace*{4pt}
   This plant has a complex conjugate pair of zeroes on the RHP. These are at: $
    0.5 \pm 0.866  \sqrt{-1}  $
   & 
   Trajectories from 
   reach a chattering set around the origin.
\\
    \vspace*{4pt}
    \centerline{%
\(  \displaystyle
        {\frac{ 
        \begin{Bmatrix}
            32 s^5 + 16 s^4 + 8 s^3 \\
            + 4 s^2 + 2 s + 1
        \end{Bmatrix}
        }
        { 
        \begin{Bmatrix}
            s^6 + s^5 + 5 s^4 + 3 s^3 \\
            + 5 s^2 +  s +  1
        \end{Bmatrix}
        }
        }  
\)
    }
    \vspace*{1pt}
    This plant has zeroes that are the zeroes of the previous plant scaled by two.
   Thus, this plant has a complex conjugate pair of zeroes on  at: $
    1 \pm 1.732  \sqrt{-1}  $
   & 
   Trajectory starting from the origin chatters there.
   Trajectories starting from $ {\left[ 0.1 , 0.1 , 0.1, 0.1, 0.1, 0.1 \right]}^T$,
   from $ {\left[ 1 , 1 , 1 , 1 , 1 , 1 \right]}^T$, and from 
   from $ {\left[ 10 , 10 , 10, 10, 10, 10 \right]}^T$ 
  all
   converge to the same symmetric unimodal limit cycle. 
\begin{center}
  \begin{tikzpicture}
\begin{axis}[%
              yscale = 1.3,
              yticklabel pos = right,
              width=5cm,
              height=2.3cm,
              ticks=none,
           ]
\addplot[ cyanTeal, thick,
         grid=both, every major grid/.style={gray, opacity=0.5}, 
    ]
    file {figuresAndData/limitCycle_sixthOrderRelay_RHPzerosBigEnough.csv};
\end{axis}
\end{tikzpicture}
\end{center}
\\
    \centerline{%
\(  \displaystyle
        {\frac{s^4 + s^3 + s^2 + s + 1}
        { \begin{Bmatrix}
            s^6 + s^5 + 5 s^4 + 3 s^3 \\
            + 5 s^2 +  s +  1
        \end{Bmatrix}
        }
        }  
\)
    }
    \vspace*{4pt}
    This plant has a complex conjugate pair of RHP zeroes at: 
    $ 0.3090 \pm 0.9511 \sqrt{-1} .$
   & 
   Trajectory starting from the origin chatters there.
   Trajectories starting from $ {\left[ 0.01 , 0.01 , 0.01, 0.01, 0.01, 0.01 \right]}^T$,
   from $ {\left[ 0.1 , 0.1 , 0.1, 0.1, 0.1, 0.1 \right]}^T$,
 and    from $ {\left[ 1 , 1 , 1 , 1 , 1 , 1 \right]}^T$ converge to a symmetric, unimodal limit cycle.
   Trajectory starting 
   from $ {\left[ 10 , 10 , 10, 10, 10, 10 \right]}^T,$ 
   from $ {\left[ 20 , 20 , 20, 20, 20, 20 \right]}^T,$ and 
   from $ {\left[ 30 , 30 , 30, 30, 30, 30 \right]}^T$ 
   all reach another symmetric, unimodal limit cycle.
\begin{center}
  \begin{tikzpicture}
\begin{axis}[%
              yscale = 1.3,
              yticklabel pos = right,
              width=5cm,
              height=2.3cm,
              ticks=none,
           ]
\addplot[ thick, cyanTeal,
         grid=both, every major grid/.style={gray, opacity=0.5}, 
    ]
    file {figuresAndData/limitCycle_6thOrder_num+1+1+1+1+1_startingAtAll0pt1.csv};
\end{axis}
\end{tikzpicture}
\end{center}
\begin{center}
  \begin{tikzpicture}
\begin{axis}[%
              yscale = 1.3,
              yticklabel pos = right,
              width=5cm,
              height=2.3cm,
              ticks=none,
           ]
\addplot[ thick, cyanTeal,
         grid=both, every major grid/.style={gray, opacity=0.5}, 
    ]
    file {figuresAndData/limitCycle_6thOrder_num+1+1+1+1+1_startingAtAllTens.csv};
\end{axis}
\end{tikzpicture}
\end{center}
\\
    \centerline{%
\(  \displaystyle
        {\frac{s^3 + s^2 + s + 1}
        { \begin{Bmatrix}
            s^6 + s^5 + 5 s^4 + 3 s^3 \\
            + 5 s^2 +  s +  1
        \end{Bmatrix}
        }
        }  
\)
    }
    \vspace*{4pt}
    This plant has no zeroes on the open RHP.
    It has a conjugate pair of zeroes on the imaginary axis:
    at $ \pm \sqrt{-1} $
   & 
   Trajectories starting from the origin ,
    from $ {\left[ 0.01 , 0.01 , 0.01, 0.01, 0.01, 0.01 \right]}^T$,
   from $ {\left[ 0.1 , 0.1 , 0.1, 0.1, 0.1, 0.1 \right]}^T$,
from $ {\left[ 1 , 1 , 1 , 1 , 1 , 1 \right]}^T$ 
 and
   from $ {\left[ 10 , 10 , 10, 10, 10, 10 \right]}^T$ 
 converge to a symmetric, unimodal limit cycle.
\begin{center}
  \begin{tikzpicture}
\begin{axis}[%
              yscale = 1.3,
              yticklabel pos = right,
              width=5cm,
              height=2.3cm,
              ticks=none,
           ]
\addplot[ thick, cyanTeal,
         grid=both, every major grid/.style={gray, opacity=0.5}, 
    ]
    file {figuresAndData/limitCycle_6thOrder_num+1+1+1+1_startingAtAllZeroes.csv};
\end{axis}
\end{tikzpicture}
\end{center}
\\
    \centerline{%
\(  \displaystyle
        {\frac{ s^4 + 2 s^3 + 3 s^2 + 3 s + 2}
        { \begin{Bmatrix}
            s^6 + s^5 + 5 s^4 + 3 s^3 \\
            + 5 s^2 +  s +  1
        \end{Bmatrix}
        }
        }  
\)
    }
    \vspace*{4pt}
    This plant has all of its zeroes on the open LHP.  
    It is minimum phase. 
   Its step response is quite oscillatory.
\begin{center}
  \begin{tikzpicture}
\begin{axis}[%
              yscale = 1.3,
              yticklabel pos = right,
              width=5cm,
              height=2.3cm,
              ticks=none,
           ]
            \begin{scope}
            \clip (axis cs:0,0) rectangle (axis cs:140,40);
        \addplot [ tolDarkYellow,] 
    file {figuresAndData/stepResponseWithInverseReaction_6thOrder_numerator+1+2+3+3+2.csv};
            \end{scope}
            \begin{scope}
            \clip (axis cs:0,0) rectangle (axis cs:140,-40);
        \addplot [ tolPaleRed, ultra thick] 
    file {figuresAndData/stepResponseWithInverseReaction_6thOrder_numerator+1+2+3+3+2.csv};
            \end{scope}
\end{axis}
\end{tikzpicture}
\end{center}
   & 
   Trajectory starting from the origin chatters there.
    Trajectories starting from $ {\left[ 0.01 , 0.01 , 0.01, 0.01, 0.01, 0.01 \right]}^T$,
   from $ {\left[ 0.1 , 0.1 , 0.1, 0.1, 0.1, 0.1 \right]}^T$,
from $ {\left[ 1 , 1 , 1 , 1 , 1 , 1 \right]}^T$ 
 and
   from $ {\left[ 10 , 10 , 10, 10, 10, 10 \right]}^T$ 
 converge to a symmetric, unimodal limit cycle.
\begin{center}
  \begin{tikzpicture}
\begin{axis}[%
              yscale = 1.3,
              yticklabel pos = right,
              width=5cm,
              height=2.3cm,
              ticks=none,
           ]
\addplot[ thick, cyanTeal,
         grid=both, every major grid/.style={gray, opacity=0.5}, 
    ]
    file {figuresAndData/limitCycle_6thOrder_num+1+2+3+3+2_startingAtAllTens.csv};
\end{axis}
\end{tikzpicture}
\end{center}
\\
    \centerline{%
\(  \displaystyle
        {\frac{ s^4 + 4 s^3 + 6 s^2 + 4 s + 1}
        { \begin{Bmatrix}
            s^6 + s^5 + 5 s^4 + 3 s^3 \\
            + 5 s^2 +  s +  1
        \end{Bmatrix}
        }
        }  
\)
    }
    \vspace*{4pt}
    This plant has four repeated zeroes at: $ - 1 .$   It is minimum phase. But yet, its step response shows an inverse response, like some minimum phase transfer functions discovered by Maciejowski~\cite{janMaciejowski2018rhpZeroesNotNecessaryForInverseResponse}.
\begin{center}
  \begin{tikzpicture}
\begin{axis}[%
              yscale = 1.3,
              yticklabel pos = right,
              width=5cm,
              height=2.3cm,
              ticks=none,
           ]
            \begin{scope}
            \clip (axis cs:0,0) rectangle (axis cs:140,40);
        \addplot [ tolDarkYellow,] 
    file {figuresAndData/stepResponseWithInverseReaction_6thOrder_numeratorFourTupleMinusOne.csv};
            \end{scope}
            \begin{scope}
            \clip (axis cs:0,0) rectangle (axis cs:140,-40);
        \addplot [ tolPaleRed, ultra thick] 
    file {figuresAndData/stepResponseWithInverseReaction_6thOrder_numeratorFourTupleMinusOne.csv};
            \end{scope}
\end{axis}
\end{tikzpicture}
\end{center}
   & 
   Trajectory starting from the origin chatters there.
    Trajectories starting from $ {\left[ 0.01 , 0.01 , 0.01, 0.01, 0.01, 0.01 \right]}^T$,
   from $ {\left[ 0.1 , 0.1 , 0.1, 0.1, 0.1, 0.1 \right]}^T$,
from $ {\left[ 1 , 1 , 1 , 1 , 1 , 1 \right]}^T$ 
 and
   from $ {\left[ 10 , 10 , 10, 10, 10, 10 \right]}^T$ 
 converge to a symmetric, unimodal limit cycle.
\begin{center}
  \begin{tikzpicture}
\begin{axis}[%
              yscale = 1.3,
              yticklabel pos = right,
              width=5cm,
              height=2.3cm,
              ticks=none,
           ]
\addplot[ thick, cyanTeal,
         grid=both, every major grid/.style={gray, opacity=0.5}, 
    ]
    file {figuresAndData/limitCycle_6thOrder_num+1+4+6+4+1_startingAtAllOnes.csv};
\end{axis}
\end{tikzpicture}
\end{center}
\\
    \centerline{%
\(  \displaystyle
        {\frac{ { \left( 5 s + 1 \right) }^4 }
        { \begin{Bmatrix}
            s^6 + s^5 + 5 s^4 + 3 s^3 \\
            + 5 s^2 +  s +  1
        \end{Bmatrix}
        }
        }  
\)
    }
    \vspace*{4pt}
    This plant has four repeated zeroes at: $ - 0.2 .$   It is minimum phase. 
    But yet, its step response shows an inverse response, like some minimum phase transfer functions discovered by Maciejowski~\cite{janMaciejowski2018rhpZeroesNotNecessaryForInverseResponse}. In fact the step response is `exceptionally bad.'
\begin{center}
  \begin{tikzpicture}
\begin{axis}[%
              yscale = 1.3,
              width=5cm,
              height=2.3cm,
           ]
            \begin{scope}
            \clip (axis cs:0,0) rectangle (axis cs:50,400);
        \addplot [ tolDarkYellow,] 
       file {figuresAndData/worstStepResponse_6ththOrder_numeratorFourtuple+5.csv};
            \end{scope}
            \begin{scope}
            \clip (axis cs:0,0) rectangle (axis cs:50,-400);
        \addplot [ tolPaleRed, ultra thick] 
        file {figuresAndData/worstStepResponse_6ththOrder_numeratorFourtuple+5.csv};
            \end{scope}
\end{axis}
\end{tikzpicture}
\end{center}
   & 
   Trajectory starting from the origin chatters there.
    Trajectories starting from $ {\left[ 0.01 , 0.01 , 0.01, 0.01, 0.01, 0.01 \right]}^T$,
   from $ {\left[ 0.1 , 0.1 , 0.1, 0.1, 0.1, 0.1 \right]}^T$,
from $ {\left[ 1 , 1 , 1 , 1 , 1 , 1 \right]}^T$ 
 and
   from $ {\left[ 10 , 10 , 10, 10, 10, 10 \right]}^T$ 
 reach the origin asymptotically. The transient has high frequency oscillations.
\begin{center}
  \begin{tikzpicture}
\begin{axis}[%
              yscale = 1.3,
              yticklabel pos = right,
              width=5cm,
              height=2.3cm,
              ticks=none,
           ]
\addplot[ very thin, tolRedMain, samples = 1000,
         grid=both, every major grid/.style={gray, opacity=0.5}, 
    ]
    file {figuresAndData/asymptoticDecay_6ththOrder_numeratorFourtuple+5_startingAtAllOnes.csv};
\end{axis}
\end{tikzpicture}
\end{center}
\\
    \centerline{%
\(  \displaystyle
        {\frac{- s^5}
        { \begin{Bmatrix}
            s^6 + s^5 + 5 s^4 + 3 s^3 \\
            + 5 s^2 +  s +  1
        \end{Bmatrix}
        }
        }  
\)
    }
   & 
   Trajectories starting from the origin, from $ {\left[ 0.1 , 0.1 , 0.1, 0.1, 0.1, 0.1 \right]}^T$,
   from $ {\left[ 1 , 1 , 1 , 1 , 1 , 1 \right]}^T$, and from 
   from $ {\left[ 10 , 10 , 10, 10, 10, 10 \right]}^T$ 
   all reach the same symmetric, unimodal limit cycle.
\begin{center}
  \begin{tikzpicture}
\begin{axis}[%
              yscale = 1.3,
              yticklabel pos = right,
              width=5cm,
              height=2.3cm,
              ticks=none,
           ]
\addplot[ thick, cyanTeal,
         grid=both, every major grid/.style={gray, opacity=0.5}, 
    ]
    file {figuresAndData/relayOutputSixthOrderNumeratorMinusS5_startingAtOrigin.csv};
\end{axis}
\end{tikzpicture}
\end{center}
\\
\midrule
    \multicolumn{2}{c}{\textcolor{cyanTeal}{\textit{Simulations of 12th order plants:}}} 
\\[6pt]
    \vspace*{1pt}
    \centerline{%
\(  \displaystyle
        {\frac{-  s + 1}
        { 
        \begin{Bmatrix}
            s^{12} + 2 s^{11} + 10 s^{10} \\
            + 15 s^9 + 35 s^8 + 37 s^7 \\
            + 53 s^6 + 37 s^5 + 35 s^4 \\
            + 15 s^3 + 10 s^2 + 2  s +  1
        \end{Bmatrix}
        }
        }  
\)
    }
   & 
   Trajectory from the origin 
   converges to a symmetric unimodal limit cycle. 
\begin{center}
    \includegraphics[width=1.1\linewidth]{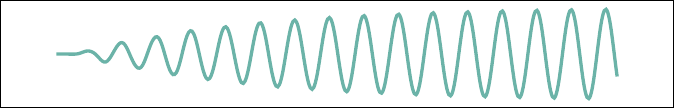}
\end{center}
\\
    \vspace*{1pt}
    \centerline{%
\(  \displaystyle
        {\frac{- s^3 + s^2 +  s + 1}
        { 
        \begin{Bmatrix}
            s^{12} + 2 s^{11} + 10 s^{10} \\
            + 15 s^9 + 35 s^8 + 37 s^7 \\
            + 53 s^6 + 37 s^5 + 35 s^4 \\
            + 15 s^3 + 10 s^2 + 2  s +  1
        \end{Bmatrix}
        }
        }  
\)
    }
   & 
   Trajectory from the origin, 
   converges to a symmetric unimodal limit cycle. 
\begin{center}
    \includegraphics[width=1.1\linewidth]{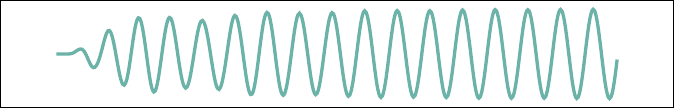}
\end{center}
\\
    \vspace*{3pt}
    \centerline{%
\(  \displaystyle
        {\frac{
        \begin{Bmatrix}
            - s^{7} + s^{6}  \\
            + s^{5} + s^{4} + s^{3}  \\
               + s^2 +  s + 1
        \end{Bmatrix}
        }
        { 
        \begin{Bmatrix}
            s^{12} + 2 s^{11} + 10 s^{10} \\
            + 15 s^9 + 35 s^8 + 37 s^7 \\
            + 53 s^6 + 37 s^5 + 35 s^4 \\
            + 15 s^3 + 10 s^2 + 2  s +  1
        \end{Bmatrix}
        }
        }  
\)
    }
   & 
   Trajectory from the origin, 
   converges to a symmetric unimodal limit cycle. 
\begin{center}
    \includegraphics[width=1.1\linewidth]{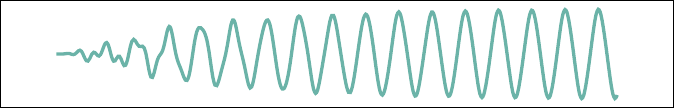}
\end{center}
\\
    \vspace*{10pt}
    \centerline{%
\(  \displaystyle
        {\frac{
        \begin{Bmatrix}
            - s^{9} \\
            + s^{8} + s^{7} + s^{6}  \\
            + s^{5} + s^{4} + s^{3}  \\
               + s^2 +  s + 1
        \end{Bmatrix}
        }
        { 
        \begin{Bmatrix}
            s^{12} + 2 s^{11} + 10 s^{10} \\
            + 15 s^9 + 35 s^8 + 37 s^7 \\
            + 53 s^6 + 37 s^5 + 35 s^4 \\
            + 15 s^3 + 10 s^2 + 2  s +  1
        \end{Bmatrix}
        }
        }  
\)
    }
   & 
   Trajectory from the origin, 
   converges to a symmetric unimodal limit cycle. 
\begin{center}
    \includegraphics[width=1.1\linewidth]{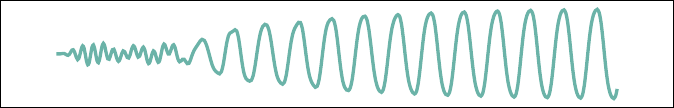}
\end{center}
\\
    \vspace*{10pt}
    \centerline{%
\(  \displaystyle
        {\frac{
        \begin{Bmatrix}
            -  s^{10} + s^{9} \\
            + s^{8} + s^{7} + s^{6}  \\
            + s^{5} + s^{4} + s^{3}  \\
               + s^2 +  s + 1
        \end{Bmatrix}
        }
        { 
        \begin{Bmatrix}
            s^{12} + 2 s^{11} + 10 s^{10} \\
            + 15 s^9 + 35 s^8 + 37 s^7 \\
            + 53 s^6 + 37 s^5 + 35 s^4 \\
            + 15 s^3 + 10 s^2 + 2  s +  1
        \end{Bmatrix}
        }
        }  
\)
    }
   & 
   Trajectory from the origin, 
   converges to a symmetric unimodal limit cycle. 
\begin{center}
    \includegraphics[width=1.1\linewidth]{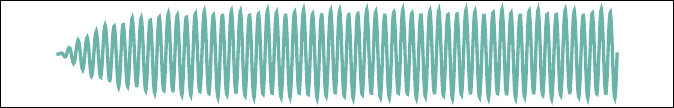}
\end{center}
\\
    \vspace*{10pt}
    \centerline{%
\(  \displaystyle
        {\frac{
        \begin{Bmatrix}
            - s^{11} + s^{10} + s^{9} \\
            + s^{8} + s^{7} + s^{6}  \\
            + s^{5} + s^{4} + s^{3}  \\
               + s^2 +  s + 1
        \end{Bmatrix}
        }
        { 
        \begin{Bmatrix}
            s^{12} + 2 s^{11} + 10 s^{10} \\
            + 15 s^9 + 35 s^8 + 37 s^7 \\
            + 53 s^6 + 37 s^5 + 35 s^4 \\
            + 15 s^3 + 10 s^2 + 2  s +  1
        \end{Bmatrix}
        }
        }  
\)
    }
    \vspace*{4pt}
    This plant is actually BRL-URF. It is listed here for comparison with the previous plants that have the same denominator, but greater pole excess.
   & 
   Trajectory from the origin, 
   converges to a symmetric unimodal limit cycle. 
\begin{center}
    \includegraphics[width=1.1\linewidth]{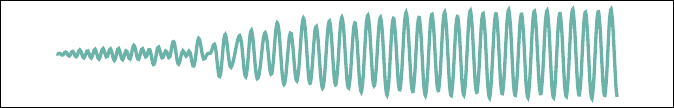}
\end{center}
\\
    \vspace*{10pt}
    \centerline{%
\(  \displaystyle
        {\frac{
            {\left( s + 1 \right)}^{10}
        }
        { 
        \begin{Bmatrix}
            s^{12} + 2 s^{11} + 10 s^{10} \\
            + 15 s^9 + 35 s^8 + 37 s^7 \\
            + 53 s^6 + 37 s^5 + 35 s^4 \\
            + 15 s^3 + 10 s^2 + 2  s +  1
        \end{Bmatrix}
        }
        }  
\)
    }
    \vspace*{4pt}
    This plant has ten repeated zeroes at: $ - 1 .$
    It is minimum phase. But yet, its step response shows an inverse response, like some minimum phase transfer functions discovered by Maciejowski~\cite{janMaciejowski2018rhpZeroesNotNecessaryForInverseResponse}.
\begin{center}
  \begin{tikzpicture}
\begin{axis}[%
              yscale = 1.3,
              width=5cm,
              height=2.3cm,
           ]
            \begin{scope}
            \clip (axis cs:0,0) rectangle (axis cs:140,40);
        \addplot [ tolDarkYellow,] 
    file {figuresAndData/stepResponseWithInverseReaction_12thOrder_numeratorTentupleMinusOne.csv};
            \end{scope}
            \begin{scope}
            \clip (axis cs:0,0) rectangle (axis cs:140,-40);
        \addplot [ tolPaleRed, ultra thick] 
    file {figuresAndData/stepResponseWithInverseReaction_12thOrder_numeratorTentupleMinusOne.csv};
            \end{scope}
\end{axis}
\end{tikzpicture}
\end{center}
   & 
   Trajectory starting from the origin chatters there.
   Trajectories starting from $ {\left[ 0.1 , 0.1 , 0.1, 0.1, 0.1, 0.1 \right]}^T$,
   from $ {\left[ 1 , 1 , 1 , 1 , 1 , 1 \right]}^T$, and from 
   from $ {\left[ 10 , 10 , 10, 10, 10, 10 \right]}^T$ 
   all reach the origin.
\begin{center}
  \begin{tikzpicture}
\begin{axis}[%
              yscale = 1.3,
              yticklabel pos = right,
              width=5cm,
              height=2.3cm,
              ticks=none,
           ]
\addplot[ thick, tolRedMain,
         grid=both, every major grid/.style={gray, opacity=0.5}, 
    ]
    file {figuresAndData/asymptoticStability_12thOrder_numeratorTentupleMinusOne_startingAtAllOnes.csv};
\end{axis}
\end{tikzpicture}
\end{center}
\\
 \end{xtabular}
 \bibliographystyle{IEEEtran}
 \bibliography{mabenRefs}

\begin{thebibliography}{10}
\providecommand{\url}[1]{#1}
\csname url@samestyle\endcsname
\providecommand{\newblock}{\relax}
\providecommand{\bibinfo}[2]{#2}
\providecommand{\BIBentrySTDinterwordspacing}{\spaceskip=0pt\relax}
\providecommand{\BIBentryALTinterwordstretchfactor}{4}
\providecommand{\BIBentryALTinterwordspacing}{\spaceskip=\fontdimen2\font plus
\BIBentryALTinterwordstretchfactor\fontdimen3\font minus
  \fontdimen4\font\relax}
\providecommand{\BIBforeignlanguage}[2]{{%
\expandafter\ifx\csname l@#1\endcsname\relax
\typeout{** WARNING: IEEEtran.bst: No hyphenation pattern has been}%
\typeout{** loaded for the language `#1'. Using the pattern for}%
\typeout{** the default language instead.}%
\else
\language=\csname l@#1\endcsname
\fi
#2}}
\providecommand{\BIBdecl}{\relax}
\BIBdecl

\bibitem{astromHagglund1984relayAutotuning}
\BIBentryALTinterwordspacing
K.~J. {\AA}str{\"o}m and T.~H{\"a}gglund, ``{A}utomatic tuning of simple
  regulators with specifications on phase and amplitude margins,''
  \emph{{A}utomatica}, vol.~20, no.~5, pp. 645--651, 1984. [Online]. Available:
  \url{https://www.sciencedirect.com/science/article/abs/pii/0005109884900141}
\BIBentrySTDinterwordspacing

\bibitem{boiko2021asymmetricOscillationsRelayFeedback}
\BIBentryALTinterwordspacing
I.~Boiko, N.~Kuznetsov, R.~Mokaev, and E.~Akimova, ``On asymmetric periodic
  solutions in relay feedback systems,'' \emph{Journal of the Franklin
  Institute}, vol. 358, no.~1, pp. 363--383, 2021. [Online]. Available:
  \url{https://www.sciencedirect.com/science/article/pii/S0016003220307158}
\BIBentrySTDinterwordspacing

\bibitem{flugge-lotz1953discontinuousControlSystems}
I.~Fl{\"u}gge-Lotz, \emph{Discontinuous control systems}.\hskip 1em plus 0.5em
  minus 0.4em\relax Princeton University press, 1952.

\bibitem{anosov1959relayChaos}
D.~V. Anosov, ``{O}n stability of equilibrium points of relay systems,''
  \emph{{A}utomation and {R}emote control}, vol.~20, no.~2, pp. 135--149, 1959.

\bibitem{cook1985relayChaos}
\BIBentryALTinterwordspacing
P.~Cook, ``Simple feedback systems with chaotic behaviour,'' \emph{Systems \&
  Control Letters}, vol.~6, no.~4, pp. 223 -- 227, 1985. [Online]. Available:
  \url{http://www.sciencedirect.com/science/article/pii/0167691185900714}
\BIBentrySTDinterwordspacing

\bibitem{holmberg1991thesisRelay}
U.~Holmberg, ``\BIBforeignlanguage{English}{Relay feedback of simple
  systems},'' Ph.D. dissertation, Department of Automatic Control, 1991.
[Online]. Available:
        \url{https://lucris.lub.lu.se/ws/portalfiles/portal/4722582/8568298.pdf}
\BIBentrySTDinterwordspacing

\bibitem{johanssonBarabanovAstrom2002chattering}
K.~H. Johansson, A.~E. Barabanov, and K.~J. {\AA}str{\"o}m, ``Limit cycles with
  chattering in relay feedback systems,'' \emph{IEEE Transactions on Automatic
  Control}, vol.~47, no.~9, pp. 1414--1423, Sep 2002.
[Online]. Available:
\url{https://lucris.lub.lu.se/ws/portalfiles/portal/4765935/625671.pdf}

\bibitem{sieber2010relayWithDelay}
\BIBentryALTinterwordspacing
J.~Sieber, P.~Kowalczyk, S.~Hogan, and M.~D. Bernardo, ``Dynamics of symmetric
  dynamical systems with delayed switching,'' \emph{Journal of Vibration and
  Control}, vol.~16, no. 7-8, pp. 1111--1140, 2010. [Online]. Available:
  \url{https://doi.org/10.1177/1077546309341124}
\BIBentrySTDinterwordspacing

\bibitem{jeffrey2018discontinuous}
\BIBentryALTinterwordspacing
M.~R. Jeffrey, \emph{Hidden dynamics}.\hskip 1em plus 0.5em minus 0.4em\relax
  Springer, Cham, 2018, the mathematics of switches, decisions and other
  discontinuous behaviour. [Online]. Available:
  \url{https://doi.org/10.1007/978-3-030-02107-8}
\BIBentrySTDinterwordspacing

\bibitem{astrom1995oscillationsRelay}
K.~J. {\AA}str{\"o}m, ``Oscillations in systems with relay feedback,'' in
  \emph{Adaptive Control, Filtering, and Signal Processing}, K.~J.
  {\AA}str{\"o}m, G.~C. Goodwin, and P.~R. Kumar, Eds.\hskip 1em plus 0.5em
  minus 0.4em\relax New York, NY: Springer New York, 1995, pp. 1--25.

\bibitem{chernysheva2014hopfBifurcationForRelay}
\BIBentryALTinterwordspacing
O.~A. Chernysheva, ``Andronov-{H}opf bifurcation theorem for relay systems,''
  \emph{J. Math. Sci. (N.Y.)}, vol. 200, no.~1, pp. 134--142, 2014, translated
  from Sovrem. Mat. Prilozh. No. 85 (2012). [Online]. Available:
  \url{https://doi.org/10.1007/s10958-014-1911-2}
\BIBentrySTDinterwordspacing

\bibitem{mabenRabi2021secondOrderRelayOscillations}
M.~Rabi, ``Relay self-oscillations for second order, stable, nonminimum phase
  plants,'' \emph{IEEE Transactions on Automatic Control}, vol.~66, no.~8, pp.
  1--1, 2021.

\bibitem{rabi2018relay}
\BIBentryALTinterwordspacing
------, ``Relay self-oscillations for second order, stable, nonminimum phase
  plants,'' 2018. [Online]. Available: \url{https://arxiv.org/abs/1810.11371}
\BIBentrySTDinterwordspacing

\bibitem{gillePelegrinDecaulne1959feedbackControlSystems}
J.-C. Gille, M.~J. Pelegrin, and P.~Decaulne, \emph{{F}eedback Control
  Systems}.\hskip 1em plus 0.5em minus 0.4em\relax MacGraw-Hill book company,
  Inc., 1959.

\bibitem{leMaitre1970graphicalAnalysisRelayOscillations}
\BIBentryALTinterwordspacing
J.-F. {Le Maitre}, J.-G. Paquet, and J.-C. Gille, ``A general approach for
  on-off control systems oscillations,'' \emph{Automatica}, vol.~6, no.~4, pp.
  609--613, 1970. [Online]. Available:
  \url{https://www.sciencedirect.com/science/article/pii/0005109870900154}
\BIBentrySTDinterwordspacing

\bibitem{tsypkin1984relayControlSystems}
Y.~Z. Tsypkin, \emph{Relay control systems}.\hskip 1em plus 0.5em minus
  0.4em\relax Cambridge: Cambridge University Press, 1984. 

\bibitem{bohn1961relay}
E.~{Bohn}, ``Stability margins and steady-state oscillations of on-off feedback
  systems,'' \emph{IRE Transactions on Circuit Theory}, vol.~8, no.~2, pp.
  127--130, June 1961.

\bibitem{bergen1962tsypkin}
A.~{Bergen}, ``A note on {T}sypkin's locus,'' \emph{IRE Transactions on
  Automatic Control}, vol.~7, no.~3, pp. 78--80, April 1962.

\bibitem{weischedel1973tsypkin}
H.~{Weischedel}, ``An exact method for the analysis of limit cycles in on-off
  control systems,'' \emph{IEEE Transactions on Automatic Control}, vol.~18,
  no.~1, pp. 40--44, February 1973.

\bibitem{juddChirlian1974graphicalAnalysisRelayLimitCycle}
\BIBentryALTinterwordspacing
F.~F. Judd and P.~M. Chirlian, ``Graphical analysis and design of limit cycles
  in autonomous relay control systems,'' \emph{International Journal of
  Control}, vol.~20, no.~2, pp. 321--334, 1974. [Online]. Available:
  \url{https://doi.org/10.1080/00207177408932740}
\BIBentrySTDinterwordspacing

\bibitem{judd1977errorBoundsRelaySystemAnalysis}
\BIBentryALTinterwordspacing
F.~F. Judd, ``Error bounds for approximate analysis of self-oscillations in
  autonomous relay control systes,'' \emph{International Journal of Control},
  vol.~25, no.~4, pp. 557--574, 1977. [Online]. Available:
  \url{https://doi.org/10.1080/00207177708922253}
\BIBentrySTDinterwordspacing

\bibitem{atherton1981relayOscillations}
\BIBentryALTinterwordspacing
D.~Atherton, ``Oscillations in relay systems,'' \emph{Transactions of the
  Institute of Measurement and Control}, vol.~3, no.~4, pp. 171--184, 1981.
  [Online]. Available: \url{https://doi.org/10.1177/014233128100300401}
\BIBentrySTDinterwordspacing

\bibitem{mees1981dynamicsOfFeedbackSystems}
A.~I. Mees, \emph{{D}ynamics of {F}eedback systems}, 1st~ed.\hskip 1em plus
  0.5em minus 0.4em\relax Chichester: A Wiley-Interscience publication, 1981.

\bibitem{kuznetsov2021tsypkinMethod}
\BIBentryALTinterwordspacing
N.~Kuznetsov, E.~Akimova, R.~Mokaev, and M.~Morozova, ``The study of periodic
  oscillations and global stability in the {T}al' model via the {T}sypkin
  method and the {LPRS} method,'' \emph{Journal of Physics: Conference Series},
  vol. 1864, no.~1, p. 012064, may 2021. [Online]. Available:
  \url{https://doi.org/10.1088/1742-6596/1864/1/012064}
\BIBentrySTDinterwordspacing

\bibitem{judd1975relationshipHamelTsypkin}
\BIBentryALTinterwordspacing
F.~F. Judd, ``Relationships between {T}sypkin, {H}amel and approximate limit
  cycle analyses†,'' \emph{International Journal of Control}, vol.~21, no.~4,
  pp. 641--653, 1975. [Online]. Available:
  \url{https://doi.org/10.1080/00207177508922018}
\BIBentrySTDinterwordspacing

\bibitem{andronovVittKhaikin1966theoryOfOscillations}
A.~A. Andronov, A.~A. Vitt, and S.~E. Khaikin, \emph{{T}heory of
  {O}scillations}, 2nd~ed., ser. International Series of Monographs in Physics,
  Vol. 4.\hskip 1em plus 0.5em minus 0.4em\relax Oxford: Pergamon press, 1966,
  {R}ussian edition published as {\textit{Teoriya kolebanii,}} Moscow:
  Gostekhizdat, in 1937.

\bibitem{minorsky1962nonlinearOscillations}
N.~Minorsky, \emph{{N}onlinear {O}scillations}, 1st~ed.\hskip 1em plus 0.5em
  minus 0.4em\relax Van Nostrand, 1962.

\bibitem{goncalvesMegretskiDahleh2001relay}
J.~M. Gon\c{c}alves, A.~Megretski, and M.~A. Dahleh, ``Global stability of
  relay feedback systems,'' \emph{IEEE Transactions on Automatic Control},
  vol.~46, no.~4, pp. 550--562, Apr 2001.

\bibitem{marioDiBernardoKarlJohansson2001ijbc}
\BIBentryALTinterwordspacing
M.~{di Bernardo}, K.~H. Johansson, and F.~Vasca, ``Self-oscillations and
  sliding in {R}elay feedback systems: symmetry and bifurcations,''
  \emph{International Journal of Bifurcation and Chaos}, vol.~11, no.~04, pp.
  1121--1140, 2001. [Online]. Available:
  \url{https://doi.org/10.1142/S0218127401002584}
\BIBentrySTDinterwordspacing

\bibitem{megretski1996globalStabilityOfRelayOscillations}
\BIBentryALTinterwordspacing
A.~Megretski, ``Global stability of oscillations induced by a relay feedback,''
  \emph{IFAC Proceedings Volumes}, vol.~29, no.~1, pp. 1931 -- 1936, 1996, 13th
  World Congress of IFAC, 1996, San Francisco USA, 30 June - 5 July. [Online].
  Available:
  \url{http://www.sciencedirect.com/science/article/pii/S1474667017579537}
\BIBentrySTDinterwordspacing

\bibitem{blimanKrasnoselskii1997periodicSolutions}
\BIBentryALTinterwordspacing
P.-A. Bliman and A.~Krasnosel'skii, ``Periodic solutions of linear systems
  coupled with relay,'' \emph{Nonlinear Analysis: Theory, Methods \&
  Applications}, vol.~30, no.~2, pp. 687 -- 696, 1997, proceedings of the
  Second World Congress of Nonlinear Analysts. [Online]. Available:
  \url{http://www.sciencedirect.com/science/article/pii/S0362546X96003720}
\BIBentrySTDinterwordspacing

\bibitem{varigondaGeorgiou2001relay}
S.~Varigonda and T.~T. Georgiou, ``Dynamics of relay relaxation oscillators,''
  \emph{IEEE Transactions on Automatic Control}, vol.~46, no.~1, pp. 65--77,
  Jan 2001.

\bibitem{johanssonRantzer1996globalAnalysisOfThirdOrderRelay}
\BIBentryALTinterwordspacing
K.~H. Johansson and A.~Rantzer, ``Global analysis of third-order relay feedback
  systems,'' \emph{IFAC Proceedings Volumes}, vol.~29, no.~1, pp. 1937 -- 1942,
  1996, 13th World Congress of IFAC, 1996, San Francisco USA, 30 June - 5 July.
  [Online]. Available:
  \url{http://www.sciencedirect.com/science/article/pii/S1474667017579549}
\BIBentrySTDinterwordspacing

\bibitem{allwright1977harmonicBalanceAndHopfBifurcation}
\BIBentryALTinterwordspacing
D.~J. Allwright, ``Harmonic balance and the {H}opf bifurcation,'' \emph{Math.
  Proc. Cambridge Philos. Soc.}, vol.~82, no.~3, pp. 453--467, 1977. [Online].
  Available: \url{https://doi.org/10.1017/S0305004100054128}
\BIBentrySTDinterwordspacing

\bibitem{stanSepulchre2007analysisOfInterconnectedOscillators}
G.-B. Stan and R.~Sepulchre, ``Analysis of interconnected oscillators by
  dissipativity theory,'' \emph{IEEE Transactions on Automatic Control},
  vol.~52, no.~2, pp. 256--270, 2007.

\bibitem{muratArcakSalama2020ringOScillatorsTanh}
X.~Ge, M.~Arcak, and K.~N. Salama, ``Nonlinear analysis of ring oscillator and
  cross-coupled oscillator circuits,'' \emph{Dyn. Contin. Discrete Impuls.
  Syst. Ser. B Appl. Algorithms}, vol.~17, no.~6, pp. 959--977, 2010.

\bibitem{alexanderYorke1978globalBifurcationsOfPeriodicOrbits}
\BIBentryALTinterwordspacing
J.~C. Alexander and J.~A. Yorke, ``Global bifurcations of periodic orbits,''
  \emph{Amer. J. Math.}, vol. 100, no.~2, pp. 263--292, 1978. [Online].
  Available: \url{https://doi.org/10.2307/2373851}
\BIBentrySTDinterwordspacing

\bibitem{balasubramanian1981stabilityOfLimitCycle}
\BIBentryALTinterwordspacing
R.~Balasubramanian, ``\BIBforeignlanguage{English}{Stability of limit cycles in
  feedback systems containing a relay},''
  \emph{\BIBforeignlanguage{English}{IEE Proceedings D (Control Theory and
  Applications)}}, vol. 128, pp. 24--29(5), January 1981. [Online]. Available:
  \url{https://digital-library.theiet.org/content/journals/10.1049/ip-d.1981.0005}
\BIBentrySTDinterwordspacing

\bibitem{majhiAtherton1998stabilityOfLimitCyclesInRelay}
\BIBentryALTinterwordspacing
S.~Majhi and D.~Atherton, ``Stability of limit cycles in relay control
  systems,'' \emph{IFAC Proceedings Volumes}, vol.~31, no.~17, pp. 91--96,
  1998, 4th IFAC Symposium on Nonlinear Control Systems Design 1998
  (NOLCOS'98), Enschede, The Netherlands, 1-3 July. [Online]. Available:
  \url{https://www.sciencedirect.com/science/article/pii/S1474667017403168}
\BIBentrySTDinterwordspacing

\bibitem{johanssonRantzerAstrom1999fastSwitchesInRFS}
\BIBentryALTinterwordspacing
K.~H. Johansson, A.~Rantzer, and K.~J. Åström, ``Fast switches in relay
  feedback systems,'' \emph{Automatica}, vol.~35, no.~4, pp. 539--552, 1999.
  [Online]. Available:
  \url{https://www.sciencedirect.com/science/article/pii/S0005109898001605}
\BIBentrySTDinterwordspacing

\bibitem{krantz2002implicitFunctionTheorem}
\BIBentryALTinterwordspacing
S.~G. Krantz and H.~R. Parks, \emph{The implicit function theorem}.\hskip 1em
  plus 0.5em minus 0.4em\relax Birkh\"{a}user Boston, Inc., Boston, MA, 2002,
  history, theory, and applications. [Online]. Available:
  \url{https://doi.org/10.1007/978-1-4612-0059-8}
\BIBentrySTDinterwordspacing

\bibitem{arnold1992odes}
V.~I. Arnol'd, \emph{Ordinary differential equations}, 3rd~ed.\hskip 1em plus
  0.5em minus 0.4em\relax Berlin: Springer-Verlag, 1992.

\bibitem{desoerVidyasagar2009}
\BIBentryALTinterwordspacing
C.~A. Desoer and M.~Vidyasagar, \emph{Feedback systems}, ser. Classics in
  Applied Mathematics.\hskip 1em plus 0.5em minus 0.4em\relax Society for
  Industrial and Applied Mathematics (SIAM), Philadelphia, PA, 2009, vol.~55,
  input-output properties, Reprint of the 1975 original [ MR0490289]. [Online].
  Available: \url{https://doi.org/10.1137/1.9780898719055.ch1}
\BIBentrySTDinterwordspacing

\bibitem{shihWu1998studiaMathematica}
\BIBentryALTinterwordspacing
M.-H. Shih and J.-W. Wu, ``Asymptotic stability in the {S}chauder fixed point
  theorem,'' \emph{Studia Mathematica}, vol. 131, no.~2, pp. 143--148, 1998.
  [Online]. Available: \url{http://matwbn.icm.edu.pl/ksiazki/sm/sm131/sm13124.pdf}
\BIBentrySTDinterwordspacing

\bibitem{vidyasagar1986undershoot}
M.~Vidyasagar, ``On undershoot and nonminimum phase zeros,'' \emph{IEEE
  Transactions on Automatic Control}, vol.~31, no.~5, pp. 440--440, May 1986.

\bibitem{vidyasagar1987authorsReply}
------, ``Author's reply,'' \emph{IEEE Transactions on Automatic Control},
  vol.~32, no.~3, pp. 272--272, March 1987.


\bibitem{irwin2001smoothDynamicalSystems}
\BIBentryALTinterwordspacing
M.~C. Irwin, \emph{Smooth dynamical systems}, 
\hskip 1em plus 0.5em minus 0.4em\relax World Scientific Publishing
  Co., Inc., River Edge, NJ, 2001, vol.~17, reprint of the 1980 original, With
  a foreword by R. S. MacKay. [Online]. Available:
  \url{https://doi.org/10.1142/9789812810120}
\BIBentrySTDinterwordspacing

\bibitem{coddingtonLevinson1955theoryOfODEs}
E.~A. Coddington and N.~Levinson, \emph{Theory of ordinary differential
  equations}.\hskip 1em plus 0.5em minus 0.4em\relax McGraw-Hill Book Co.,
  Inc., New York-Toronto-London, 1955.

\bibitem{jacquesWillems1968stabilityOfNonlinearOscillations}
J.~L. Willems, ``The stability of oscillations in nonlinear networks,''
  \emph{IEEE Transactions on Circuit Theory}, vol.~15, no.~3, pp. 284--286,
  1968.

\bibitem{khatskevichShoiykhet1994differentiableOperatorsAndNonlinearEquations}
\BIBentryALTinterwordspacing
V.~Khatskevich and D.~Shoiykhet, \emph{Differentiable operators and nonlinear
  equations}, ser. Operator Theory: Advances and Applications.\hskip 1em plus
  0.5em minus 0.4em\relax Birkh\"{a}user Verlag, Basel, 1994, vol.~66, 
  [Online]. Available:
  \url{https://doi.org/10.1007/978-3-0348-8512-6}
\BIBentrySTDinterwordspacing

\bibitem{wikipediaWeierstrassNonDifferentiableFunction}
\BIBentryALTinterwordspacing
``Weierstrass function,'' Apr 2021. [Online]. Available:
  \url{https://en.wikipedia.org/wiki/Weierstrass_function}
\BIBentrySTDinterwordspacing

\bibitem{krantz2013limitsOfHolomorphicSequences}
\BIBentryALTinterwordspacing
S.~G. Krantz, ``On limits of sequences of holomorphic functions,'' \emph{The
  Rocky Mountain Journal of Mathematics}, vol.~43, no.~1, pp. 273--283, 2013.
  [Online]. Available: \url{https://doi.org/10.1216/RMJ-2013-43-1-273}
\BIBentrySTDinterwordspacing

\bibitem{janMaciejowski2018rhpZeroesNotNecessaryForInverseResponse}
J.~M. Maciejowski, ``Right-half plane zeros are not necessary for inverse
  response,'' in \emph{2018 European Control Conference (ECC)}, 2018, pp.
  2488--2492.
  [Online]. Available: \url{https://doi.org/10.4064/sm-131-2-143-148}

\end{thebibliography}

\def\polhk#1{\setbox0=\hbox{#1}{\ooalign{\hidewidth
  \lower1.5ex\hbox{`}\hidewidth\crcr\unhbox0}}} \def\cprime{$'$}

\end{document}